\documentclass[11pt,reqno,a4paper]{amsart}

\usepackage{amsmath,textcomp,  microtype, lmodern, mathrsfs,csquotes,fbb}

\usepackage{hyperref}
\hypersetup{
  colorlinks=true,
  linkcolor=blue!50!green,
  citecolor=blue!50!green,
  urlcolor=blue!50!green,
  pdfauthor={Mayukh Mukherjee, Soumyadeb Samanta and Soumyadip Thandar},
  pdftitle={Harmonic Rigidity and Random Walk Asymptotics on Finitely Generated Groups}
}
\usepackage[nameinlink,capitalize,noabbrev]{cleveref}
\usepackage{marginnote}
\usepackage{url}
\usepackage[utf8]{inputenc} 
\usepackage[T1]{fontenc}    
\usepackage{url}            
\usepackage{booktabs}       
\usepackage{amsfonts}       
\usepackage{nicefrac}       
\usepackage{microtype}      
\usepackage[pdftex]{graphicx}
\usepackage{todonotes}
\usepackage{pdfpages}
\usepackage{amsmath}
\usepackage{tikz-cd}
\usepackage{amsfonts}
\usepackage{mathrsfs}  
\usepackage{centernot}
\usepackage{amssymb}
\usepackage{amsthm}
 
\usepackage[capitalize]{cleveref}
\usepackage[noadjust]{cite}
\usepackage{caption}
 
\usepackage{dutchcal}
\usepackage{tabularx}
\usepackage{bbm}
\usepackage{tikz}
\usetikzlibrary{decorations.pathreplacing}
 
\usepackage{mathtools}

\newtheorem*{theorem*}{Theorem}
\newtheorem{theorem}{Theorem}[section]
\newtheorem{remark}[theorem]{Remark}

\newtheorem{corollary}[theorem]{Corollary}

\newtheorem{question}[theorem]{Question}
\newtheorem{prop}[theorem]{Proposition}

\newtheorem{lemma}[theorem]{Lemma}

\newtheorem{definition}[theorem]{Definition}
\newtheorem{defn}[theorem]{Definition}

\theoremstyle{definition}

\theoremstyle{plain}  
\newcommand{\thistheoremname}{}
\newtheorem{generictheorem}[theorem]{\thistheoremname}

\usepackage{thmtools}

\declaretheoremstyle[
  headfont=\normalfont\bfseries,
  bodyfont=\itshape,
  headpunct={.},
  spaceabove=6pt, spacebelow=6pt,
  notebraces={(}{)},
  headformat=\NAME\NOTE   
]{introthm_nonum}

\declaretheorem[name=Theorem A, style=introthm_nonum, within=section]{theoremA}
\declaretheorem[name=Theorem B, style=introthm_nonum, within=section]{theoremB}
\declaretheorem[name=Theorem C, style=introthm_nonum, within=section]{theoremC}
\declaretheorem[name=Theorem D, style=introthm_nonum, within=section]{theoremD}
\declaretheorem[name=Theorem E, style=introthm_nonum, within=section]{theoremE}
\declaretheorem[name=Theorem F, style=introthm_nonum, within=section]{theoremF}
\declaretheorem[name=Theorem G, style=introthm_nonum, within=section]{theoremG}
\declaretheorem[name=Theorem H, style=introthm_nonum, within=section]{theoremH}



\crefname{theoremA}{Theorem~A}{Theorem~A}
\Crefname{theoremA}{Theorem~A}{Theorem~A}
\crefname{theoremB}{Theorem~B}{Theorem~B}
\Crefname{theoremB}{Theorem~B}{Theorem~B}
\crefname{theoremC}{Theorem~C}{Theorem~C}
\Crefname{theoremC}{Theorem~C}{Theorem~C}
\crefname{theoremD}{Theorem~D}{Theorem~D}
\Crefname{theoremD}{Theorem~D}{Theorem~D}
\crefname{theoremE}{Theorem~E}{Theorem~E}
\Crefname{theoremE}{Theorem~E}{Theorem~E}
\crefname{theoremF}{Theorem~F}{Theorem~F}
\Crefname{theoremF}{Theorem~F}{Theorem~F}
\crefname{theoremG}{Theorem~G}{Theorem~G}
\Crefname{theoremG}{Theorem~G}{Theorem~G}
\crefname{theoremH}{Theorem~H}{Theorem~H}
\Crefname{theoremH}{Theorem~H}{Theorem~H}
\Crefname{theoremI}{Theorem~I}{Theorem~I}
 
\makeatletter

\newcommand{\EE}{\mathbb E}
\newcommand{\E}{\mathbb E}
\newcommand{\PP}{\mathbb P}
\newcommand{\RR} {\mathbb R}
\newcommand{\CC} {\mathbb C}
\newcommand{\ZZ} {\mathbb Z}
\newcommand{\NN} {\mathbb N}

\newcommand{\HH} {\mathbb H}

\newcommand{\GG} {\mathbb G}
\newcommand{\Ind} {\mathbf 1}

\newcommand{\dist}{\operatorname{dist}}

\newcommand{\pa}{\partial}        

\newcommand{\Z}{\mathbb{Z}}

\DeclareMathOperator{\supp}{supp}
\DeclareMathOperator{\Vol}{Vol}

\newcommand{\Cal} {\mathcal}

\newcommand{\beq} {\begin{equation}}
\newcommand{\eeq} {\end{equation}}

\usepackage[top=2.5cm,bottom=2.4cm,left=2.8cm,right=2.8cm,headsep=0.2in]{geometry}

\usepackage{fancyhdr}

\begin{document}

\title[Strong Liouville, Green geometry and random walks]{Green geometry, Martin boundary and random walk asymptotics on groups}

\author{Mayukh Mukherjee}
\address{MM: Department of Mathematics, Indian Institute of Technology Bombay, Powai, Mumbai. 400 076. INDIA.}
\email{mathmukherjee@gmail.com}

\author{Soumyadeb Samanta}
\address{SS: Department of Mathematics, Indian Institute of Technology Bombay, Powai, Mumbai. 400 076. INDIA.}
\email{soumyadeb\texttt{@}iitb.ac.in}

\author{Soumyadip Thandar}
\address{ST: Theoretical Statistics and Mathematics unit, Indian Statistical Institute, Kolkata 700108,
INDIA.}
\email{stsoumyadip@gmail.com}


\begin{abstract}

We identify a single \emph{computationally checkable} analytic quantity interlacing Martin boundary collapse, Green geometry, and linear escape for transient random walks on finitely generated groups: the \emph{Green-variation functional}
\[
\Delta(S;a,b):=\max_{x\in\partial S}\frac{|G(a,x)-G(b,x)|}{G(a,x)}.
\]
We prove that $\Delta\to0$ along exhaustions characterises the \emph{strong Liouville property} (under mild, verifiable hypotheses on the ``strong Liouville $\Rightarrow \Delta\to0$'' direction), turning boundary oscillation estimates for Green kernels into potential-theoretic rigidity.

We then give two general criteria for $\Delta$-vanishing.  The first one derives quantitative bounds on $\Delta$ from coarse \emph{heat-kernel envelopes} at an intrinsic scale together with a Tauberian comparability, covering Gaussian/sub-Gaussian and stable-like regimes; and the second one is  purely \emph{elliptic}: an ``elliptic H\"older exhaustion'' criterion. Conversely, on groups of exponential growth, $\Delta$ fails to decay along balls already under stretched-exponential on-diagonal upper bounds, yielding a quantitative obstruction to strong Liouville.

As consequences, trivial Martin boundary forces linear-scale collapse of Green geometry ($d_G(e,x)=o(|x|)$) and vanishing Green speed (in probability), without any entropy hypothesis. On the non-Liouville side we prove an \emph{abundance} principle: the existence of a single minimal positive harmonic function at a prescribed growth scale forces infinitely many. Finally, we clarify the role of moment assumptions in speed theory: any linear-speed law of large numbers on a set of positive probability forces $\mathbb E|X_1|<\infty$, while on torsion-free nilpotent groups one can have $\mathbb E|X_1|=\infty$ yet $|X_n|/n\to0$ in probability.

\end{abstract}

\maketitle

\tableofcontents

\section{Introduction}

\subsection{Green analysis as a bridge between Martin boundary, geometry, and speed}

This paper investigates the interplay between \emph{Green geometry}, the \emph{Martin boundary}, and \emph{random-walk asymptotics} for transient random walks on finitely generated groups. Our guiding theme is the \emph{collapse (or non-collapse) of the Martin boundary}: on one hand, collapse corresponds to potential-theoretic rigidity (strong Liouville), and on the other hand, non-collapse comes with a rich supply of minimal harmonic functions and boundary structure.

We develop a single, computable analytic object that serves as a certificate for this collapse/non-collapse dichotomy and, at the same time, interfaces with large-scale geometry and speed. Concretely, we introduce the \emph{Green-variation functional}
\begin{equation}\label{def:Green_functional}
    \Delta(S;a,b):=\max_{x\in\partial S}\frac{|G(a,x)-G(b,x)|}{G(a,x)},
\end{equation}
evaluated on a finite set $S\subset\GG$ and two basepoints $a,b\in\GG$ (assuming $G(\cdot,x)>0$).
It measures the relative oscillation of the Green function $G(\cdot,x)$ on the boundary $\partial S$ when one shifts the starting point from $a$ to $b$.
Our first main observation is that, under very mild hypotheses, the decay of $\Delta$ along exhaustions $S\nearrow\GG$ is \emph{equivalent} to Liouville-type rigidity (strong Liouville), and hence provides an analytic certificate for Martin boundary collapse in the natural regimes of interest.

Our second contention is that $\Delta\to0$ is not merely a formal criterion: we give two general mechanisms that allow one to \emph{verify} $\Delta$-decay from coarse input.
One mechanism is parabolic, deriving quantitative $\Delta$-estimates from heat-kernel envelopes (as provided, for example, by a parabolic Harnack inequality);
the other is purely elliptic, based on an interior H\"older regularity along a suitable exhaustion and designed to remain stable under time changes (e.g., trap models).

Our third observation is that the certificate is somewhat sharp: on groups of exponential growth, balls cannot realise $\Delta\to0$ (indeed $\Delta$ stays bounded away from $0$ along a subsequence), even under very weak on-diagonal assumptions.
This provides a robust obstruction mechanism on the non-Liouville side.

Once the collapse/non-collapse picture is in place, we develop two families of consequences.
On the collapse side, trivial Martin boundary forces \emph{collapse of Green geometry at linear scale}: the Green metric $d_G$ becomes sublinear in the word metric and the Green speed vanishes, without any entropy hypothesis.
On the non-collapse side, once a non-trivial minimal harmonic function exists at a given growth scale, we prove an \emph{abundance} phenomenon: there are infinitely many distinct minimals at the same prescribed scale.

Finally, we revisit the role of moment assumptions in the study of speed and rigidity.
We show that any genuinely linear speed phenomenon (on a positive-probability event) forces integrability of the increments in complete generality, and we complement this with heavy-tailed constructions on nilpotent groups exhibiting infinite first moment but vanishing speed in probability.
We also include additional rigidity results for bounded harmonic functions on certain nilpotent groups under a dispersion hypothesis (see Proposition \ref{prop:Margulis_class2} and Definition \ref{def:dispersion_center}).

Throughout, $\GG$ is an infinite finitely generated group with a fixed finite symmetric generating set, and $|\cdot|$ denotes the associated word length.
A probability measure $\mu$ on $\GG$ is \emph{non-degenerate} if its support generates $\GG$.
We write $X_n=S_1\cdots S_n$ for the right random walk with i.i.d.\ increments $(S_i)$ having common law $\mu$.

\subsection{Main results}

We place the Green-variation functional $\Delta$ at the centre of our collapse/non-collapse analysis.

\medskip

\noindent(I) Collapse certificate: $\Delta\to0$ and strong Liouville.
Our first main theorem shows that $\Delta$-decay along exhaustions is equivalent to strong Liouville, with the direction ``strong Liouville $\Rightarrow \Delta\to0$'' requiring only mild hypotheses that hold in the principal cases of interest.

\begin{theoremA}[$\boldsymbol{\Delta\to0}$ characterises strong Liouville; cf. Proposition \ref{prop:Delta_from_strong_Liouville}, Theorem \ref{lem:Pol_2.1_conv}]
Let $\mu$ be non-degenerate and generate a transient walk on $\GG$.
\begin{itemize}
\item[(i)] (\emph{Strong Liouville $\Rightarrow \Delta\to0$.})
Assume $\mu$ is symmetric and there is a finite symmetric generating set $S\subset\supp\mu$ with $\rho:=\min_{s\in S}\mu(s)>0$ and
\begin{equation}\label{eq:rho_exp_mom}
\sum_{g\in\GG}\mu(g)\rho_0^{-|g|}<\infty
\quad\text{for some }\rho_0<\rho.
\end{equation}
If $(\GG,\mu)$ is strong Liouville, then for every $a,b\in\GG$ and every exhaustion $S_n\nearrow\GG$ one has
$$
\Delta(S_n;a,b)\longrightarrow 0.
$$
\item[(ii)] (\emph{$\Delta\to0$ $\Rightarrow$ strong Liouville.})
Conversely, if for all $a,b\in\GG$ and every exhaustion $S_n\nearrow\GG$ one has $\Delta(S_n;a,b)\to0$, then $(\GG,\mu)$ is strong Liouville.
\end{itemize}
\end{theoremA}

\emph{Comment.}
Observe that \eqref{eq:rho_exp_mom} subsumes two important (and natural) classes of measures: namely, symmetric finitely supported measures, and more generally, measures with superexponential moments. Part~(ii) does not assume symmetry or moment conditions: $\Delta\to0$ in our sense cannot hold unless every positive harmonic function is constant.

\medskip

\noindent (II) A practical calculus: two verification mechanisms for $\Delta$-decay.
The next two results explain how to \emph{check} $\Delta\to0$ from analytic or geometric input.

\begin{theoremB}[Heat-kernel envelopes $\Rightarrow$ $\boldsymbol{\Delta\to0}$; cf. Theorem \ref{thm:Delta-IH}]
Assume that the transition probabilities $p_n(x,y)$ satisfy:
\begin{itemize}
\item global on-diagonal upper and near-diagonal lower bounds at an intrinsic scale $\rho(n)$ with reference volume $\Vol(\rho(n))$;
\item interior H\"{o}lder regularity in the starting point at scale $\rho(n)$ (e.g.\ from a parabolic Harnack inequality);
\item a Tauberian comparability (TR$_\alpha$) between $\sum_{n<m}\Vol(\rho(n))^{-1}\rho(n)^{-\alpha}$ and $\rho(m)^{-\alpha}\sum_{n\ge m}\Vol(\rho(n))^{-1}$ for some $\alpha\in(0,1]$.
\end{itemize}
Then for every finite $S\subset\GG$ and $a,b\in\GG$,
$$
\Delta(S;a,b)\le C\left(\frac{d(a,b)}{R(S;a,b)}\right)^{\alpha},
$$
where $R(S;a,b)=\dist(\{a,b\},\partial S)$.
In particular, $\Delta(S_k;a,b)\to0$ along any exhaustion $S_k\nearrow\GG$.
\end{theoremB}

\emph{Comment.}
This ``envelope-to-$\Delta$'' scheme applies to symmetric random walks on virtually nilpotent groups with Gaussian or sub-Gaussian bounds, to many fractal graphs with sub-Gaussian heat kernels, and to stable-like/non-local processes on Ahlfors-regular spaces; see Remark \ref{rem:applications}.

\medskip

We also introduce a complementary criterion based on \emph{elliptic} regularity, which applies even when the time-parametrization of the walk is distorted.

\begin{theoremC}[Geometric $\Delta$-decay via Elliptic H\"older Exhaustions; cf. Theorem \ref{thm:EHE_to_Delta}]
Let $\mu$ be non-degenerate and transient. Let $(F_k)$ be an \emph{elliptic H\"older exhaustion} (Definition \ref{def:EHE}), meaning that positive $\mu$-harmonic functions on $F_k$ satisfy a uniform interior H\"older condition relative to the boundary distance $R_k=R(F_k;a,b)$.
Then for every $a,b\in\GG$,
$$
\Delta(F_k;a,b) \le C\left(\frac{d(a,b)}{R_k}\right)^{\alpha}\to 0 \text{ as } k\to\infty.
$$
In particular, $(\GG,\mu)$ is strong Liouville.
\end{theoremC}

\emph{Comment.}
Theorem \ref{thm:EHE_to_Delta} is decisive in regimes where global parabolic theory is unavailable but static potential theory remains robust.

\medskip

\noindent (III) Obstruction on balls: exponential growth forces non-decay of $\Delta$.
The next theorem shows that balls on exponentially growing groups \emph{cannot} realise $\Delta\to0$, even under very weak hypotheses.

\begin{theoremD}[Exponential growth $\Rightarrow$ no decay of $\boldsymbol{\Delta}$ on balls; cf. Theorem \ref{prop:exp_growth_delta_zero}]
Let $\GG$ be a finitely generated group of exponential growth with a fixed finite symmetric generating set, and let $\mu$ be symmetric and non-degenerate.
Assume that the on-diagonal heat kernel satisfies a stretched-exponential upper bound
$$
p_{2m}(e,e) \le C\exp(-cm^{\beta})
\qquad\text{for all large $m$, some }\beta>0.
$$
Then there exist a generator $t$ and $\delta_0>0$ such that
$$
\limsup_{n\to\infty}\ \Delta\big(B(e,n);e,t\big)\ \ge\ \delta_0.
$$
In particular, along the exhausting sequence $S_n=B(e,n)$, the quantities $\Delta(S_n;e,t)$ do not converge to $0$.
\end{theoremD}

\emph{Comment.}
On amenable groups of exponential growth, the stretched-exponential bound is automatic (Proposition \ref{prop:HK13_infinite_range}), with no moment assumptions.
Thus, on amenable exponential-growth groups, balls never achieve $\Delta\to0$ in the sense above, and strong Liouville fails; in this context, our result here extends previous work in  \cite{AK18, Po21}.

\medskip

\noindent\textbf{Consequences and applications.}
The remaining results develop two complementary directions: consequences of collapse for Green geometry and speed, and structure on the non-Liouville side.

\medskip

\noindent (IV) Collapse $\Rightarrow$ Green geometry collapse.
We next connect trivial Martin boundary (strong Liouville) to the Green metric and Green speed.

\begin{theoremE}[Speed-Liouville bridge via Green metric; cf. Theorem \ref{thm:form_entropy_harmonic}]
Let $\mu$ be symmetric, non-degenerate, and generate a transient random walk $(X_n)$ on $\GG$.
If the Martin boundary of $(\GG,\mu)$ is trivial, then
$$
\lim_{r\to\infty}\ \sup_{|x|=r}\ \frac{d_G(e,x)}{r}=0,
$$
and for every path $\omega$ in the Wiener space along which the word-metric speed $\lim_{n\to\infty}|X_n(\omega)|/n$ exists and is finite, the Green speed
$$
l_G(\omega):=\lim_{n\to\infty}\frac{d_G(e,X_n(\omega))}{n}
$$
exists and equals $0$. In particular, if the word speed is finite in probability, then $d_G(e,X_n)/n\to0$ in probability.
\end{theoremE}

\emph{Comment.}
Under trivial Martin boundary, the Green geometry collapses at linear scale purely from potential-theoretic rigidity, without any entropy or moment assumptions.

\medskip

\noindent (V) Non-collapse $\Rightarrow$ abundance at prescribed growth scales.
We now describe the structure of the \emph{non-Liouville} side, once non-trivial minimals exist at a given growth scale.
Recall that a function $\psi:\NN\to(0,\infty)$ is a \emph{growth gauge} if it is non-decreasing and translation-stable (Definition \ref{def:growth_gauge}).

\begin{theoremF}[Abundance at prescribed growth; cf.\ Theorem \ref{thm:pos_harm_alph_growtth_symm}]
Let $\mu$ be symmetric, non-degenerate, finitely supported on a group $\GG$ which is not virtually $\ZZ$ or $\ZZ^2$, and let $\psi$ be a growth gauge.
If there exists a non-constant minimal $\mu$-harmonic function whose growth is at least (respectively, at most) $\psi$, then there are \emph{infinitely many} distinct normalised minimal $\mu$-harmonic functions with growth at least (respectively, at most) $\psi$.
\end{theoremF}

\medskip

\noindent (VI) Random walk asymptotics beyond first moment.
We turn to the role of moment assumptions in the study of speed. Many classical linear-scale statements for random walks (laws of large numbers for $|X_n|/n$, drift formulas, and many boundary--speed comparisons) are formulated under a finite first moment hypothesis.
The next two results clarify what changes when $\mathbb E|X_1|=\infty$: we show that heavy tails do \emph{not} force ballistic behaviour (indeed one can still have vanishing speed in a probabilistic sense), while conversely any genuine linear-speed limit on a set of positive probability automatically forces integrability of the increments.

\begin{theoremG}[Infinite expectation, yet speed $\to0$ in probability; cf. Theorem \ref{thm:nilpotent-wlln}]\label{thmG}
If $\GG$ is a finitely generated torsion free nilpotent group, then for every symmetric measure $\mu$ satisfying \eqref{eq:shell-condition} we have
$$
\mathbb{E}|X_1|=\infty\qquad\text{and}\qquad\frac{|X_n|}{n}\to 0 \text{ in probability}.
$$
\end{theoremG}

\begin{theoremH}[Finite speed with positive probability forces finite moment; cf. Theorem \ref{thm:finite-speed-finite-moment}]\label{thmH}
Let $\GG$ be any group endowed with a symmetric length function, and let $(X_n)$ be i.i.d.\ $\GG$-valued increments with law~$\mu$.
If there exists an event $E$ with $\mathbb P(E)>0$ and $v<\infty$ such that
$$
\frac{|X_n|}{n}\longrightarrow v
\quad\text{on }E,
$$
then $\mathbb E|X_1|<\infty$.
\end{theoremH}

\emph{Comment.}
\Cref{thmH} shows that as soon as one has a genuine linear-speed law of large numbers (even with limit $v=0$) on a set of positive probability, the finite first moment condition is \emph{automatic}.
In particular, in the Liouville/collapse regime, Theorem \ref{thm:form_entropy_harmonic} yields vanishing Green speed along any trajectory with finite word speed, and Theorem \ref{thm:finite-speed-finite-moment} guarantees that whenever such a word-speed limit exists (on a positive-probability event) one is already in the finite-moment regime.
Thus the infinite-moment phenomena in Theorem \ref{thm:nilpotent-wlln} necessarily lie beyond the almost-sure linear-speed framework, and the corresponding speed collapse can only hold in weaker modes of convergence (here, in probability).

\medskip

\noindent\emph{Further rigidity via dispersion.}
We also introduce a dispersion property along the centre (Definition \ref{def:dispersion_center}) and show that, on nilpotent groups of class at most $2$, it forces bounded harmonic functions to factor through the abelianisation; see Proposition \ref{prop:Margulis_class2}.

\subsection{Organisation of the paper}

Section \ref{sec:preliminaries} collects preliminaries on random walks, harmonic functions, and the Green metric.

Section \ref{sec:del_char} is devoted to the $\Delta$-framework.
We prove the strong-Liouville $\leftrightarrow\Delta\to0$ equivalence (Theorem A) and establish some comparison results between $\Delta$ and the Amir-Kozma functional. 

Section \ref{sec:del_framework} establishes the pointwise criteria for $\Delta$-decay: the heat-kernel envelope method (Theorem B) and the Elliptic H\"older Exhaustion method (Theorem C).

Section \ref{sec:del_obstruction} treats exponential-growth groups, proving the failure of $\Delta\to0$ on balls under stretched-exponential assumptions (Theorem D).

Section \ref{sec:poly_growing_harmonic_functions} develops the abundance theory for positive harmonic functions (Theorem F), with some related construction of a Martin kernel.

Section \ref{sec:speed} develops the speed theory: we first establish the Green-speed collapse under trivial Martin boundary (Theorem E). We then prove the ``finite speed $\Rightarrow$ integrability'' principle (Theorem H) and construct heavy-tailed symmetric measures with infinite first moment and vanishing speed in probability on virtually nilpotent groups (Theorem G). We end this section by introducing the dispersion property and prove the abelian factorization rigidity for class-2 nilpotent groups.

\section{Preliminaries}\label{sec:preliminaries}

\subsection{Notations} For a set $D$ and functions $f,g:D \to \RR$ we say 
\begin{enumerate}
    \item $f \lesssim g$ if there exists $C>0$ such that $f(x) \le C g(x)$ for all $x\in D$,
    \item $f \gtrsim g$ if $g \lesssim f$, and
    \item $f\asymp g$ if there exists $C,D>0$ such that $Cf(x) \le g(x) \le Df(x)$ for all $x\in D$.
\end{enumerate} 

\subsection{Harmonic functions and SAS measures}

Throughout the paper, $\GG$ is a finitely generated group with identity element $e$, and $\mu$ is a probability measure (not assumed to be finitely supported unless mentioned explicitly) on $\GG$. We denote by $\mu^{(n)}$ the $n$-fold convolution power of $\mu$.

Fix once and for all a finite symmetric generating set $S$ of $\GG$, and let $|x|$ denote the corresponding word length of $x \in \GG$. (Different choices of $S$ give equivalent notions of exponential moments.)

Typically, we will need $\mu$ to satisfy certain additional conditions; otherwise the corresponding harmonic functions may not have nice properties. There are standard conditions of this kind in the existing literature, and we start by recalling the ones that are relevant to us.

\begin{defn}[Symmetric measure]
The probability measure $\mu$ is called \emph{symmetric} if $\mu(g)=\mu(g^{-1})$ for all $g \in \GG$.
\end{defn}

\begin{defn}[Non-degenerate/adapted measure]\label{defn:non_degenerate_measure}
The probability measure $\mu$ is called \emph{non-degenerate} (or \emph{adapted}) if $\supp(\mu)$ generates $\GG$ as a semigroup.
\end{defn}

\begin{defn}[Smooth measure]\label{def:smooth_measure}
A measure $\mu$ on a group $\GG$ is called \emph{smooth} if the generating function 
\beq\label{eq:exp_moment}
\Psi(\zeta) := \sum_{x \in \GG}\mu(x) e^{\zeta|x|} < \infty
\eeq  
for some positive real number $\zeta$. The measure $\mu$ is said to have \emph{superexponential moments} if \eqref{eq:exp_moment} holds for all real $\zeta > 0$.
\end{defn}

Clearly, the above definition provides control on all the moments of $X$, where $X$ is a random variable taking values in $\GG$ with law $\mu$. Since the definition looks somewhat unintuitive, we include a few words on the importance of smooth measures and why they are natural. Recall that a random variable $X$ is said to have \emph{exponential tail} if
$$
\PP\bigl( |X| > t\bigr) \le c_1 e^{-c_2 t}, \quad t \ge 0,
$$
for some positive constants $c_j$. One can calculate that this guarantees $\EE\bigl( e^{\alpha |X|}\bigr) < \infty$ for some $\alpha > 0$. Conversely, if $\EE\bigl( e^{\alpha |X|}\bigr) < \infty$ for some $\alpha > 0$, then $|X|$ has an exponential tail. In other words, a measure $\mu$ on $\GG$ is smooth if and only if the length of a $\mu$-random element of $\GG$ has an exponential tail.

\begin{defn}[Harmonic function]
Let $\GG$ be a group and $\mu$ be a probability measure on $\GG$. A function $f: \GG \to \CC$ is \emph{$\mu$-harmonic} at $k \in \GG$ if 
\beq\label{eq:def_harmonic}
f(k)=\sum_{g\in \GG}\mu(g) f(kg),
\eeq
and the above sum converges absolutely. We say that $f$ is $\mu$-harmonic on $\GG$ if \eqref{eq:def_harmonic} holds at all $k \in \GG$.
\end{defn}

By slight abuse of notation, we shall call a function $f$ simply \emph{harmonic} when the corresponding measure $\mu$ is tacitly understood in context, or not critical to the discussion.

\begin{remark}
(a) The group $\GG$ acts naturally on the set of harmonic functions by
$$
(g \cdot f)(k)=f(g^{-1}k), \quad g,k\in \GG.
$$
If $f$ is $\mu$-harmonic, then so is $g \cdot f$ for every $g\in\GG$.

\smallskip
(b) It can be proved by induction that if $f: \GG \to \CC$ is $\mu$-harmonic then it is $\mu^{(n)}$-harmonic for every $n\in\NN$, i.e.
$$
f(k)=\sum_{g\in\GG} \mu^{(n)}(g) f(kg), \quad k\in\GG.
$$
\end{remark}

\begin{defn}\label{def:gradient}
Let $S$ be the fixed symmetric generating set for $\GG$, and let $f : \GG \to \CC$ be a function. We define the (discrete) gradient of $f$ at $x\in\GG$ by
$$
    |\nabla f(x)| := \sup_{s \in S} | f(xs) - f(x) |, \quad x \in \GG.
$$
\end{defn}

\subsection{Martin boundary and the Poisson-Martin representation theorem}\label{subsec:Martin_bdry_representation}

Let $\GG$ be a finitely generated group and $\mu$ be a probability measure on $\GG$ such that the associated random walk on $\GG$ is irreducible and transient. Let $p_n$ denote the $n$-step transition probabilities of the $\mu$-random walk $(X_k)_{k\ge 0}$ on $\GG$, i.e. 
\begin{align*}
    p_n(x, y) & = \mu^{(n)}(x^{-1} y), \quad n \ge 1,\\
    p_0(x, y) &= \begin{cases}
        1, & x = y, \\
        0, & x \neq y,
    \end{cases}
\end{align*}
for all $x, y \in \GG$.

The \emph{Green function} $G(x, y)$ on $\GG \times \GG$ is defined by 
\beq\label{eqn:green_function}
G(x, y) := \sum_{k= 0}^\infty p_k(x, y).
\eeq
Since the random walk is transient, the sum in \eqref{eqn:green_function} is finite for all $x,y\in\GG$. The study of the asymptotics of the Green's function corresponding to a random walk is an active area of research in probability theory (for instance, see \cite{Dussaule2023,GekhtmanGerasimovPotyagailoYang2021,DenisovWachtel2024} and references therein). In this paper, we use the Green's function to construct a variant of the functional $\epsilon(S)$ defined in \cite[pp. 2]{AK18}; the idea of this functional seems implicit in several earlier works (e.g. see \cite[Chapter 2, Sections 2, 3]{ScYa94}). For a general background, also see \cite{Anderson1983, Sullivan1983, AndersonSchoen1985}.

Fix $y \in \GG$, and define
$$
K_y(x, z) := \frac{G(x, z)}{G(y, z)}, \qquad x, z \in \GG.
$$
By irreducibility and transience, $0<G(y,z)<\infty$ for all $z\in\GG$, so $K_y$ is well defined. The function $K_y$ is called the Martin kernel on $\GG \times \GG$ based at $y$. Fix a basepoint $o \in \GG$ (usually we take $o = e$, the identity of $\GG$), and set $K(x,z) := K_o(x,z)$.

Recall that the \emph{Martin compactification} of $\GG$ is the unique smallest compactification $\hat{\GG}(\mu)$ of $\GG$ to which the Martin kernels $K(x, \cdot)$ extend continuously (as functions of the second variable). The \emph{Martin boundary} is then defined as
$$
\Cal M(\GG) = \Cal M(\GG,\mu) := \hat{\GG}(\mu) \setminus \GG,
$$
which is well-defined (up to homeomorphism) under a change of the basepoint $o$; see \cite[Section 24]{Woess2000} for a detailed construction. For interested readers, we also suggest \cite{Martin1941, Dynkin1969}.

The \emph{minimal Martin boundary} $\Cal M_{\mathrm{min}}(\GG)$ is the set of all $\xi \in \Cal M(\GG)$ such that $K(\cdot,\xi)$ is a minimal positive harmonic function (in the first variable).

We now state one of the main theorems in this context that will be useful for us:

\begin{theorem}[Poisson-Martin representation theorem]\label{thm:Poisson-Martin-repn-thm}
Let $\GG$ be a finitely generated group and $\mu$ be a probability measure on $\GG$ such that the associated random walk on $\GG$ is irreducible and transient. For every positive $\mu$-harmonic function $h$ on $\GG$ there exists a positive Borel measure $\nu_h$ on $\mathcal{M}(\GG)$
such that
\beq\label{eq:Poisson-Martin-repn-thm}
h(x) = \int_{\mathcal{M}(\GG)} K(x,\xi) d\nu_h(\xi), \qquad x\in\GG.
\eeq
\end{theorem}

Moreover, the measure $\nu_h$ in Theorem \ref{thm:Poisson-Martin-repn-thm} can be uniquely chosen so that
$$
\nu_h\bigl(\Cal M(\GG) \setminus \Cal M_{\mathrm{min}}(\GG)\bigr)=0.
$$

\section{A $\Delta$-characterisation of strong Liouville}\label{sec:del_char}

Fix a symmetric generating set $T$ of $\GG$. Let $S \subseteq \GG$ and $\mu$ be a non-degenerate probability measure on $\GG$. Recall that $x\in \partial S$ if $x\notin S$ and there exists $y\in S$ such that $x=yt$ for some $t\in T$ (i.e.\ $x$ and $y$ are adjacent in the Cayley graph of $\GG$ corresponding to $T$).

The following discussion is intended to be a variant of the harmonic measure construction in \cite{AK18}, \cite[Chapter~2]{ScYa94}. Recall that Amir-Kozma define the functional
\beq\label{def:Amir_Kozma_functional_org}
\epsilon(S; a, b) := \max_{x \in \partial S} \frac{|\mu_S(a, x) - \mu_S(b, x)|}{|\mu_S(a, x)|},
\eeq
where $\mu_S(p, x)$ denotes the probability that the $\mu$-random walk starting at $p$ exits $S$ at $x \in \partial S$. Here, we work instead with the functional $\Delta(S;a,b)$ defined in \eqref{def:Green_functional}, because we believe it is analytically somewhat more amenable. This raises the following question.

\begin{question}\label{ques:eps_del_compare}
    Are the functionals $\epsilon(S)$ and $\Delta(S)$ comparable as $S \nearrow \GG$? 
\end{question}

We make some remarks on the above question in the latter part of this section. 

Recall that \cite[Proposition~2.1]{Po21} states a result of the form below; we remark that \cite{Po21} does not seem to assume transience of $\mu$ or superexponential moments (see also the discussion below \cite[Lemma~1]{AK18}).

\begin{lemma}[Limits of normalised Green kernels]\label{lem:limit-green}
Suppose $\mu$ is a non-degenerate probability measure on a finitely generated group $\GG$ which generates a transient random walk on $\GG$.
Fix $a\in \GG$ and let $(x_k)$ be a sequence in $\GG$ which leaves every finite subset of $\GG$ eventually (i.e.\ $x_k\to\infty$). Define
$$
\psi_k(v):=\frac{G(v,x_k)}{G(a,x_k)},\qquad v\in\GG,
$$
and suppose $\psi_k\to\psi$ pointwise along a subsequence. Then:

\emph{(a)} $\psi$ is $\mu$-superharmonic: for all $y\in\GG$,
$$
\psi(y) \ge \sum_{s\in\GG}\mu(s)\psi(ys).
$$

\emph{(b)} If, for each $y\in\GG$, there exists $H_y:\GG\to[0,\infty)$ with
$$
\psi_k(ys)\le H_y(s)\quad\forall k,s
\qquad\text{and}\qquad
\sum_{s\in\GG}\mu(s)H_y(s)<\infty,
$$
then $\psi$ is $\mu$-harmonic:
$$
\psi(y) = \sum_{s\in\GG}\mu(s)\psi(ys)\qquad\forall y\in\GG.
$$
\end{lemma}

\begin{proof}
For all $x,y\in \GG$ we have the following identity
\begin{equation}\label{eq:res}
G(y,x) = \sum_{s\in \GG}\mu(s)G(ys,x)+\delta_{y,x},
\end{equation}
Dividing \eqref{eq:res} by $G(a,x_k)$ yields
\begin{equation}\label{eq:res-norm}
\psi_k(y) = \sum_{s\in \GG}\mu(s)\psi_k(ys)+\frac{\delta_{y,x_k}}{G(a,x_k)}\qquad\text{for all }k\in\NN.
\end{equation}
Let $(\psi_{n_k})$ be the subsequence of $(\psi_k)$ which converges to $\psi$ pointwise on $\GG$.

\medskip\noindent
\emph{(a) Superharmonicity.}
Fix $y\in \GG$. Since $x_k\to\infty$, for all sufficiently large $k$ we have $x_k\neq y$, so the last term in \eqref{eq:res-norm} vanishes and
$$
\psi_k(y) = \sum_{s\in \GG}\mu(s)\psi_k(ys)\qquad\text{for all large $k$}.
$$
By Fatou's lemma,
$$
\sum_{s\in\GG}\mu(s)\psi(ys)
\le
\liminf_{k\to\infty}\sum_{s\in\GG}\mu(s)\psi_{n_k}(ys)
=
\liminf_{k\to\infty}\psi_{n_k}(y)
=
\psi(y),
$$
where the last equality uses the assumed pointwise convergence $\psi_{n_k}(y)\to\psi(y)$
along the subsequence. Thus $\psi$ is $\mu$-superharmonic at $y$, and since $y$ was
arbitrary, on all of $\GG$.

\medskip\noindent
\emph{(b) Harmonicity under domination.}
Assume now the domination hypothesis. For fixed $y\in\GG$, the family
$$
\bigl\{\psi_k(ys):k\in\NN\bigr\}_{s\in\GG}
$$
is dominated by $(H_y(s))_{s\in\GG}$, and $\sum_s\mu(s)H_y(s)<\infty$ by assumption. Hence,
by dominated convergence,
$$
\sum_{s\in\GG}\mu(s)\psi_{n_k}(ys) \longrightarrow \sum_{s\in\GG}\mu(s)\psi(ys)
\qquad\text{as }k\to\infty.
$$
On the other hand, for all large $k$ we have 
$\psi_k(y) = \sum_{s\in\GG}\mu(s)\psi_k(ys)$.
Passing to the limit along the subsequence $(\psi_{n_k})$ gives
$$
\psi(y) = \sum_{s\in\GG}\mu(s)\psi(ys).
$$
Thus $\psi$ is $\mu$-harmonic at $y$, and since $y$ was arbitrary, on all of $\GG$.
\end{proof}

\begin{corollary}[Easy domination criteria]\label{cor:easy-domination}
Assume that $\operatorname{supp}\mu$ contains a finite generating set $S$ of $\GG$, and set
$$
\rho:=\min_{s\in S}\mu(s)\in(0,1].
$$
If
$$
\sum_{g \in \GG}\mu(g)\rho^{-|g|}<\infty
\quad\text{(\emph{an exponential moment})},
$$
then for any sequence $(x_k)$ as in Lemma \ref{lem:limit-green} the domination condition in Lemma \ref{lem:limit-green}(b) holds with
$$
H_y(s) := \rho^{-|a^{-1}y|-|s|},
$$
and hence any pointwise subsequential limit $\psi$ of $\psi_k(v):=G(v,x_k)/G(a,x_k)$ is $\mu$-harmonic.
In particular, if $\mu$ is finitely supported then every such limit $\psi$ is harmonic.
\end{corollary}

\begin{proof}
We first prove a basic comparison estimate for the Green function. Fix $z\in\GG$. Since $S$ generates $\GG$, we can choose $t_1,t_2,\dots t_m \in S$ with $m = |a^{-1}z|$ such that $a t_1\cdots t_m = z$.
Let
$$
p(z) := \prod_{i=1}^m \mu(t_i) \ >0.
$$
Consider the $\mu$-random walk started at $a$. On the event that the first $m$ steps follow
the word $t_1,\dots,t_m$ (which has probability $p(z)$), the position at time $m$ is $z$, and the expected number of future visits to a point $x$ is $G(z,x)$. Therefore, by the Markov property,
$$
G(a,x)  \ge p(z)G(z,x)\qquad\text{for all }x\in\GG.
$$
Since each $t_i\in S$ and $\mu(t_i)\ge\rho$, we have
$$
p(z) = \prod_{i=1}^m \mu(t_i) \ge \rho^m = \rho^{|a^{-1}z|},
$$
and thus
\begin{equation}\label{eq:Green-comparison}
G(z,x) \le \rho^{-|a^{-1}z|}G(a,x)\qquad\text{for all }z,x\in\GG.
\end{equation}

\noindent
Fix $y\in\GG$ and set, for $s\in\GG$,
$$
H_y(s)\ :=\ \rho^{-|a^{-1}y|-|s|}.
$$
Given a sequence $(x_k)$ with $x_k\to\infty$, as in Lemma \ref{lem:limit-green}, we have
for each $k$ and each $s\in\GG$:
$$
\psi_k(ys)
=\frac{G(ys,x_k)}{G(a,x_k)}
\ \le\ \rho^{-|a^{-1}(ys)|}
\ \le\ \rho^{-|a^{-1}y|-|s|},
$$
where the first inequality is \eqref{eq:Green-comparison} with $z=ys$, and the second uses
the triangle inequality
$
|a^{-1}(ys)| \le |a^{-1}y| + |s|
$
together with the fact that $\rho\in(0,1]$ makes $r\mapsto\rho^r$ decreasing in $r$.
Hence
$$
\psi_k(ys)\ \le\ H_y(s)\qquad\forall k,s.
$$
By the exponential moment assumption,
$$
\sum_{s\in\GG}\mu(s)H_y(s)
 = \rho^{-|a^{-1}y|}\sum_{s\in\GG}\mu(s)\rho^{-|s|}
 < \infty.
$$
Thus the domination condition in Lemma \ref{lem:limit-green}(b) is satisfied for each fixed
$y$, and any pointwise subsequential limit of $(\psi_k)$ is $\mu$-harmonic.

\medskip\noindent
\emph{Finitely supported case.}
If $\mu$ is finitely supported, the same conclusion is immediate even without the
exponential moment assumption: for each fixed $y$ and $s$, the comparison
\eqref{eq:Green-comparison} yields
$$
\psi_k(ys) \le H_y(s) := \rho^{-|a^{-1}(ys)|},
$$
and $H_y$ satisfies $\sum_s \mu(s)H_y(s)<\infty$ because $\supp\mu$ is finite.
\end{proof}

\begin{prop}\label{prop:Delta_from_strong_Liouville}
Let $\GG$ be a finitely generated group and let $\mu$ be a symmetric, non-degenerate
probability measure on $\GG$ which generates a transient random walk on $\GG$.
Let $S \subset\supp\mu$ be a finite symmetric generating set for $\GG$, and set
$\rho:=\min_{s\in S}\mu(s)\in(0,1]$. Assume there exists $\rho_0<\rho$ such that
$$
\sum_{g \in \GG}\mu(g)\rho_0^{-|g|}<\infty.
$$
If $(\GG,\mu)$ has the \emph{strong Liouville} property, then for all $a,b\in\GG$ and for every exhaustion $S_n\nearrow\GG$ we have
$$
\Delta(S_n;a,b)\longrightarrow 0.
$$
In particular, if $\mu$ has superexponential moments, then the conclusion holds.
\end{prop}

\begin{proof}
First note that the assumption with $\rho_0<\rho$ implies
$$
\sum_{g\in\GG}\mu(g)\rho^{-|g|}
 \le \sum_{g\in\GG}\mu(g)\rho_0^{-|g|}<\infty,
$$
since $0<\rho_0<\rho<1$ gives $\rho^{-|g|}\le\rho_0^{-|g|}$ for all $g$.
Thus the hypothesis of Corollary \ref{cor:easy-domination} is satisfied (with the same
$\rho$ coming from the generating set $S$).

Suppose, toward a contradiction, that there exist $a,b\in\GG$ and an exhaustion
$(S_n)$ of $\GG$ with $\Delta(S_n;a,b)\not\to 0$. Then there exist $\varepsilon>0$,
a subsequence $(S_{n_k})$, and points $x_k\in\partial S_{n_k}$ such that
\begin{equation}\label{eq:Delta-counter}
\frac{|G(a,x_k)-G(b,x_k)|}{G(a,x_k)}\ \ge\ \varepsilon\qquad\text{for all }k.
\end{equation}
Since $S_n\nearrow\GG$ with $S_n$ finite, the boundary points $x_k\in\partial S_{n_k}$
escape every finite subset of $\GG$, so $x_k\to\infty$.

Define the normalised kernels
$$
\psi_k(v):=\frac{G(v,x_k)}{G(a,x_k)},\qquad v\in\GG.
$$
By \eqref{eq:Green-comparison} in the proof of Corollary \ref{cor:easy-domination} (with $z=v$), for each fixed $v\in\GG$ there exists a constant $C_v>0$ such that
$$
0\ \le\ \psi_k(v)\ =\ \frac{G(v,x_k)}{G(a,x_k)}\ \le\ C_v\qquad\text{for all }k.
$$
Since $\GG$ is countable, a diagonal subsequence argument now yields a subsequence of $(\psi_k)$
(still denoted as $(\psi_k)$) converging pointwise on $\GG$ to some function
$\psi:\GG\to[0,\infty)$.

By Lemma \ref{lem:limit-green}(a), $\psi$ is $\mu$-superharmonic. By Corollary \ref{cor:easy-domination} $\psi$ is $\mu$-harmonic. Moreover, by construction,
$$
\psi(a)
=\lim_{k\to\infty}\psi_k(a)
=\lim_{k\to\infty}\frac{G(a,x_k)}{G(a,x_k)}=1.
$$

From \eqref{eq:Delta-counter} we obtain
$$
|\psi_k(b)-1|
=\left|\frac{G(b,x_k)}{G(a,x_k)}-1\right|
\ \ge\ \varepsilon\qquad\text{for all }k,
$$
so passing to the limit gives
$$
|\psi(b)-1|\ \ge\ \varepsilon>0.
$$
Thus $\psi$ is a non-constant positive harmonic function on $(\GG,\mu)$, contradicting
the strong Liouville property. This contradiction shows that for every $a,b\in\GG$ and
every exhaustion $S_n\nearrow\GG$ we must have $\Delta(S_n;a,b)\to 0$.

Finally, if $\mu$ has superexponential moments, then for every $\zeta>0$ we have
$\sum_{g\in\GG}\mu(g)e^{\zeta|g|}<\infty$ by Definition \ref{def:smooth_measure}. Choose
$\zeta>0$ so small that $e^{-\zeta}<\rho$ and set $\rho_0:=e^{-\zeta}$. Then
$$
\sum_{g\in\GG}\mu(g)\rho_0^{-|g|} = \sum_{g\in\GG}\mu(g)e^{\zeta|g|}<\infty,
$$
so the hypothesis of the proposition is satisfied.
\end{proof}

We now prove the following converse to Proposition \ref{prop:Delta_from_strong_Liouville}.

\begin{theorem}[$\Delta\to 0$ along exhaustions $\Rightarrow$ strong Liouville]\label{lem:Pol_2.1_conv}
Let $\mu$ be a non-degenerate probability measure on $\GG$ which generates a transient
random walk. Assume that for all $a,b\in\GG$ and for every exhaustion $S\nearrow\GG$ we have
$$
\Delta(S;a,b)\longrightarrow 0.
$$
Then $(\GG,\mu)$ has the strong Liouville property.
\end{theorem}

\begin{proof}
Let $\xi$ be a point in the (minimal) Martin boundary of $\GG$, and fix a representing
sequence $(y_n)\subset\GG$ such that
$$
h(x)
:= K_e(x,\xi)
= \lim_{n\to\infty} K_e(x,y_n)
= \lim_{n\to\infty} \frac{G(x,y_n)}{G(e,y_n)}
\quad\text{for all }x\in\GG.
$$
We first show that $h$ is constant. Let $S=(S_n)_{n\ge1}$ be the exhausting sequence given by
$$
S_n:=\{x\in\GG:\ |x|<|y_n|\}.
$$
Then $y_n\notin S_n$, and since $y_n$ lies on a geodesic from $e$, its predecessor on this
geodesic lies in $S_n$, so $y_n\in\partial S_n$.

Fix arbitrary $a,b\in\GG$ and $\varepsilon>0$. By convergence of the Martin kernels in $x$,
there exists $n_1$ such that for all $n\ge n_1$,
$$
|h(a)-K_e(a,y_n)|<\varepsilon/6
\quad\text{and}\quad
|h(b)-K_e(b,y_n)|<\varepsilon/6.
$$
By the hypothesis $\Delta(S_n;e,a)\to0$ and $\Delta(S_n;e,b)\to0$, there exists $n_2$ such
that for all $n\ge n_2$,
$$
\Delta(S_n;e,a)<\varepsilon/3
\quad\text{and}\quad
\Delta(S_n;e,b)<\varepsilon/3.
$$
Take $m\ge\max\{n_1,n_2\}$. Then $y_m\in\partial S_m$, and
\begin{align*}
|h(a)-h(b)|
&\le |h(a)-K_e(a,y_m)| + |K_e(a,y_m)-K_e(b,y_m)| + |K_e(b,y_m)-h(b)|\\
&\le \frac{\varepsilon}{6}
 + \left|\frac{G(a,y_m)}{G(e,y_m)} - \frac{G(b,y_m)}{G(e,y_m)}\right|
 + \frac{\varepsilon}{6}\\
&\le \frac{\varepsilon}{3}
 + \left|\frac{G(a,y_m)}{G(e,y_m)} - 1\right|
 + \left|\frac{G(b,y_m)}{G(e,y_m)} - 1\right|.
\end{align*}
Since $y_m\in\partial S_m$, we can bound the last two terms using $\Delta$:
$$
\left|\frac{G(a,y_m)}{G(e,y_m)} - 1\right|
= \frac{|G(a,y_m)-G(e,y_m)|}{G(e,y_m)}
\le \Delta(S_m;e,a),
$$
and similarly for $b$. Hence
$$
|h(a)-h(b)|
\le \frac{\varepsilon}{3} + \Delta(S_m;e,a)+\Delta(S_m;e,b)
< \frac{\varepsilon}{3} + \frac{\varepsilon}{3} + \frac{\varepsilon}{3}
= \varepsilon.
$$
Since $\varepsilon>0$ and $a,b\in\GG$ were arbitrary, this shows $h(a)=h(b)$ for all
$a,b\in\GG$, i.e.\ $h$ is constant.

Thus every Martin kernel $K_e(\cdot,\xi)$ corresponding to a boundary point $\xi$ is
constant. By the Poisson-Martin representation theorem
(see, for example, \cite[Theorem~24.7]{Woess2000}), every positive $\mu$-harmonic function on
$\GG$ can be represented as an integral of such kernels against a finite measure on the
minimal Martin boundary. Since each kernel is constant, every positive harmonic function
must be constant. Hence $(\GG,\mu)$ has the strong Liouville property.
\end{proof}

\begin{remark}[Ancona inequalities and failure of $\Delta\to0$ on hyperbolic groups]\label{rem:Prop_4.2_free_hyp}
For non-elementary Gromov-hyperbolic groups with finitely supported, symmetric,
non-degenerate measures, one can see directly that $\Delta$ \emph{cannot} decay along
balls, providing a concrete obstruction to strong Liouville.

Indeed, suppose $\GG$ is non-elementary hyperbolic and $\mu$ is symmetric, finitely
supported, and non-degenerate. Ancona inequalities (see, e.g., \cite{Ancona1990,BHM2011})
assert that there exists $C>0$ such that for all points $x,y,z$ lying on a common
geodesic with $y$ between $x$ and $z$,
$$
G(x,z) \le CG(x,y)G(y,z).
$$
As recalled earlier (and in fact true for any transient irreducible random walk), we have
$G(e,x)\to 0$ as $|x|\to\infty$; see, for example, \cite[Lemma 3.2]{BHM2011} for a proof.

Choose $b\in\GG$ such that
$$
\varepsilon := 1 - CG(e,b)  > 0.
$$
Let $x$ lie on the geodesic ray from $e$ through $b$, with $b$ between $e$ and $x$. Applying
Ancona with $(x,y,z)=(e,b,x)$ gives
$$
G(e,x) \le CG(e,b)G(b,x),
$$
and hence
$$
\left|1 - \frac{G(e,x)}{G(b,x)}\right| \ge 1 - \frac{G(e,x)}{G(b,x)}\ge1 -CG(e,b)=\varepsilon.
$$
Taking $S=B(e,|x|)$, this shows that
$$
\Delta\big(B(e,|x|);e,b\big) \ge \varepsilon
$$
whenever $x$ lies sufficiently far along the geodesic ray through $b$. In particular,
$\Delta(B(e,r);e,b)$ does \emph{not} tend to $0$ as $r\to\infty$ on such groups, in line
with the presence of many non-trivial positive harmonic functions in the
hyperbolic/exponential-growth setting.
\end{remark}

Although this paper focuses exclusively on positive harmonic functions and the strong Liouville property, we must 
mention that the study of the Liouville property (that is, {\em bounded} harmonic functions being constant) is a very active and important area of research. To get an eclectic feel for the ideas, questions and methods involved in the latter, we refer the reader to \cite{Erschler2004, AErschler2004, FrischHartmanTamuzVahidi19,Gournay2016,Shalom2004} and references therein.

In a related vein, we also refer to the work of Lyons and Sullivan \cite{LyonsSullivan1984} on the Liouville and the strong Liouville properties of the covering space. It should be emphasised that it is not known if these properties depend only on the topology of the base manifold, or if the Riemannian metric plays a role. Let $p: M \to N$ be a normal Riemannian covering of a closed manifold, with deck transformation group $\Gamma$. The authors show that if $\Gamma$ is non-amenable, then there exist non-constant, bounded harmonic functions on $M$ \cite[Theorem 3]{LyonsSullivan1984}.
About the strong Liouville property, the authors show that if $\Gamma$ is virtually nilpotent, then any positive harmonic function on $M$ is constant \cite[Theorem 1]{LyonsSullivan1984}. In \cite[page 305]{LyonsSullivan1984}, the authors conjecture that $\Gamma$ is of exponential growth if and only if $M$ admits non-constant, positive harmonic functions. This is proved in \cite{BBE94}, \cite{BougerolElie1995}, under the assumption that $\Gamma$ is linear, that is, a closed subgroup of $GL_n(\RR)$, for some $n \in \NN$.

\subsection{Comparability of $\Delta$ with $\epsilon$ (the Amir-Kozma functional)}\label{section:Appendix_A}

We now make some comments about Question \ref{ques:eps_del_compare}, at least in some concrete settings.

Let $\GG$ be a finitely generated group, and fix a symmetric generating set $T$.
Fix once and for all a basepoint $o\in\GG$ (typically the identity element).
For a finite set $S\subset \GG$, define the (outer) $T$-boundary
$$
\pa S := \{ x \in \GG \setminus S \mid \exists y \in S, \exists t\in T: x=yt\}.
$$
Let $(X_n)_{n\ge 0}$ be the $\mu$-random walk on $\GG$
and assume it is \emph{transient}. Then the Green function
$$
G(a,x):=\EE_a\left[\sum_{n\ge0}\mathbf 1_{\{X_n=x\}}\right]
=\sum_{n\ge0}\PP_a[X_n=x]
$$
is finite for all $a,x \in \GG$.
For a finite $S\subset\GG$ and $a\in S$, let
$$
\tau_S:=\inf\{n\ge 1:\ X_n\notin S\}
$$
be the first exit time from $S$, and define the exit distribution
$$
\mu_S(a,z):=\PP_a[X_{\tau_S}=z],\qquad z\in\GG\setminus S.
$$

\begin{prop}[Green decomposition via the exit distribution]\label{prop:green_exit_decomp}
Assume the $\mu$-random walk on $\GG$ is transient.
Fix a finite set $S\subset \GG$, a starting point $a\in S$, and a point $x\in \GG\setminus S$.
Then
$$
G(a,x)=\sum_{z\in \GG\setminus S} \mu_S(a,z) G(z,x)
=\EE_a\big[G(X_{\tau_S},x)\big].
$$
\end{prop}

\begin{remark}[Finite range and support of the exit law]\label{rem:finite_range_exit_support}
If $\mu$ is finitely supported, $\supp(\mu)$ generates $\GG$, and we take $T=\supp(\mu)$, then $X_{\tau_S}\in \pa S$ almost surely
(the walk exits by a single $T$-step). Hence $\mu_S(a,\cdot)$ is supported on $\pa S$ and the sum in Proposition \ref{prop:green_exit_decomp} reduces to $\sum_{z\in\pa S}$.

More precisely, $Z:=\{z\notin S:\ \PP_a[X_{\tau_S}=z]>0\}$ is exactly the subset of $\pa S$ reachable as an exit point from $a$.
In particular, if $S$ is connected in the Cayley graph of $(\GG,T)$ (or if one replaces $S$ by the connected component of $a$ inside $S$),
then $Z=\pa S$.
\end{remark}

\begin{proof}[Proof of Proposition \ref{prop:green_exit_decomp}]
Since the $\mu$-random walk on $\GG$ is transient, $\tau_S<\infty$ almost surely.
Moreover, for $x\notin S$ and $n<\tau_S$ we have $X_n\in S$, hence $\mathbf 1_{\{X_n=x\}}=0$.
Therefore, 
$$
\sum_{n\ge 0}\mathbf 1_{\{X_n=x\}} = \sum_{n\ge \tau_S}\mathbf 1_{\{X_n=x\}}
= \sum_{k\ge 0}\mathbf 1_{\{X_{\tau_S+k}=x\}}.
$$
Taking $\EE_a$ gives
$$
G(a,x)= \EE_a\left[ \sum_{k\ge 0}\mathbf 1_{\{X_{\tau_S+k}=x\}} \right].
$$
By the strong Markov property at time $\tau_S$, for each $z\notin S$,
$$
\EE_a\left[ \sum_{k\ge 0}\mathbf 1_{\{X_{\tau_S+k}=x\}} \Big| X_{\tau_S}=z \right]
= \EE_z\left[ \sum_{k\ge 0}\mathbf 1_{\{X_k=x\}} \right]
= G(z,x).
$$
Hence, by the tower property,
$$
G(a,x)=\EE_a\left[\EE_a\left[ \sum_{k\ge 0}\mathbf 1_{\{X_{\tau_S+k}=x\}} \Big| X_{\tau_S} \right]\right]
=\sum_{z\in\GG\setminus S}\mu_S(a,z)G(z,x).
$$
\end{proof}

\subsubsection{The boundary Green matrix and invertibility (finite range, symmetric case)}

Assume now that $\mu$ is \emph{symmetric}, \emph{finitely supported}, and that $\supp(\mu)$ generates $\GG$; set $T=\supp(\mu)$.
For a finite $S\subset \GG$, define the linear operator $M_S:\RR^{\pa S}\to \RR^{\pa S}$ by
$$
(M_S c)(x):=\sum_{z\in \pa S} G(z,x)c(z),\qquad x\in\pa S.
$$
For $a\in S$ define the vectors
$$
g_S(a)(x):=G(a,x)\quad (x\in\pa S),
\qquad
\mu_S(a)(z):=\mu_S(a,z)\quad (z\in\pa S).
$$
Then Proposition \ref{prop:green_exit_decomp} and Remark \ref{rem:finite_range_exit_support} give the vector identity
\begin{equation}\label{eq:vector_green_exit}
g_S(a)=M_S\mu_S(a).
\end{equation}

\begin{lemma}[Why $M_S$ is symmetric positive definite]\label{lem:MS_spd}
Under the above assumptions (symmetric, transient walk), $M_S$ is symmetric and strictly positive definite, hence invertible.
\end{lemma}

\begin{proof}
Symmetry of $\mu$ implies $G(z,x)=G(x,z)$, hence $M_S$ is symmetric.

For strict positive definiteness, let $c\neq 0$ in $\RR^{\pa S}$ and define a finitely supported function
$f:\GG\to\RR$ by $f(z)=c(z)$ for $z\in\pa S$ and $f\equiv 0$ on $\GG\setminus \pa S$.
Let $P$ be the Markov operator
$$
(P\varphi)(y):=\sum_{t\in T}\mu(t)\varphi(yt).
$$
Define the potential $u:\GG\to\RR$ by
$$
u(y):=\sum_{z\in\GG} G(y,z)f(z)=\sum_{z\in\pa S} G(y,z)c(z).
$$
A standard computation shows that $(I-P)u=f$.

Now compute
$$
\langle f,u\rangle
:=\sum_{y\in\GG} f(y)u(y)
=\sum_{z,x\in\pa S} c(z)G(z,x)c(x),
$$
which is exactly the quadratic form associated to $M_S$.
Using symmetry of $\mu$, 
$$
\langle (I-P)u,u\rangle
=\frac12\sum_{y\in\GG}\sum_{t\in T}\mu(t)\big(u(y)-u(yt)\big)^2 \ge0.
$$
If the right-hand side were $0$, then $u(y)=u(yt)$ for all $y$ and all $t$ with $\mu(t)>0$; since $T=\supp(\mu)$ generates $\GG$,
this forces $u$ to be constant on $\GG$, hence $f=(I-P)u=0$, contradicting $c\neq 0$.
Therefore $\sum_{z,x\in\pa S} c(z)G(z,x)c(x)>0$ for all $c\neq 0$.
\end{proof}

\noindent
By Lemma \ref{lem:MS_spd}, \eqref{eq:vector_green_exit} yields
$$
\mu_S(a)=M_S^{-1}g_S(a).
$$
Using the $\ell^\infty$ norm on vectors indexed by $\pa S$ and the induced operator norm $\|\cdot\|_{\infty\to\infty}$,
we obtain for $a,b\in S$:
$$
\|\mu_S(a)-\mu_S(b)\|_\infty
\le \|M_S^{-1}\|_{\infty\to\infty}\|g_S(a)-g_S(b)\|_\infty,
\qquad
\|g_S(a)-g_S(b)\|_\infty
\le \|M_S\|_{\infty\to\infty}\|\mu_S(a)-\mu_S(b)\|_\infty.
$$
Hence, if along an exhaustion $(S_k)$ one has a uniform conditioning bound
$$
\sup_k  \|M_{S_k}\|_{\infty\to\infty}\|M_{S_k}^{-1}\|_{\infty\to\infty}  <  \infty,
$$
then the max-absolute differences of $\mu_{S_k}$ and $g_{S_k}$ are comparable (up to uniform multiplicative constants).

To pass to the \emph{relative} functionals $\epsilon,\Delta$, one typically wants a \emph{pointwise Poisson/Green factorisation}:

\begin{lemma}[Cancellation in the relative functionals under Poisson/Green factorisation]\label{lem:eps_Delta_cancellation}
Fix a finite set $S\subset \GG$, basepoints $a,b\in S$, and a boundary point $x\in \pa S$.

\smallskip
\noindent\textbf{(i) Exact factorisation.}
Assume that for all $a\in S$ and all $x\in\pa S$ one has the pointwise identity
\begin{equation}\label{eq:exact_factorisation}
\mu_S(a,x)=c_S(x)\frac{G(a,x)}{G(o,x)},
\end{equation}
where $c_S(x)>0$ depends only on $(S,x)$ (in particular, not on $a$). Then for every $x\in\pa S$,
\begin{equation}\label{eq:pointwise_eps_equals_Delta}
\frac{|\mu_S(a,x)-\mu_S(b,x)|}{\mu_S(a,x)}
=\frac{|G(a,x)-G(b,x)|}{G(a,x)}.
\end{equation}
Consequently,
\begin{equation}\label{eq:eps_equals_Delta}
\epsilon(S;a,b)=\Delta(S;a,b),
\end{equation}
where $\epsilon$ and $\Delta$ are evaluated over the same boundary index set (here $\pa S$).

\smallskip
\noindent\textbf{(ii) Approximate factorisation.}
Assume instead that there exists $\eta\in[0,1)$ such that for all relevant $a\in S$ and all $x\in\pa S$,
\begin{equation}\label{eq:approx_factorisation}
\mu_S(a,x)=c_S(x)\frac{G(a,x)}{G(o,x)}(1+e_a(x)),
\qquad |e_a(x)|\le \eta,
\end{equation}
with $c_S(x)>0$ depending only on $(S,x)$. For each $x\in\pa S$ set
$$
\Delta_x(S;a,b):=\frac{|G(a,x)-G(b,x)|}{G(a,x)}.
$$
Then for every $x\in\pa S$ one has the two-sided bounds
\begin{equation}\label{eq:pointwise_eps_vs_Delta}
\frac{(1-\eta)\Delta_x(S;a,b)-2\eta}{1+\eta}\le\frac{|\mu_S(a,x)-\mu_S(b,x)|}{\mu_S(a,x)}\le\frac{(1+\eta)\Delta_x(S;a,b)+2\eta}{1-\eta}.
\end{equation}
In particular, taking maxima over $x\in\pa S$ yields
\begin{equation}\label{eq:eps_vs_Delta_global}
\frac{(1-\eta)\Delta(S;a,b)-2\eta}{1+\eta}
 \le
\epsilon(S;a,b)
 \le
\frac{(1+\eta)\Delta(S;a,b)+2\eta}{1-\eta},
\end{equation}
where $\Delta(S;a,b)=\max_{x\in\pa S}\Delta_x(S;a,b)$.
\end{lemma}

\begin{proof}
\textbf{(i)} Starting from \eqref{eq:exact_factorisation},
$$
\mu_S(a,x)-\mu_S(b,x)
=c_S(x)\frac{G(a,x)-G(b,x)}{G(o,x)}.
$$
Dividing by $\mu_S(a,x)=c_S(x)G(a,x)/G(o,x)$ gives \eqref{eq:pointwise_eps_equals_Delta}, and taking maxima over $x\in\pa S$ gives \eqref{eq:eps_equals_Delta}.

\smallskip
\noindent\textbf{(ii)} Under \eqref{eq:approx_factorisation}, the common prefactor $c_S(x)/G(o,x)$ cancels:
$$
\frac{|\mu_S(a,x)-\mu_S(b,x)|}{\mu_S(a,x)}
=
\frac{\big|G(a,x)(1+e_a(x))-G(b,x)(1+e_b(x))\big|}{G(a,x)(1+e_a(x))}.
$$
Write $G_a:=G(a,x)$, $G_b:=G(b,x)$, $e_a:=e_a(x)$, $e_b:=e_b(x)$, and $\Delta_x:=|G_a-G_b|/G_a$.
Using $|1+e_a|\ge 1-\eta$ and $|1+e_a|\le 1+\eta$, together with
$$
\frac{G_b}{G_a}\le 1+\Delta_x \qquad\text{(since $G_b=G_a+(G_b-G_a)$),}
$$
we obtain the upper bound
\begin{align*}
\frac{|\mu_S(a,x)-\mu_S(b,x)|}{\mu_S(a,x)}
&\le \frac{|G_a-G_b|+\eta G_a+\eta G_b}{G_a(1-\eta)}\\
&\le \frac{\Delta_x+\eta+\eta(1+\Delta_x)}{1-\eta}
= \frac{(1+\eta)\Delta_x+2\eta}{1-\eta}.
\end{align*}
For the lower bound, use $|A+B|\ge |A|-|B|$ and $|1+e_a|\le 1+\eta$:
\begin{align*}
\frac{|\mu_S(a,x)-\mu_S(b,x)|}{\mu_S(a,x)}
&\ge \frac{|G_a-G_b|-\eta G_a-\eta G_b}{G_a(1+\eta)}\\
&\ge \frac{\Delta_x-\eta-\eta(1+\Delta_x)}{1+\eta}
= \frac{(1-\eta)\Delta_x-2\eta}{1+\eta}.
\end{align*}
This yields \eqref{eq:pointwise_eps_vs_Delta}. Taking maxima over $x\in\pa S$ gives \eqref{eq:eps_vs_Delta_global}.
\end{proof}

\subsection{Settings where comparability holds }

\subsubsection{Cayley balls in $\Z^d$ ($d\ge 3$)}
Let $\GG=\Z^d$ with the standard nearest-neighbour generating set and let $\mu$ be simple random walk.
Let $S=B(o,R)$ be the discrete Euclidean ball and consider fixed $a,b\in\Z^d$ (independent of $R$).

\smallskip
\noindent
For uniformly elliptic symmetric nearest-neighbour kernels on $\Z^d$, one has Lipschitz dependence
of the exit distribution on the starting point; see \cite[Lemma~4.1]{BaurBolthausen2015}.
Write $\mu_R(x,z):=\mu_{B(o,R)}(x,z)$ for $z\in\pa B(o,R)$.
Lemma~4.1(i)-(ii) of \cite{BaurBolthausen2015} (applied to the symmetric kernel $p=p_0$ of simple random walk and their discrete ball $V_L$)
gives: for any fixed $\eta\in(0,1)$ there exists $C_\eta >0$ such that for all large $R$,
for all $x,x'\in B(o,\eta R)$ and all $z\in\pa B(o,R)$,
\begin{equation}\label{eq:BB_Lip_checked}
|\mu_R(x,z)-\mu_R(x',z)|  \le C_\eta|x-x'|R^{-d},
\end{equation}
and moreover there exists $c_\eta>0$ such that
\begin{equation}\label{eq:BB_bounds_checked}
c_\eta R^{-(d-1)}  \le \mu_R(x,z) \le C_\eta R^{-(d-1)}
\qquad (x\in B(o,\eta R), z\in\pa B(o,R)).
\end{equation}
Fix $a\in\Z^d$. For all sufficiently large $R$, we have $a,o\in B(o,\eta R)$, hence by \eqref{eq:BB_Lip_checked}-\eqref{eq:BB_bounds_checked},
\begin{equation}\label{eq:exit-ratio-estimate}
\left|\frac{\mu_R(a,z)}{\mu_R(o,z)}-1\right|=\frac{|\mu_R(a,z)-\mu_R(o,z)|}{\mu_R(o,z)}\le\frac{C_\eta |a|R^{-d}}{c_\eta R^{-(d-1)}}=\frac{C_1}{R},
\end{equation}
uniformly in $z\in\pa B(o,R)$. 

\medskip
\noindent
For simple random walk on $\Z^d$ ($d\ge 3$), let $G$ denote the Green's function and (with slight abuse of notation) let $G(z):=G(o,z)=G(z,o)$, $o$ being the identity element in $\ZZ^d$. Then, $G(a,z)=G(z-a)$ for all $z,a\in\ZZ^d$.
Fix $a\in\Z^d$ and let $z\in\pa B(o,R)$, so $|z|\asymp R$. Choose a nearest-neighbour lattice path
$(z_0,\dots,z_m)$ from $z_0=z$ to $z_m=z-a$ with $m\le |a|_1$.
By telescoping,
$$
|G(z-a)-G(z)|
\le \sum_{k=0}^{m-1} |G(z_{k+1})-G(z_k)|
\le \sum_{k=0}^{m-1}\max_{1\le j\le d}|\nabla_j G(z_k)|,
$$
where $\nabla_j G(x):=G(x+e_j)-G(x)$.
By \cite[Corollary 4.3.3]{LawlerLimic2010}, $|\nabla_j G(x)|=O(|x|^{1-d})$ as $|x|\to\infty$, uniformly in $j$.
Moreover, for $R\ge 2|a|_1$ we have $|z_k|\ge |z|-|a|_1\ge R/2$ for all $k$, hence
$$
|G(z-a)-G(z)| \le C|a|R^{1-d}.
$$
On the other hand, \cite[Theorem 4.3.1]{LawlerLimic2010} gives
$$
G(z)=c_d|z|^{2-d}+O(|z|^{-d}) \asymp R^{2-d}.
$$
Therefore,
\begin{equation}\label{eq:Green-ratio-estimate}
\left|\frac{G(a,z)}{G(o,z)}-1\right|
=\left|\frac{G(z-a)-G(z)}{G(z)}\right|
\le \frac{C|a|R^{1-d}}{cR^{2-d}}
 =\frac{C_2}{R},
\end{equation}
uniformly in $z\in\pa B(o,R)$ (for fixed $a$). 

For each $z\in\pa B(o,R)$ define
$$
c_R(z):=\mu_R(o,z)\quad e_R(z):=\frac{\mu_R(a,z)}{\mu_R(o,z)}-1, \qquad d_R(z):=\frac{G(a,z)}{G(o,z)}-1.
$$
Then \eqref{eq:exit-ratio-estimate}-\eqref{eq:Green-ratio-estimate} give
$$
\sup_{z\in\pa S_R}|e_R(z)|\le \frac{C_1}{R},\qquad \sup_{z\in\pa S_R}|d_R(z)|\le \frac{C_2}{R}.
$$
For each $z$ we have the exact identity
\begin{align*}
\frac{\mu_R(a,z)}{c_R(z)G(a,z)/G(o,z)}
&=
\frac{\mu_R(a,z)/\mu_R(o,z)}{G(a,z)/G(o,z)}
=
\frac{1+e_R(z)}{1+d_R(z)}.
\end{align*}
Hence
$$
\theta_R(z):=\frac{\mu_R(a,z)}{c_R(z)G(a,z)/G(o,z)}-1=\frac{e_R(z)-d_R(z)}{1+d_R(z)}.
$$
For $R>C_2$ we have $|d_R(z)|\le C_2/R<1$ uniformly in $z$, hence
$$
|1+d_R(z)|\ge 1-|d_R(z)|\ge 1-\frac{C_2}{R}=\frac{R-C_2}{R}.
$$
Therefore, uniformly in $z\in\pa S_R$,
$$
|\theta_R(z)|\le \frac{|e_R(z)|+|d_R(z)|}{|1+d_R(z)|}\le\frac{(C_1+C_2)R^{-1}}{(R-C_2)R^{-1}}=\frac{C_1+C_2}{R-C_2}.
$$
Multiplying by $c_R(z)G(a,z)/G(o,z)$ gives $\mu_R(a,z)\le \frac{C_1+C_2}{R-C_2}\cdot c_R(z)\frac{G(a,z)}{G(o,z)}$.
Consequently, $\epsilon(B(o,R);a,b)$ and $\Delta(B(o,R);a,b)$ are asymptotically comparable as $R\to\infty$ (for fixed $a,b$).

\subsubsection{Hyperbolic groups} Let $\GG$ be a non-elementary word-hyperbolic group and let $\mu$
be symmetric and finitely supported. Fix $a\in\GG$,  $R\ge 1$ and set $S=B(e,R)$. In this case we believe that the uniform Ancona inequalities \cite[Definition 2.2, Theorem 2.3]{Gouezel2014} will be useful in obtaining a bound of the form 
\begin{equation}\label{eq:Martin_vs_finite_ratio}
K(a,\xi_x) \asymp \frac{G(a,x)}{G(e,x)}.
\end{equation}
where $\gamma_x:[0,\infty)\to\GG$ is a geodesic ray extending the word geodesic segment $[e,x]$, $x\in\partial S$. We think that the boundary Harnack principle \cite{Ancona1990}, \cite{KemperLohkamp2022} can be used to obtain  
\begin{equation}\label{eq:Poisson_RN_Martin}
\frac{\mu_{B(o,R)}(a,x)}{\mu_{B(e,R)}(e,x)}\ \asymp\ K(a,\xi_x),
\end{equation}
uniformly in $R$ and $x\in\pa B(o,R)$.
This should give us comparability between $\Delta$ and $\epsilon$, though we have not checked the details.

\subsection{Settings where comparability fails}

\begin{itemize}
\item \emph{Recurrence} ($\Z,\Z^2$): $G(a,x)=\infty$, so $\Delta$ (as defined from the full Green function) is undefined,
though $\epsilon$ still makes sense as an exit-probability functional.
\smallskip
\item \emph{Infinite support measures}: If $\mu$ has infinite support, then in general $X_{\tau_S}$ need not lie in $\pa S$ (it may jump far outside $S$), so the finite-dimensional boundary matrix $M_S$
is no longer the right object. Comparability of $\epsilon$ and $\Delta$ is then governed instead by Poisson/Martin kernel factorisations and boundary Harnack principles
for the corresponding jump process.
\smallskip
\item \emph{Hyperbolic groups}: If $\mu$ has \emph{superexponential tails}, then strong Ancona inequalities hold and the Martin boundary coincides with the Gromov boundary;
see \cite{Gouezel2015}. In this regime one again obtains Green/Martin factorisations implying comparability of $\epsilon$ and $\Delta$ along balls. If $\mu$ has only exponential tails, Ancona inequalities may fail and the Martin boundary can be more complicated; comparability of $\epsilon$ and $\Delta$
need not hold in general (see \cite{Gouezel2015} for examples/pathologies).

\end{itemize}

\section{Martin boundary collapse: a unified, pointwise framework for $\Delta\to0$}\label{sec:del_framework}

Let $(\GG,d)$ be a connected bounded-degree graph, with adjacency relation $x\sim y$ and graph distance $d$.
Let $(X_n)_{n\ge 0}$ be a \emph{lazy} symmetric irreducible Markov chain on $\GG$ with one-step transition kernel
$$
p(x,y):=\PP_x[X_1=y],\qquad p(x,y)=p(y,x),\qquad \inf_{x\in\GG}p(x,x)>0.
$$
Write
$$
p_0(x,y)=\delta_{xy},\qquad p_n(x,y)=\PP_x[X_n=y]\quad(n\ge 1),\qquad
G(x,y)=\sum_{n\ge 0}p_n(x,y).
$$
(In the Cayley graph case with a symmetric nondegenerate step distribution $\mu$, one implies this setting by
replacing $\mu$ with its lazy version $\mu_{\mathrm{lazy}}:=\tfrac12(\delta_e+\mu)$; then $p_n(x,y)=\mu_{\mathrm{lazy}}^{*n}(x^{-1}y)$.)

For $a,b \in \GG$ and a finite set $S\subset\GG$, define the (outer) vertex boundary
$$
\pa S :=\{x\notin S: \exists y\in S \text{ with } x\sim y\},
$$
and
\begin{align*}
R(S;a,b)&=\dist(\{a,b\},\pa S),\\
\Delta(S;a,b)&=\max_{x\in\pa S}\frac{|G(a,x)-G(b,x)|}{G(a,x)}.
\end{align*}

\smallskip
Fix a \emph{near-diagonal parameter} $\kappa>0$.
Fix a strictly increasing intrinsic scale $\rho:\{1,2,\dots\}\to(0,\infty)$ such that $\lim_{n\to\infty}\rho(n)=\infty$,
and an increasing reference volume $\Vol:(0,\infty)\to(0,\infty)$.
Set the near-diagonal tail
$$
A(m):=\sum_{n\ge m}\Vol(\rho(n))^{-1}\qquad(m\ge 1),
$$
and the near-diagonal threshold
$$
n_-(r):=\min\{n\in \NN:\ \rho(n)\ge r/\kappa\}\qquad(r\ge 0),
$$
which is well-defined and finite since $\rho(n)\uparrow\infty$.

\begin{remark}[Transience vs.\ resolvent]\label{rem:transient-resolvent}
Under \textup{(HK)} below, for every $x\in\GG$ we have
$$
c_0A(1)\ \le\ G(x,x)-1\ \le\ C_0A(1).
$$
In particular, $G$ is finite (the chain is transient) if and only if $A(1)<\infty$.
Our main estimates are stated in this transient regime.
\end{remark}

We will invoke the following scale-robust assumptions (see the discussion preceding Theorem \ref{thm:Delta-IH}):

\begin{itemize}
\item[(HK)] \emph{Global upper and near-diagonal lower bounds:}
There exists $c_0,C_0>0$ such that for all $n\ge 1$ and all $x,y\in \GG$,
\begin{equation}\label{eq:HK-ud}
p_n(x,y) \le \frac{C_0}{\Vol(\rho(n))},\qquad
d(x,y)\le \kappa\rho(n)  \Rightarrow p_n(x,y) \ge \frac{c_0}{\Vol(\rho(n))}.
\end{equation}
\emph{(Remark.)} In nearest-neighbor examples one typically normalises $\rho(1)$ and chooses $\kappa$ so that
the constraint $d(x,y)\le \kappa\rho(1)$ lies within the one-step jump range.

\item[(OD)] \emph{Off-diagonal upper envelope at the $\rho$-scale:}
There exist $C_{\mathrm{OD}}>0$ and a nonincreasing function $\Phi:[0,\infty)\to(0,1]$
with $\Phi(t)\to0$ as $t\to\infty$ such that for all $n\ge 1$ and all $x,y\in\GG$,
\begin{equation}\label{eq:OD}
p_n(x,y) \le \frac{C_{\mathrm{OD}}}{\Vol(\rho(n))}\Phi\Big(\frac{d(x,y)}{\rho(n)}\Big).
\end{equation}

\item[(IH)] \emph{Interior H\"older in the start variable at the $\rho$-scale:}
There exist constants $\theta\in(0,\kappa/4]$, $\alpha\in(0,1]$ and $C_{\mathrm H}>0$
such that for every $n\ge 1$, every $z_0,y\in\GG$, and all $z,z'\in B(z_0,\theta\rho(n))$,
\begin{equation}\label{eq:IH}
|p_n(z,y)-p_n(z',y)|
\ \le\ C_{\mathrm H}\Big(\frac{d(z,z')}{\rho(n)}\Big)^{\alpha}
\sup_{w\in B(z_0,\theta\rho(n))}p_n(w,y).
\end{equation}
By \eqref{eq:HK-ud}, $\sup_{w\in B(z_0,\theta\rho(n))} p_n(w,y)\le C_0/\Vol(\rho(n))$.

\item[(TR$_\alpha$)] \emph{Off-diagonal Tauberian comparability at exponent $\alpha$:}
There exists $C_{\mathrm T}>0$ such that for all $r\ge 1$, with $m:=n_-(r)$,
\begin{equation}\label{eq:TRa}
\sum_{1\le n<m}\frac{1}{\Vol(\rho(n))\rho(n)^\alpha}\Phi\Big(\frac{r}{2\rho(n)}\Big)\le C_{\mathrm T}\rho(m)^{-\alpha}A(m),\qquad A(m):=\sum_{n\ge m}\frac{1}{\Vol(\rho(n))}.
\end{equation}
\end{itemize}

\begin{lemma}[Universal H\"{o}lder-from-(IH)]\label{lem:UH2}
Under the assumptions \eqref{eq:HK-ud} and \eqref{eq:IH} there exists $C_{\mathrm{UH}}>0$ such that for all $n\in \NN$ and all $x,a,b\in\GG$,
\begin{equation}\label{eq:UH}
|p_n(a,x)-p_n(b,x)| \le C_{\mathrm{UH}}\Big(\frac{d(a,b)}{\rho(n)}\Big)^{\alpha}\frac{1}{\Vol(\rho(n))}.
\end{equation}
\end{lemma}

\begin{proof}
Fix $n\in \NN$ and $a,b,x\in\GG$, and write $r:=d(a,b)$.

If $r\le \theta\rho(n)$, then $a,b\in B(a,\theta\rho(n))$ and we may apply \eqref{eq:IH}
with $z=a$, $z'=b$, $z_0=a$, $y=x$, obtaining
$$
|p_n(a,x)-p_n(b,x)|
\le C_{\mathrm H}\Big(\frac{r}{\rho(n)}\Big)^{\alpha}
\sup_{w\in B(a,\theta\rho(n))}p_n(w,x)
\le C_{\mathrm H}C_0\Big(\frac{r}{\rho(n)}\Big)^{\alpha}\frac{1}{\Vol(\rho(n))},
$$
where the last inequality uses \eqref{eq:HK-ud}.

If $r> \theta\rho(n)$, then by \eqref{eq:HK-ud},
$$
|p_n(a,x)-p_n(b,x)|\le p_n(a,x)+p_n(b,x)\le \frac{2C_0}{\Vol(\rho(n))}.
$$
Since $r/\rho(n)>\theta$, we have $1\le \theta^{-\alpha}(r/\rho(n))^\alpha$, hence
$$
|p_n(a,x)-p_n(b,x)|
\le \frac{2C_0}{\theta^\alpha}\Big(\frac{r}{\rho(n)}\Big)^{\alpha}\frac{1}{\Vol(\rho(n))}.
$$

Taking $C_{\mathrm{UH}}:=\max\{C_{\mathrm H}C_0,\ 2C_0/\theta^\alpha\}$ yields \eqref{eq:UH}.
\end{proof}

\subsection{Quantitative $\Delta$-decay: the pointwise theorem}

\begin{theorem}[Quantitative $\Delta$-decay from (HK)+(OD)+(IH)+(TR$_\alpha$)]\label{thm:Delta-IH}
Assume $A(1)<\infty$ (equivalently, the chain is transient under \textup{(HK)}),
and assume \eqref{eq:HK-ud}, \eqref{eq:OD}, \eqref{eq:IH}, and \eqref{eq:TRa} for some
$\alpha\in(0,1]$. Then there exists $C>0$ such that for every finite
$S\subset\GG$ and all $a,b\in\GG\setminus \pa S$ with $d(a,b)\le R(S;a,b)$,
\begin{equation}\label{eq:Delta-bound}
\Delta(S;a,b)
\ \le\ C\left(\frac{d(a,b)}{R(S;a,b)}\right)^{\alpha}.
\end{equation}
In particular, for each fixed $a,b$ we have $\Delta(S;a,b)\to0$ along any exhaustion
$S\nearrow\GG$ (since $R(S;a,b)\to\infty$ and eventually $d(a,b)\le R(S;a,b)$).
\end{theorem}

\begin{proof}
Fix a finite set $S\subset\GG$ and $a,b\in\GG\setminus\pa S$ with $d(a,b)\le R(S;a,b)$, and write
$$
R:=R(S;a,b)=\dist(\{a,b\},\pa S)\ge 1.
$$
Let $x\in\pa S$ be arbitrary and set $r:=d(a,x)\ge R$. Let $m:=n_-(r)$.

Since $x\in\pa S$ and $a,b\notin\pa S$, we have $x\neq a$ and $x\neq b$, hence
$p_0(a,x)=p_0(b,x)=0$ and therefore
$$
|G(a,x)-G(b,x)|\ \le\ \sum_{n\ge 1}|p_n(a,x)-p_n(b,x)|.
$$
Split at $n=m$:
$$
|G(a,x)-G(b,x)| \le \underbrace{\sum_{1\le n<m}|p_n(a,x)-p_n(b,x)|}_{I_1} + \underbrace{\sum_{n\ge m}|p_n(a,x)-p_n(b,x)|}_{I_2}.
$$

\medskip\noindent
For $n\ge m$ we have $r\le\kappa\rho(n)$, so Lemma \ref{lem:UH2} yields
$$
|p_n(a,x)-p_n(b,x)|
 \le C_{\mathrm{UH}}\Big(\frac{d(a,b)}{\rho(n)}\Big)^{\alpha}\frac{1}{\Vol(\rho(n))}.
$$
Since $\rho$ is increasing and $\rho(n)\ge\rho(m)\ge r/\kappa$ for $n\ge m$, we get
$$
\frac{d(a,b)^\alpha}{\rho(n)^\alpha}
 \le \kappa^\alpha\left(\frac{d(a,b)}{r}\right)^{\alpha}.
$$
Hence
$$
I_2
\le C_{\mathrm{UH}}\kappa^\alpha\left(\frac{d(a,b)}{r}\right)^\alpha
\sum_{n\ge m}\frac{1}{\Vol(\rho(n))}.
$$
Using the near-diagonal lower bound in \eqref{eq:HK-ud} (valid since $d(a,x)=r\le\kappa\rho(n)$ for all $n\ge m$),
$$
p_n(a,x) \ge \frac{c_0}{\Vol(\rho(n))}\qquad(n\ge m),
$$
so
$$
I_2
\le \frac{C_{\mathrm{UH}}\kappa^\alpha}{c_0}\left(\frac{d(a,b)}{r}\right)^\alpha
\sum_{n\ge m}p_n(a,x)
\le \frac{C_{\mathrm{UH}}\kappa^\alpha}{c_0}\left(\frac{d(a,b)}{r}\right)^\alpha G(a,x).
$$

\medskip\noindent
For $n<m$ we have $\rho(n)<r/\kappa$. For any $w\in B(a,\theta\rho(n))$,
$$
d(w,x) \ge d(a,x)-d(a,w) \ge r-\theta\rho(n)
 \ge r-\theta(r/\kappa) \ge \frac r2,
$$
using $\theta\le\kappa/4$.
Hence by \eqref{eq:OD} and monotonicity of $\Phi$,
$$
\sup_{w\in B(a,\theta\rho(n))}p_n(w,x)
 \le \frac{C_{\mathrm{OD}}}{\Vol(\rho(n))}
\Phi\Big(\frac{r}{2\rho(n)}\Big).
$$
If $d(a,b)\le \theta\rho(n)$, apply \eqref{eq:IH} on $B(a,\theta\rho(n))$ to get
$$
|p_n(a,x)-p_n(b,x)|
\le C_{\mathrm H}\Big(\frac{d(a,b)}{\rho(n)}\Big)^{\alpha}
\sup_{w\in B(a,\theta\rho(n))}p_n(w,x).
$$
If $d(a,b)>\theta\rho(n)$, then $1\le \theta^{-\alpha}(d(a,b)/\rho(n))^\alpha$, and since
$d(a,b)\le R$ and $x\in\pa S$ we have $d(b,x)\ge R$ and moreover $d(b,x)\ge r/2$
(indeed, if $r\le 2R$ then $r/2\le R\le d(b,x)$; while if $r>2R$ then
$d(b,x)\ge d(a,x)-d(a,b)\ge r-R>r/2$). Thus \eqref{eq:OD} yields
$$
|p_n(a,x)-p_n(b,x)|
\le p_n(a,x)+p_n(b,x)\le \frac{2C_{\mathrm{OD}}}{\Vol(\rho(n))}\Phi\Big(\frac{r}{2\rho(n)}\Big)\le \frac{2C_{\mathrm{OD}}}{\theta^\alpha}\Big(\frac{d(a,b)}{\rho(n)}\Big)^{\alpha}\frac{1}{\Vol(\rho(n))}\Phi\Big(\frac{r}{2\rho(n)}\Big).
$$
In both cases we obtain
$$
|p_n(a,x)-p_n(b,x)|
 \le C_*\Big(\frac{d(a,b)}{\rho(n)}\Big)^{\alpha}\frac{1}{\Vol(\rho(n))}
\Phi\Big(\frac{r}{2\rho(n)}\Big),
$$
for $C_*:=\max\{C_{\mathrm H}C_{\mathrm{OD}}, 2C_{\mathrm{OD}}/\theta^\alpha\}$.
Summing over $1\le n<m$ and using \eqref{eq:TRa} yields
$$
I_1
\le C_*d(a,b)^\alpha\sum_{1\le n<m}\frac{1}{\Vol(\rho(n))\rho(n)^\alpha}
\Phi\Big(\frac{r}{2\rho(n)}\Big)
\le C_*C_{\mathrm T}d(a,b)^\alpha\rho(m)^{-\alpha}A(m).
$$
Using $\rho(m)\ge r/\kappa$ gives
$$
I_1\le C_*C_{\mathrm T}\kappa^\alpha\Big(\frac{d(a,b)}{r}\Big)^\alpha A(m).
$$
On the other hand, by \eqref{eq:HK-ud},
$$
G(a,x)\ \ge\ \sum_{n\ge m}p_n(a,x)\ge c_0\sum_{n\ge m}\frac{1}{\Vol(\rho(n))}=c_0 A(m).
$$
Thus
$$
\frac{I_1}{G(a,x)} \le \frac{C_*C_{\mathrm T}\kappa^\alpha}{c_0}\left(\frac{d(a,b)}{r}\right)^\alpha.
$$

\medskip\noindent
Combining the estimates on $I_1$ and $I_2$, we get
$$
\frac{|G(a,x)-G(b,x)|}{G(a,x)}\le \frac{I_1}{G(a,x)}+\frac{I_2}{G(a,x)}\le C\left(\frac{d(a,b)}{r}\right)^\alpha\le C\left(\frac{d(a,b)}{R}\right)^\alpha,
$$
for a constant $C>0$ depending only on the constants in (HK), (OD), (IH), (TR$_\alpha$).
Taking the maximum over $x\in\pa S$ yields \eqref{eq:Delta-bound}.
\end{proof}

\subsubsection{When \textup{(TR$_\alpha$)} is automatic}

On Ahlfors/polynomial scales, \textup{(TR$_\alpha$)} follows from an explicit off-diagonal envelope
$\Phi$ by a dyadic decomposition.

\begin{lemma}[Ahlfors scales + stretched-exponential envelope $\Rightarrow$ \textup{(TR$_\alpha$)}]\label{lem:TR-Ahlfors-subG}
Assume $\Vol(r)\asymp r^{d_*}$ and $\rho(n)\asymp n^{1/\gamma}$ with $d_*/\gamma>1$.
Assume \eqref{eq:OD} holds with an envelope satisfying
$$
\Phi(t) \le C\exp\big(-ct^{\eta}\big)\qquad (t\ge 0)
$$
for some $c,C,\eta>0$. Then \textup{(TR$_\alpha$)} holds for every $\alpha\in(0,1]$.
\end{lemma}

\begin{proof}
Fix $\alpha\in(0,1]$ and let $r\ge 1$, $m:=n_-(r)$. Under the Ahlfors hypotheses,
$$
\frac{1}{\Vol(\rho(n))\rho(n)^\alpha} \asymp n^{-p},
\qquad p:=\frac{d_*+\alpha}{\gamma}>1,
$$
and
$$
\rho(m)^{-\alpha}A(m) \asymp m^{-{\alpha}/{\gamma}}\sum_{n\ge m}n^{-d_*/\gamma}
 \gtrsim m^{-{\alpha}/{\gamma}}m^{1-d_*/\gamma}
 = m^{1-p}.
$$
Write
$$
S:=\sum_{1\le n<m}n^{-p}\Phi\Big(\frac{r}{2\rho(n)}\Big).
$$

Since $\rho(n)\asymp n^{1/\gamma}$, fix constants $0<c_\rho\le C_\rho$ such that
$c_\rho n^{1/\gamma}\le \rho(n)\le C_\rho n^{1/\gamma}$ for all $n\ge 1$.
Let $m=n_-(r)$, so $\rho(m-1)<r/\kappa\le \rho(m)$ (for $m\ge 2$; if $m=1$ then $S=0$ and there is nothing to prove).
These inequalities imply $m\asymp r^\gamma$ and hence $\rho(m)\asymp r$.
Therefore, for $n<m$,
$$
\frac{r}{\rho(n)} \asymp \frac{\rho(m)}{\rho(n)}\asymp \Big(\frac{m}{n}\Big)^{1/\gamma}.
$$

Decompose dyadically: for $j\ge0$, let $I_j:=\{n:\ 2^{-(j+1)}m\le n<2^{-j}m\}$.
For $n\in I_j$ we have $m/n\ge 2^j$, hence
$$
\Phi\Big(\frac{r}{2\rho(n)}\Big)\ \le\ C\exp\big(-c'2^{j\eta/\gamma}\big),
$$
and also $\sum_{n\in I_j}n^{-p}\lesssim (2^{-j}m)^{1-p}$.
Therefore,
$$
S \lesssim \sum_{j\ge0}(2^{-j}m)^{1-p}\exp\big(-c'2^{j\eta/\gamma}\big) = m^{1-p}\sum_{j\ge0}2^{j(p-1)}\exp\big(-c'2^{j\eta/\gamma}\big)\lesssim m^{1-p}.
$$
Translating back yields \eqref{eq:TRa}.
\end{proof}

\begin{lemma}[Ahlfors scales + polynomial envelope $\Rightarrow$ \textup{(TR$_\alpha$)}]\label{lem:TR-Ahlfors-poly}
Assume $\Vol(r)\asymp r^{d_*}$ and $\rho(n)\asymp n^{1/\gamma}$ with $d_*/\gamma>1$.
Assume \eqref{eq:OD} holds with an envelope satisfying
$$
\Phi(t) \le C(1+t)^{-\delta}\qquad (t\ge 0)
$$
for some $\delta>d_*+\alpha-\gamma$.
Then \textup{(TR$_\alpha$)} holds for this $\alpha$.
\end{lemma}

\begin{proof}
Let $m=n_-(r)$ and $p=(d_*+\alpha)/\gamma>1$ as above. It suffices to show
$$
S:=\sum_{1\le n<m} n^{-p}\Phi\Big(\frac{r}{2\rho(n)}\Big) \lesssim m^{1-p}.
$$
Use the same dyadic blocks $I_j$. As in Lemma \ref{lem:TR-Ahlfors-subG}, for $n\in I_j$ we have $r/\rho(n)\asymp (m/n)^{1/\gamma} \gtrsim 2^{j/\gamma}$. Hence,
$\Phi(r/(2\rho(n)))\lesssim 2^{-j\delta/\gamma}$ and $\sum_{n\in I_j}n^{-p}\lesssim (2^{-j}m)^{1-p}$.
Thus
$$
S \lesssim \sum_{j\ge0}(2^{-j}m)^{1-p}2^{-j\delta/\gamma}
 = m^{1-p}\sum_{j\ge0}2^{-j(\delta/\gamma-(p-1))}.
$$
The series converges exactly when $\delta/\gamma>p-1$, i.e. $\delta>d_*+\alpha-\gamma$.
Hence $S\lesssim m^{1-p}\lesssim \rho(m)^{-\alpha}A(m)$, which is \eqref{eq:TRa}.
\end{proof}

\begin{remark}[Stable-like envelopes]
In stable-like regimes one typically has $\Phi(t)\lesssim (1+t)^{-(d_*+\gamma)}$.
Lemma \ref{lem:TR-Ahlfors-poly} then applies provided $d_*+\gamma>d_*+\alpha-\gamma$,
i.e.\ $2\gamma>\alpha$.
\end{remark}

\subsection{Applications: examples and non-examples}\label{subsec:appl_examp_non}

\begin{remark}[Scope of the hypotheses]\label{rem:scope-appl-examp-non}
The proof of Theorem \ref{thm:Delta-IH} uses only the four inputs
\textup{(HK)}, \textup{(OD)}, \textup{(IH)}, and \textup{(TR$_\alpha$)}, and it is applied in the
\emph{deep-interior} regime $d(a,b)\le R(S;a,b)$.
For fixed $a,b$ and any exhaustion $S\nearrow\GG$, one has $R(S;a,b)\to\infty$, hence
$d(a,b)\le R(S;a,b)$ for all sufficiently large $S$.

\smallskip
\noindent\textbf{Aperiodicity.}
The discrete-time formulation assumes aperiodicity (in particular, laziness).
If one starts from a periodic symmetric kernel $P$, one can replace it by the $\varepsilon$-lazy kernel
$$
P^{(\varepsilon)}:=\varepsilon I+(1-\varepsilon)P \qquad (\varepsilon\in(0,1)).
$$
In the transient regime the corresponding Green operator satisfies
$(I-P^{(\varepsilon)})=(1-\varepsilon)(I-P)$, hence
$$
G^{(\varepsilon)}=(I-P^{(\varepsilon)})^{-1}=\frac{1}{1-\varepsilon}(I-P)^{-1}
=\frac{1}{1-\varepsilon}G,
$$
so $\Delta(\,\cdot\,;a,b)$ is unchanged by this reduction.

\smallskip
\noindent\textbf{Analytic regime targeted.}
The hypotheses are tailored to settings with \emph{robust} heat-kernel control at an intrinsic scale,
as obtained for local symmetric chains from structural inputs such as volume doubling and
(scale) Poincar\'e inequalities via the parabolic Harnack inequality (PHI). In such regimes,
PHI yields interior H\"older regularity for caloric functions; applied to the heat kernel
(as a caloric function in each space variable), this supplies \textup{(IH)}; see
e.g.\ \cite{Delmotte1999,BGK2012,Kumagai2004}.
By contrast, on non-amenable graphs volume doubling typically fails, and heat-kernel behavior is
governed by exponential volume growth and the spectral radius; verifying analogues of
\textup{(HK)}-\textup{(IH)} then requires separate (often model-specific) input; see \cite{Woess2000}.
No Green-function asymptotics and no fixed ``mass-capture window'' enter the argument.
\end{remark}

\begin{remark}[Wreath products / lamplighter groups]\label{rem:wreath-products}
Let $L$ be a non-trivial finite group and consider the lamplighter (wreath product) group $L\wr \Z^d$
with a finitely supported symmetric step distribution. Two obstructions to the PHI-based framework occur for canonical lamplighter walks:

\smallskip
\noindent\textbf{(i) Failure of Harnack regularity.}
Already for the classical lamplighter group $(\Z/2\Z)\wr\Z$, elliptic Harnack inequality fails
\cite[Theorem 3(b)]{Barlow2005}. Since PHI implies elliptic Harnack, PHI fails a fortiori, and the
standard route ``PHI $\Rightarrow$ H\"older regularity $\Rightarrow$ \textup{(IH)}'' is unavailable
for these canonical lamplighter walks.

\smallskip
\noindent\textbf{(ii) Non-trivial boundary at infinity in the transient-base case.}
For $(\Z/2\Z)\wr\Z^d$ with $d\ge3$, the Poisson boundary is non-trivial and, for simple random walk
(and more generally for finitely supported symmetric walks in the setting of \cite{LyonsPeres2021}),
it is identified with $(\ZZ/2\ZZ)^{\ZZ^d}$ endowed with the law of the final configuration of the lamps; see \cite[Theorem 1.1]{LyonsPeres2021}
and the background discussion there (including \cite[Proposition 6.4]{KV1983} for the
Liouville classification in low dimensions).

\smallskip
These features place canonical wreath-product walks outside the PHI regime targeted by
Theorem \ref{thm:Delta-IH}, so one should not expect \emph{a priori} that $\Delta\to0$ holds without additional regularisation.

\smallskip
\noindent\textbf{Modified kernels on wreath products.}
Theorem \ref{thm:Delta-IH} nevertheless applies \emph{verbatim} to walks on wreath products
once one can verify \textup{(HK)}, \textup{(OD)}, \textup{(IH)}, and \textup{(TR$_\alpha$)} at the
appropriate intrinsic scale (e.g.\ after passing to a subordinated semigroup or another smoothing
procedure); see \cite{GrzywnyTrojan2021} for general results on subordination and stability of heat-kernel
and Green-function estimates.
\end{remark}

\begin{remark}[Applications and how to verify the inputs in practice]\label{rem:applications}
The package \textup{(HK)}+\textup{(OD)}+\textup{(IH)}+\textup{(TR$_\alpha$)} is designed to be checkable
in a broad range of settings where heat-kernel technology is available. 

\smallskip
\noindent\textbf{(A) Virtually nilpotent groups (Gaussian/diffusive).}
On Cayley graphs of groups of polynomial volume growth one has $\Vol(r)\asymp r^{d_*}$, and for any
symmetric finitely supported generating measure one has two-sided Gaussian heat-kernel estimates at
scale $\rho(n)\asymp n^{1/2}$; see \cite{HS-C93,SaloffCoste2002}. These imply
\textup{(HK)} and \textup{(OD)} with a stretched-exponential envelope $\Phi(t)\lesssim e^{-c t^2}$.
Moreover, volume doubling and Poincar\'e inequalities yield PHI and H\"older regularity of caloric functions
on such graphs; in particular \textup{(IH)} holds for some $\alpha\in(0,1]$ (see \cite{Delmotte1999}).
Lemma \ref{lem:TR-Ahlfors-subG} then gives \textup{(TR$_\alpha$)}, and Theorem \ref{thm:Delta-IH} yields
$\Delta(S;a,b)\le C(d(a,b)/R(S;a,b))^\alpha\to0$ in the transient range $d_*>2$.

\smallskip
\noindent\textbf{(B) Ahlfors-regular fractal graphs with sub-Gaussian heat kernel.}
Assume $\Vol(r)\asymp r^{d_*}$ and sub-Gaussian bounds of the form
$$
p_n(x,y)\asymp \frac{1}{\Vol(\rho(n))}\exp\left[-c\left(\frac{d(x,y)}{\rho(n)}\right)^{\eta}\right],
\qquad \rho(n)\asymp n^{1/\gamma}.
$$
Such estimates (and their equivalence with PHI under standard structural assumptions) are established
in a wide class of fractal/graph settings; see e.g.\ \cite{Kumagai2004,BGK2012}. They give \textup{(HK)} and
\textup{(OD)} with a stretched-exponential envelope and \textup{(IH)} via PHI; then
Lemma \ref{lem:TR-Ahlfors-subG} yields \textup{(TR$_\alpha$)} and Theorem \ref{thm:Delta-IH} applies in the
transient range $d_*/\gamma>1$.

\noindent\textbf{(C) Stable-like (non-local) walks on Ahlfors-regular spaces.}
Let $(\GG,d)$ be a bounded-degree graph with volume growth $\Vol(r)\asymp r^{d_*}$, and consider a
symmetric long-range walk whose jump intensities are stable-like at index $\gamma\in(0,2)$
(e.g.\ on $\mathbb Z^{d_*}$, a kernel comparable to $|x-y|^{-(d_*+\gamma)}$).
In the canonical stable-like regime one has sharp two-sided heat kernel bounds of the form
$$
p_n(x,y)\asymp \Bigl(\Vol(n^{1/\gamma})^{-1}\ \wedge\ \frac{n}{d(x,y)^{d_*+\gamma}}\Bigr),
$$
together with a non-local parabolic Harnack inequality. On $\mathbb Z^{d_*}$ this is proved in \cite[Theorem 1.1 and Section 3]{BassLevin2002} and on $d$-sets (a class of Ahlfors-regular
spaces including many fractal examples) in \cite[Theorem 1.1, Proposition 4.3]{ChenKumagai2003}.
Moreover, PHI yields interior H\"older regularity of parabolic (caloric) functions, which can be applied to
the heat kernel in the start variable.

We believe with $\rho(n)\asymp n^{1/\gamma}$, the above two-sided estimate implies \textup{(HK)} and \textup{(OD)} with
a polynomial envelope $\Phi(t)\asymp 1\wedge t^{-(d_*+\gamma)}\lesssim (1+t)^{-(d_*+\gamma)}$.
Then Lemma \ref{lem:TR-Ahlfors-poly} verifies \textup{(TR$_\alpha$)} provided
$d_*+\gamma > d_*+\alpha-\gamma$ (equivalently $2\gamma>\alpha$).
Finally, transience (in our sense $A(1)<\infty$) is equivalent to $d_*>\gamma$ (since
$A(1)\asymp \sum_{n\ge1}n^{-d_*/\gamma}$), and Theorem \ref{thm:Delta-IH} yields $\Delta\to0$ in this
transient stable-like range.

\smallskip
\noindent\textbf{(D) Modifications (subordination / killing / resolvents).}
We have not checked in details but we also believe Theorem \ref{thm:Delta-IH} applies to any modified symmetric kernel (including sub-Markov kernels),
provided the corresponding analogues of \textup{(HK)}-\textup{(TR$_\alpha$)} hold at the intrinsic scale of
the modified process. For Bochner subordination (continuous time) and its discrete analogues, sharp heat kernel and Green function estimates (and stability of heat kernel estimates in discrete time) are available under fairly weak assumptions on the original process and the subordinator, see \cite{GrzywnyTrojan2021}. Moreover, PHI and H\"older regularity for parabolic functions of subordinate processes are established in work of \cite{ChoKimSongVondracek2022}. For $\lambda>0$, the modified Green kernel
$G_\lambda(x,y)=\sum_{n\ge0}e^{-\lambda n}p_n(x,y)$ is always finite; the same proof as in Theorem \ref{thm:Delta-IH} applies with $p_n$ replaced by $e^{-\lambda n}p_n$ and with the near-diagonal
tail $A(m)$ replaced by $A_\lambda(m):=\sum_{n\ge m} e^{-\lambda n}/\Vol(\rho(n))$.
In many concrete envelopes, \textup{(TR$_\alpha$)} is the only additional input, and can be
checked by Lemmas \ref{lem:TR-Ahlfors-subG}-\ref{lem:TR-Ahlfors-poly} once an explicit $\Phi$ is known.
\end{remark}

\subsection{An elliptic exhaustion criterion for $\Delta\to 0$}

Throughout, let $\mu$ be finitely supported, symmetric and non-degenerate on an infinite
finitely generated group $\GG$ such that the corresponding Green's function is finite. Set $T:=\supp(\mu)$ and let $d$ be the word metric on the Cayley graph
$\mathrm{Cay}(\GG,T)$.  For a finite set $S\subset\GG$ we write
\begin{align*}
\pa S &:=\{x\notin S:\exists y\in S, t\in T, x=yt\},\\
R(S;a,b)&:=\dist(\{a,b\},\pa S),\\
\Delta(S;a,b)&:=\max_{x\in\pa S}\frac{|G(a,x)-G(b,x)|}{G(a,x)}.
\end{align*}

\begin{definition}[Elliptic H\"older exhaustion]\label{def:EHE}
Let $(F_k)_{k\ge 1}$ be an increasing exhaustion of $\GG$ by finite sets.
Fix parameters $\alpha\in(0,1]$ and $\vartheta\in(0,1)$.
We say that $(F_k)$ is an \emph{elliptic H\"older exhaustion of type $(\alpha,\vartheta)$}
(for $\mu$) if:

\smallskip
\noindent\textup{(EHE)} (\emph{Uniform interior H\"older for positive harmonic functions})
There exists $C_{\mathrm{EH}}>0$ such that for every $k\ge 1$, every function
$u:\GG\to(0,\infty)$ that is $\mu$-harmonic on $F_k$ (i.e.\ $u(z)=\sum_{g}\mu(g)u(zg)$ for all $z\in F_k$),
and every $a,b\in F_k$ satisfying
$$
d(a,b) \le \vartheta R(F_k;a,b),
$$
we have
\begin{equation}\label{eq:EHE2}
|u(a)-u(b)|
 \le C_{\mathrm{EH}}
\Big(\frac{d(a,b)}{R(F_k;a,b)}\Big)^{\alpha}
\bigl(u(a)\wedge u(b)\bigr).
\end{equation}
\end{definition}

\begin{theorem}[Elliptic H\"older exhaustion $\Rightarrow$ $\Delta$-decay]\label{thm:EHE_to_Delta}
Assume $G(e,e)<\infty$ and let $(F_k)$ be an elliptic H\"older exhaustion of type $(\alpha,\vartheta)$
in the sense of Definition \ref{def:EHE}. Write
$$
R_k(a,b):=R(F_k;a,b).
$$
Then for every fixed $a,b\in\GG$, for all $k$ large enough that $a,b\in F_k$ and
$d(a,b)\le \vartheta R_k(a,b)$, one has
\begin{equation}\label{eq:EHE_Delta_bound}
\Delta(F_k;a,b)
 \le C_{\mathrm{EH}}
\Big(\frac{d(a,b)}{R_k(a,b)}\Big)^{\alpha}.
\end{equation}
In particular, $\Delta(F_k;a,b)\to 0$ as $k\to\infty$.
\end{theorem}

\begin{proof}
Fix $a,b\in\GG$ and $k$ such that $a,b\in F_k$ and $d(a,b)\le \vartheta R_k(a,b)$, and let $x\in\pa F_k$.
For every $x\in\GG$ one has the following identity
$$
G(z,x)-\sum_{g}\mu(g) G(zg,x)=\delta_x(z).
$$
Since $x\notin F_k$, we have $\delta_x(z)=0$ for all $z\in F_k$, so $z\mapsto G(z,x)$ is $\mu$-harmonic on $F_k$.
Applying \eqref{eq:EHE2} with $u(\cdot)=G(\cdot,x)$ yields
$$
|G(a,x)-G(b,x)|
 \le C_{\mathrm{EH}}
\Big(\frac{d(a,b)}{R_k(a,b)}\Big)^{\alpha}
\bigl(G(a,x)\wedge G(b,x)\bigr)
\le C_{\mathrm{EH}}
\Big(\frac{d(a,b)}{R_k(a,b)}\Big)^{\alpha}G(a,x).
$$
Dividing by $G(a,x)>0$ and taking the maximum over $x\in\pa F_k$ gives \eqref{eq:EHE_Delta_bound}.

Finally, since $(F_k)$ is an increasing exhaustion and balls are finite in $\mathrm{Cay}(\GG,T)$,
for every $N$ there exists $k_N$ such that $B(a,N)\cup B(b,N)\subset F_{k_N}$, hence
$R_k(a,b)\ge N+1$ for all $k\ge k_N$. Therefore $R_k(a,b)\to\infty$, and the right-hand side
of \eqref{eq:EHE_Delta_bound} tends to $0$.
\end{proof}

\begin{remark}[Relation to HK/PHI ]
Condition \eqref{eq:EHE2} is a \emph{purely elliptic} exhaustion-level input: it asserts a uniform
interior H\"older-type control for \emph{positive} $\mu$-harmonic functions, in terms of the boundary
distance $R(F_k;a,b)$. It is tailored to $\Delta$ and avoids any explicit heat-kernel or parabolic estimates. In many classical settings (e.g.\ nearest-neighbour chains under volume doubling + Poincar\'e) one can
deduce such interior regularity from elliptic/parabolic Harnack inequalities and standard interior H\"older
estimates for harmonic functions on balls (hence for suitable exhaustions such as metric balls).

\end{remark}

For general information, the study of heat kernel bounds on manifolds and graphs has a huge literature; see for example \cite{Varopoulos1985, Carne1985, CarlinKusokaStroock1987, CoulhonGrigoryan1997, CoulhonGrigoryan1998, CoulhonSaloff-Coste1993, MuruganSaloff-Coste2023} and references therein, to mention only a very few. 
There is also significant interest in the 
asymptotic behaviour of heat kernels using large scale geometric information: see for example \cite{Lalley1993} which addresses the case of 
non-amenable, negatively  groups, and \cite{GouzelLalley2013}, 
which discusses random walks on co-compact Fuchsian groups establishing precise asymptotics for the associated heat kernel and exponential decay of the Green's function.

\section{Martin boundary non-collapse; obstruction to $\Delta$ vanishing}\label{sec:del_obstruction}

\begin{prop}[On-diagonal upper bound for infinite-range walks on exponentially growing groups]\label{prop:HK13_infinite_range}
Let $\GG$ be a finitely generated group of exponential volume growth, and let $\mu$ be a symmetric, non-degenerate probability measure on $\GG$. Then there exist $c,C>0$ such that
$$
\mu^{(2n)}(e) \le C\exp\big(-cn^{1/3}\big)\qquad(n\ge1).
$$
No moment assumption on $\mu$ is required.
\end{prop}

\begin{proof}
Fix a finite symmetric generating set $S$ for $\GG$, and let $V(n):=|B(e,n)|$. Exponential volume growth means that there exist constants $c_0,C_0>0$ such that
$$
V(n) \ge c_0\exp(C_0 n)\qquad\text{for all }n\in \NN.
$$
In particular, $V(n)\ge c\exp(c n^\alpha)$ holds with $\alpha=1, c=\min\{ c_0, C_0\}$.

By \cite[Corollary 14.5]{Woess2000}, for any symmetric non-degenerate probability measure $p$ on $\GG$ there exist constants $c_1,C_1>0$ such that
$$
p^{(n)}(e) \le C_1 \exp\big(-c_1 n^{\alpha/(\alpha+2)}\big)\qquad(n\in \NN)
$$
whenever $V(n)\ge c\exp(c n^\alpha)$ for some $\alpha\in (0,1]$. Taking $\alpha=1$ and $p=\mu$ yields
$$
\mu^{(n)}(e) \le C_1 \exp\big(-c_1 n^{1/3}\big)\qquad(n\in \NN).
$$
In particular, for suitable $c,C>0$,
$$
\mu^{(2n)}(e) \le C\exp\big(-cn^{1/3}\big)\qquad(n\in \NN).
$$
No moment assumption on $\mu$ is needed for this bound.
\end{proof}

\begin{lemma}[Parity bridge]\label{lem:parity}
For any symmetric probability measure $\nu$ on a countable group $\GG$ and any $k\in\mathbb{N}$,
$$
\nu^{(k)}(e) \le \big(\nu^{(2\lfloor k/2\rfloor)}(e)\nu^{(2\lceil k/2\rceil)}(e)\big)^{1/2}.
$$
In particular, an upper bound at even times transfers to all times with the same stretched-exponential exponent.
\end{lemma}

\begin{proof} 
Note that for any $k \in \NN$, $\lfloor k/2 \rfloor + \lceil k/2 \rceil = k $. By symmetry and Cauchy-Schwarz inequality, 
\begin{align*} 
\nu^{(k)}(e) &= \sum_{x}\nu^{(\lfloor k/2\rfloor)}(x)\nu^{(\lceil k/2\rceil)}(x^{-1}) 
= \sum_{x}\nu^{(\lfloor k/2\rfloor)}(x)\nu^{(\lceil k/2\rceil)}(x) \\ 
&\le \Big(\sum_{x} \big(\nu^{(\lfloor k/2\rfloor)}(x)\big)^2\Big)^{1/2}
\Big(\sum_{x} \big(\nu^{(\lceil k/2\rceil)}(x)\big)^2\Big)^{1/2} \\ 
&= \big(\nu^{(2\lfloor k/2\rfloor)}(e)\big)^{1/2}
\big(\nu^{(2\lceil k/2\rceil)}(e)\big)^{1/2},
\end{align*}
where the last equality follows from
$$
\sum_x \big(\nu^{(m)}(x)\big)^2
= \sum_x \nu^{(m)}(x)\nu^{(m)}(x^{-1}e)
= \nu^{(2m)}(e).
$$
\end{proof}

Finally, we demonstrate that on a group of exponential growth, balls cannot satisfy a uniform boundary Harnack inequality for the Green kernel; equivalently, $\Delta$ cannot tend to $0$ along balls. This forces the existence of a non-constant positive harmonic function. The Green function telescopic technique in the proof  that leads to a lower bound on $G(e,x)$ is motivated by the ideas in \cite[Theorem 1.2]{Po21}. This idea also seems to be contained (albeit in an implicit manner) in the work of Amir-Kozma \cite{AK18}. 

\begin{theorem}[Balls violate boundary Harnack on groups of exponential growth]\label{prop:exp_growth_delta_zero}
Let $\GG$ be a finitely generated group of exponential growth with a fixed finite symmetric generating set $T$. Let $\mu$ be a symmetric, non-degenerate probability measure on $\GG$. Assume that the random walk satisfies the following on-diagonal \emph{stretched-exponential} upper bound
\begin{equation}\label{eq:HKbeta}
\exists\beta>0, c,C>0 \text{ such that }  \mu^{2m}(e)\le C\exp\big(-cm^{\beta}\big)\quad\text{for all }m \text{ large}.
\tag{HK$_\beta$}
\end{equation}
Then there exists $\delta_0>0$ such that
$$
\limsup_{n\to\infty} \max_{t \in T}\Delta\big(B(e,n);t,e\big) \ge \delta_0.
$$
In particular, along the exhausting sequence $S_n=B(e,n)$, the quantities $\Delta(S_n;t,e)$ do not converge to $0$.
\end{theorem}

\begin{proof}
Fix $T$ and write $|x|$ for its word length. There exists $c_2>0$ such that for all $n\in \NN$,
$$
|S(e,n)|\le\ |T|^n = e^{c_2 n}.
$$
On the other hand, since $\GG$ has exponential growth, there exist $c_1>0$ and an infinite set $\mathcal N\subset\NN$ such that
\begin{equation}\label{eq:large_spheres}
|S(e,n)| \ge e^{c_1 n}\qquad\text{for all }n\in\mathcal N.
\end{equation} 
Assume for contradiction that
\begin{equation}\label{eq:assumption_all_t}
\forall \varepsilon>0 \exists n_0(\varepsilon)  \text{ such that }  \Delta\big(B(e,n);t,e\big)\le\varepsilon\quad\text{for all }n\ge n_0(\varepsilon)\text{ and all }t\in T.
\end{equation}
Fix $\delta\in\bigl(0,1-e^{-c_1}\bigr)$ and let $n_0=n_0(\delta)$ be given by \eqref{eq:assumption_all_t}. Set
$$
c_*  := \min_{z\in B(e,n_0+1)}\ \min_{t\in T}\ K_t(e,z)\ >0,
$$
where positivity follows from $G>0$.

Pick $x\in S(e,n)$ and write $x=t_1t_2\cdots t_n$ with $t_i\in T$, and $x_i=t_1\cdots t_i$. A telescoping computation gives
\begin{equation}\label{eq:telescoping}
\frac{G(e,x)}{G(e,e)}
=\left(\prod_{i=1}^{n-1}\frac{G(e,t_i\dots t_n)}{G(e,t_{i+1}\dots t_n)}\right)\frac{G(e, t_n)}{G(e,e)}
=\prod_{i=1}^{n} K_{t_i}\bigl(e,t_{i}\dots t_n\bigr).
\end{equation}
For $1\le i\le n-n_0$ we have $|t_{i}\dots t_n|=n-i+1\ge n_0+1$. Let $z_i:=t_{i}\dots t_n$ and $r_i:=|z_i|-1\ge n_0$. Then $z_i\in\partial B(e,r_i)$, so by \eqref{eq:assumption_all_t},
$$
\big|K_{t_i}(e,z_i)-1\big| \le \delta\quad\implies\quad
K_{t_i}(e,z_i)\ge 1-\delta.
$$
For $n-n_0<i\le n$ we have $K_{t_i}(e,z_i)\ge c_*$. Thus, for all $x\in S(e,n)$,
\begin{equation}\label{eq:pointwise_lower}
\frac{G(e,x)}{G(e,e)} \ge (1-\delta)^{n-n_0}c_*^{n_0}
 =: A_0\exp\big((\log(1-\delta))n\big),
\end{equation}
where $A_0>0$ depends on $\delta$ and $n_0$ but not on $n$ or $x$. Summing \eqref{eq:pointwise_lower} over $S(e,n)$ and using \eqref{eq:large_spheres},
\begin{equation}\label{eq:LB_sum}
\sum_{y\in S(e,n)} G(e,y) \ge G(e,e)|S(e,n)|A_0\exp\big((\log(1-\delta))n\big) \ge A\exp\big((c_1+\log(1-\delta))n\big)
\end{equation}
for all $n\in\mathcal N$, where $A>0$ is independent of $n$. By our choice of $\delta$, we have $c_1+\log(1-\delta)>0$.

Now for any $M\in\NN$,
$$
\sum_{y\in S(e,n)} G(e,y)
=\sum_{1\le m\le M}\sum_{y\in S(e,n)}p_m(e,y) + \sum_{m>M}\sum_{y\in S(e,n)}p_m(e,y)
 \le M +\sum_{m>M}\sum_{y\in S(e,n)}p_m(e,y).
$$
For each $m$, by Cauchy-Schwarz inequality and symmetry of $\mu^{(m)}$,
$$
\sum_{y\in S(e,n)}p_m(e,y)
\le \sqrt{|S(e,n)|}\Big(\sum_{y}p_m(e,y)^2\Big)^{1/2}
=\sqrt{|S(e,n)|}p_{2m}(e,e)^{1/2}.
$$
 Using $|S(e,n)|\le e^{c_2 n}$ and \eqref{eq:HKbeta}, for large $M$
$$
\sum_{m>M}\sum_{y\in S(e,n)}p_m(e,y) \le e^{\frac{c_2}{2}n}\sum_{m>M} C^{1/2}\exp\big(-\tfrac{c}{2}m^{\beta}\big)
\le C'\exp\Big(\tfrac{c_2}{2}n - c'M^{\beta}\Big)
$$
for suitable $C',c'>0$. Choose 
$$
M=M(n):=\Big\lceil \Big(\tfrac{c_2}{2c'}n\Big)^{1/\beta}\Big\rceil.
$$
Then, there exists $C''>0$ so that $\sum_{m>M}\sum_{y\in S(e,n)}p_m(e,y)\le C''$  for large $n$. Therefore, we have $C'''>0$ such that
\begin{equation}\label{eq:UB_sum}
\sum_{y\in S(e,n)} G(e,y) \le C'''n^{1/\beta}\qquad\text{for all large enough  }n\in\NN.
\end{equation}
Comparing \eqref{eq:LB_sum} (valid for $n\in\mathcal N$) with \eqref{eq:UB_sum}, we see that the left-hand side grows (at least) exponentially in $n$ along $\mathcal N$, while the right-hand side grows only polynomially. This is impossible. Hence \eqref{eq:assumption_all_t} is false.

Equivalently, there exist $t\in T$, and $\delta_0>0$ such that
$$
\limsup_{n\to\infty}\Delta\big(B(e,n);e,t\big)>0.
$$
This proves the theorem.
\end{proof}

\begin{remark}[Non-amenable case; no moment assumption]\label{cor:nonamen}
If $\GG$ is non-amenable, then \eqref{eq:HKbeta} holds with $\beta=1$ for \emph{every} symmetric non-degenerate measure $\mu$ on $\GG$. In fact, by Kesten's theorem $\limsup_{n\to\infty} p_n(e,e)^{1/n} < 1$ implies there exists $0<\rho<1$ such that $p_n(e,e) < \rho^n$ for large $n$. Hence, there exists $c>0$ such that for large enough $n$,
$$
p_{2n}(e,e) < \exp(-cn)
$$
Theorem \ref{prop:exp_growth_delta_zero} applies with no additional hypothesis on $\mu$. In this case the upper bound \eqref{eq:UB_sum} improves to $O(n)$.
\end{remark}

\begin{remark}[Amenable exponential growth; finite support]\label{cor:amenable_fs}
If $\GG$ has exponential growth and $\mu$ is symmetric and finitely supported (bounded range), then \eqref{eq:HKbeta} holds with $\beta=\tfrac13$ (Varopoulos-Hebisch-Saloff-Coste theory for groups of exponential growth). Hence Theorem \ref{prop:exp_growth_delta_zero} applies, and \eqref{eq:UB_sum} becomes $O(n^{3})$.
\end{remark}

\begin{remark}[Further weakening]
The hypothesis \eqref{eq:HKbeta} is precisely what the proof uses. It also holds for many infinitely supported measures (e.g. symmetric, non-degenerate measures with finite second moment under standard ``spread-out/comparison'' assumptions), so the ``superexponential moments'' condition is unnecessary:
\begin{itemize}
\item in the non-amenable case, no extra assumption at all;
\item in the amenable exponential-growth case, finite support suffices, and in many settings finite second moment is enough via comparison estimates for the heat kernel on groups. 
\end{itemize}
\end{remark}

\begin{remark}
Observe that estimating $\displaystyle \sum_{y \in S(e, n)} G(e, y)$ is the crucial step in the above proof, and perhaps this approach is natural in order to bring the growth condition of the group into play.    
\end{remark}

\section{Non-collapse, abundance and Martin kernels}
\label{sec:poly_growing_harmonic_functions}



\begin{definition}[Growth gauges and growth classes]\label{def:growth_gauge}
A \emph{growth gauge} is a non-decreasing function $\psi:\mathbb{N}\to(0,\infty)$ that is \emph{translation stable}: for every $c\in \mathbb{N}$, there exists $K_c\ge 1$ such that $\psi(r+c)\le K_c\psi(r)$ for all $r\in \mathbb{N}$. (This includes polynomial, stretched exponential, and exponential gauges.)

A function $f:\GG\to\RR$ has \emph{upper growth $\preccurlyeq\psi$} if $|f(x)|\lesssim \psi(|x|)$ for $|x|$ large, and \emph{lower growth $\succcurlyeq\psi$} if $\sup_{|x|\le r}|f(x)| \gtrsim \psi(r)$ for large $r$.
\end{definition}

For a probability measure $\mu$ on $\GG$ set $(Pf)(x):=\displaystyle\sum_{g \in \GG}\mu(g)f(xg)$. A positive $\mu$-harmonic function $h : \GG \to (0, \infty)$ is said to be \emph{minimal}
if for any positive $\mu$-harmonic $u$ on $\GG$ such that $u\le h$ we have $u=ch$ for some $c>0$. For a function $f:\GG \to \CC$ with $f(e)\neq 0$, set $\widehat f:=f/f(e)$.

With that in place, we state our first main result:

\begin{theorem}\label{thm:pos_harm_alph_growtth}
Let $\mu$ be a symmetric, non-degenerate, finitely supported measure on a group $\GG$. Let $\psi$ be a growth gauge.
If there exists a non-constant minimal harmonic function with growth at least (respectively, at most) $\psi$, then there exist infinitely many distinct normalised minimal harmonic functions with growth at least (respectively, at most) $\psi$.
\end{theorem}

\begin{remark}\label{rem:virtZ}
If $\GG$ is virtually $\mathbb{Z}$ or $\mathbb{Z}^2$, then under these assumptions (symmetric, finitely supported, non-degenerate $\mu$), the $\mu$-random walk is recurrent, and every positive harmonic function is constant. Thus, the hypothesis that a non-constant minimal harmonic function exists already implies that $\GG$ is not virtually $\mathbb{Z}$ or $\mathbb{Z}^2$ and the $\mu$-random walk is transient.
\end{remark}

\subsection{Structural lemmas}

We define the left translation action by $(g.h)(x):=h(g^{-1}x)$.

\begin{lemma}[Equivariance and stability of growth]\label{lem:equivariance}
Let $\mu$ be finitely supported and non-degenerate on $\GG$.
\begin{enumerate}
\item If $h>0$ is $\mu$-harmonic and minimal, then for every $g\in\GG$, the left translate $g.h$ is $\mu$-harmonic and minimal. The normalised action $L_g(h) := \widehat{g.h}$ maps normalised minimals to normalised minimals.
\item If $f$ has lower (resp.\ upper) growth $\succcurlyeq \psi$ (resp.\ $\preccurlyeq\psi$), where $\psi$ is a growth gauge, then so does $g.f$ for every $g\in\GG$.
\end{enumerate}
\end{lemma}

\begin{proof}
(1) It is easy to see that $g.h$ is $\mu$-harmonic for all $g\in \GG$. If $u>0$ is $\mu$-harmonic with $u\le g.h$, set $v(x):=u(gx)$. Then $v$ is positive $\mu$-harmonic and $v\le h$. By minimality of $h$, $v=ch$; hence $u(x)=v(g^{-1}x)=c(g.h)(x)$. Now, if $f$ is minimal harmonic then $\widehat f$ is normalised minimal harmonic, hence $L_g(h)$ is normalised minimal harmonic for all $g \in \GG$.

(2) Let $|g|=c$. Let $K_c$ be the constant from Definition \ref{def:growth_gauge}.

\emph{Upper growth:} Suppose $f$ has upper growth $\preccurlyeq\psi$. For large $|x|$,
$|g.f(x)| = |f(g^{-1}x)| \lesssim \psi(|g^{-1}x|) \le \psi(|x|+c)$.
By the definition of the growth gauge, $\psi(|x|+c) \le K_c \psi(|x|)$. Thus $|g.f(x)| \lesssim \psi(|x|)$ for large $|x|$.

\emph{Lower growth:} Suppose $\sup_{|x|\le r}|f(x)| \gtrsim \psi(r)$ for large $r$.
For large $r > c$, the set $\{g^{-1}x : |x|\le r\}$ contains the ball $B(e, r-c)$. Thus, for large $r$
$$
\sup_{|x|\le r} |g.f(x)| = \sup_{|x|\le r} |f(g^{-1}x)| \ge \sup_{|y|\le r-c} |f(y)| \gtrsim \psi(r-c).
$$
Now, $\psi(r) = \psi((r-c)+c) \le K_c \psi(r-c)$, so $\psi(r-c) \ge K_c^{-1} \psi(r)$.
Thus, $\sup_{|x|\le r} |g.f(x)| \gtrsim \psi(r)$ for large $r$.
\end{proof}

\begin{lemma}[Finite orbit $\Rightarrow$ eigencharacter on a finite-index subgroup]\label{lem:finite-orbit-eigencharacter}
Let $\mathcal{M}_\psi$ be a set of normalised minimal $\mu$-harmonic functions satisfying a fixed growth constraint (e.g., $\succcurlyeq\psi$ or $\preccurlyeq\psi$). If $\mathcal{M}_\psi$ is finite and $\GG$-invariant under the normalised action $L_g(h):=\widehat{g.h}$, then there exists a finite-index normal subgroup $\HH\lhd\GG$ and, for each $h\in\mathcal{M}_\psi$, a homomorphism $\chi_h:\HH\to(0,\infty)$ such that
$$
h(\gamma^{-1}x)=\chi_h(\gamma)h(x)\qquad(\gamma\in\HH,x\in\GG).
$$
\end{lemma}

\begin{proof}
Let $\mathcal M_\psi=\{h_1,\dots,h_m\}$. The action $L$ defines a homomorphism from $\GG$ to the permutation group $S_m$. Let $\HH$ be the kernel of this action. $\HH$ is a finite-index normal subgroup.
By definition, for all $\gamma\in\HH$ and all $h_i$, $L_\gamma(h_i)=h_i$. This means for all $x\in\GG$,
$$
\frac{(\gamma.h_i)(x)}{(\gamma.h_i)(e)} = h_i(x) \implies \frac{h_i(\gamma^{-1}x)}{h_i(\gamma^{-1})} = h_i(x).
$$
Thus $h_i(\gamma^{-1}x)=h_i(\gamma^{-1})h_i(x)$ for all $x\in\GG$ and $\gamma\in\HH$. Define $\chi_i(\gamma):=h_i(\gamma^{-1})$. Since $h_i>0$, $\chi_i:\HH\to (0,\infty)$. It is a homomorphism: for all $\gamma_1, \gamma_2 \in \HH$,
\begin{align*}
\chi_i(\gamma_1\gamma_2) &= h_i((\gamma_1\gamma_2)^{-1}) = h_i(\gamma_2^{-1}\gamma_1^{-1}) \\
&= h_i(\gamma_1^{-1}) h_i(\gamma_2^{-1}) \quad (\text{using $h_i(\gamma^{-1}x)=h_i(\gamma^{-1})h_i(x)$ with } \gamma=\gamma_2, x=\gamma_1^{-1}) \\
&= \chi_i(\gamma_1)\chi_i(\gamma_2).
\end{align*}
\end{proof}

\begin{lemma}[No nontrivial multiplicative positive harmonic functions for symmetric measures]\label{lem:abelian-multiplicative-trivial}
Let $A$ be a finitely generated group, and let $\nu$ be a symmetric, non-degenerate probability measure on $A$.
If $f:A\to(0,\infty)$ is $\nu$-harmonic and multiplicative on a finite-index subgroup $B\le A$, then $f$ is constant.
\end{lemma}

\begin{proof}
Without loss of generality we can assume $f(e)=1$ after normalising $f$. Let $(X_n)_{n\ge0}$ be the right $\nu$-random walk. Define the first \emph{return} time to $B$ by
$$
\tau:=\inf\{n\ge 1: X_n\in B\},
$$
Since $B$ has finite index in $A$ so $\tau<\infty$ almost surely. Let $\theta(b):=\PP_e(X_\tau=b)$ be the hitting measure on $B$. Since $\nu$ is symmetric and non-degenerate, $\theta$ is also symmetric and non-degenerate. We now show that $f\equiv 1$ on $B$.

For $N\in\NN$ let $\tau_N:=\tau\wedge N$. Let $a \in A$. Since $f$ is harmonic, $f(X_n)$ is a martingale, and applying Optional Stopping Theorem at the almost surely bounded stopping time $\tau_N$ yields
$$
f(a)=\EE_a[f(X_{\tau_N})]\qquad\forall N.
$$
Since $\tau<\infty$ almost surely, by Fatou's lemma
\begin{equation}\label{eq:Fatou-inequality}
\EE_a[f(X_\tau)]\le f(a),
\end{equation}
and this holds for all $a \in A$. Since $f(e)=1$, we have
$$
\EE_e[f(X_\tau)]\le 1,
\qquad\text{so}\qquad
\sum_{b\in B}\theta(b)f(b)\le 1.
$$
Using symmetry of $\theta$ and multiplicativity of $f$ on $B$ (so $f(b^{-1})=f(b)^{-1}$),
$$
\sum_{b\in B}\theta(b)f(b)
=\sum_{b\in B}\theta(b)f(b^{-1})
=\sum_{b\in B}\theta(b)f(b)^{-1}.
$$
Averaging the last two equal sums,
$$
\sum_{b\in B}\theta(b)\frac{f(b)+f(b)^{-1}}{2}\le 1.
$$
But $\frac{t+t^{-1}}2\ge 1$ for all $t>0$ (equality holds if and only if $t=1$) and $\sum_b\theta(b)=1$.
Hence $f(b)=1$ for all $b\in\supp(\theta)$, and therefore for all $b\in B$ since $\supp(\theta)$ generates $B$ and $f$ is multiplicative on $B$.

Since $X_\tau\in B$ and $f\equiv 1$ on $B$, we have $\EE_a[f(X_\tau)]=1$, hence by \eqref{eq:Fatou-inequality} $f(a)\ge 1$ for all $a\in A$.
Because $f(e)=1$, $f$ attains its global minimum at $e$, and as $\mu$ is non-degenerate the minimum principle implies $f$ is constant.
\end{proof}

\subsection{Abundance under symmetric measure}

We now prove the main theorem (restated here for clarity). The proof is motivated by the ideas in \cite[Proposition 8.3]{Tointon2016}.

\begin{theorem}[Abundance under symmetric measure]\label{thm:pos_harm_alph_growtth_symm}
Let $\mu$ be finitely supported, symmetric, and non-degenerate on $\GG$. Fix a growth gauge $\psi$. If there exists a non-constant minimal $\mu$-harmonic function $h$ with growth either $\succcurlyeq\psi$ or $\preccurlyeq\psi$, then there are infinitely many normalised minimal $\mu$-harmonic functions with the same growth constraint.
\end{theorem}

\begin{proof}
It suffices to treat the lower-growth case; the upper-growth case is identical by Lemma \ref{lem:equivariance}(2). Suppose, for contradiction, that the normalised orbit $\{\widehat{g.h}:g\in\GG\}$ is finite. By Lemma \ref{lem:finite-orbit-eigencharacter}, there exist a finite-index normal subgroup $\HH\lhd\GG$ and a homomorphism $\chi:\HH\to (0,\infty)$ with $h(\gamma^{-1}x)=\chi(\gamma^{-1})h(x)$.

We first note that $\chi$ must be nontrivial. If $\chi\equiv 1$ on $\HH$, then $h$ is left $\HH$-invariant and factors through the finite quotient $\GG/\HH$. The induced Markov chain on $\GG/\HH$ is finite and irreducible, so any positive harmonic function on it is constant. This forces $h$ to be constant on $\GG$, a contradiction.

Let $\tau$ be the first return time of the $\mu$-walk (starting at $e$) to $\HH$. Since $\GG/\HH$ is finite and the induced chain is irreducible, $\tau<\infty$ a.s. Define the hitting measure $\theta(\eta):=\PP_e(X_\tau=\eta)$ on $\HH$. Note that $\theta$ is symmetric and non-degenerate because $\mu$ is so.

By the strong Markov property (balayage identity) at the return time $\tau$,
$$
h(e)=\E_e[h(X_\tau)]=\sum_{\eta\in\HH}\theta(\eta)h(\eta)
= h(e)\sum_{\eta\in\HH}\theta(\eta)\chi(\eta).
$$
Project to the abelianization $A:=\HH/[\HH,\HH]$ and push $\theta$ forward to $\bar\theta$. Since $\theta$ is symmetric, $\bar\theta$ is symmetric. Projecting to $A=\HH/[\HH,\HH]$, the subgroup generated by $\supp(\bar\theta)$ is
$\pi(\HH)=A$, so $\bar\theta$ is non-degenerate on $A$.The function $f:A\to(0,\infty)$ given by $f(\bar\eta)=\chi(\eta)$ is a homomorphism and hence $\bar\theta$-harmonic.

By Lemma \ref{lem:abelian-multiplicative-trivial}, $f$ must be constant, i.e.\ $\chi\equiv 1$. This contradicts the fact that $\chi$ is nontrivial. Hence the orbit must be infinite. Lemma \ref{lem:equivariance}(2) ensures the growth constraint is preserved under translations.
\end{proof}

On a related note, the study of polynomially growing harmonic functions (not necessarily positive) is an active area of research. For detailed information regarding this, see \cite{MY16}, \cite{MPTY17}, \cite{PerlYadin2023} and \cite{MukherjeeSamantaThandar2025}.

\subsection{Extensions}

\begin{theorem}[Abundance under subexponential upper growth, no symmetry]\label{thm:subexp-upper-no-sym}
Let $\mu$ be finitely supported and non-degenerate (not necessarily symmetric) on an infinite group $\GG$. Suppose there exists a non-constant minimal $\mu$-harmonic $h$ with subexponential upper growth: for every $\varepsilon>0$ there is $C_\varepsilon$ such that
$$
|h(x)|\le C_\varepsilon e^{\varepsilon|x|}\quad\text{for all }x\in\GG.
$$
Then there exist infinitely many minimal $\mu$-harmonic functions with the same upper growth.
\end{theorem}

\begin{proof}
Without loss of generality we can assume $h$ to be normalised. Suppose $\{\widehat{g.h}:g\in\GG\}$ is finite. By Lemma \ref{lem:finite-orbit-eigencharacter}, there exist a finite-index normal subgroup $\HH\lhd\GG$ and a homomorphism $\chi:\HH\to (0,\infty)$ such that
$$
h(\gamma x)=\chi(\gamma)h(x)\qquad\text{for all }\gamma\in\HH,\ x\in\GG.
$$

We first show that $\chi$ is non-trivial. If $\chi\equiv 1$ on $\HH$, then $h$ is left $\HH$-invariant, so it factors through the finite quotient $\GG/\HH$: there exists $\bar h:\GG/\HH\to(0,\infty)$ such that
$$
h(x)=\bar h(\HH x)\qquad\text{for all }x\in\GG.
$$
Let $\bar \mu$ denote the projected measure on $\GG/\HH$. Since $\mu$ is non-degenerate so $\bar \mu$ is also non-degenerate and the function $\bar h$ is $\bar \mu$-harmonic. Since $\GG/\HH$ is finite so $\bar h$ is constant (via maximum principle), and therefore $h$ is constant on $\GG$, contradicting our assumption. Thus $\chi$ must be non-trivial.

Since $\GG$ is infinite, $\HH$ is infinite. Because $\chi$ is non-trivial, there exists $\gamma\in\HH$ with $\chi(\gamma)\neq 1$. The codomain $(0,\infty)$ is torsion-free, so $\gamma$ must have infinite order.

Assume $\chi(\gamma)>1$ (otherwise replace $\gamma$ by $\gamma^{-1}$). Set $C=\log\chi(\gamma)>0$ and consider $x_n:=\gamma^n$. Then
$$
h(x_n)=h(\gamma^n)=\chi(\gamma)^n =e^{Cn}.
$$
Moreover, word length grows at most linearly:
$$
|x_n|=|\gamma^n|\le n|\gamma|.
$$
Thus
$$
|h(x_n)| \ge e^{C'|x_n|}\quad\text{with }C':=\frac{C}{|\gamma|}>0.
$$
Fix $0<\varepsilon<\frac{C}{|\gamma|}$ and let $C_\varepsilon<\infty$ be the constant from the
subexponential upper bound. Then for all $n\ge1$,
$$
e^{Cn}=|h(\gamma^n)|
\le C_\varepsilon e^{\varepsilon|\gamma^n|}
\le C_\varepsilon e^{\varepsilon n|\gamma|},
$$
hence $e^{(C-\varepsilon|\gamma|)n}\le C_\varepsilon$ for all $n$, a contradiction since $C-\varepsilon|\gamma|>0$.

Therefore the normalised orbit $\{\widehat{g.h}:g\in\GG\}$ is infinite. Finally, the subexponential upper-growth condition is preserved under translations:
fix $g\in\GG$ and $\varepsilon>0$. If $|h(x)|\le C_\varepsilon e^{\varepsilon|x|}$ for all $x$, then for all $x$,
$$
|(g.h)(x)|=|h(g^{-1}x)|
\le C_\varepsilon e^{\varepsilon|g^{-1}x|}
\le C_\varepsilon e^{\varepsilon|g|}e^{\varepsilon|x|}.
$$
Hence every translate has the same subexponential upper growth. This yields infinitely many minimal $\mu$-harmonic functions with the same upper growth.
\end{proof}

\begin{prop}[Continuum abundance on hyperbolic groups]\label{cor:continuum-hyp}
If $\GG$ is non-elementary Gromov hyperbolic and $\mu$ is symmetric, finitely supported, and non-degenerate, then there exist $0<c\le C<\infty$ such that for every $\xi\in\partial\GG$ (Martin boundary of $\GG$), the Martin kernel $K(\cdot,\xi)$ is minimal and, for all large $r$,
$$
\exp(c r)\ \preccurlyeq\ \sup_{|x|\le r}K(x,\xi)\ \preccurlyeq\ \exp(C r),
$$
with constants independent of $\xi$. Consequently, within any fixed exponential band there are uncountably many distinct normalised minimals.
\end{prop}

\begin{proof}
For symmetric finitely supported measures on a non-elementary Gromov-hyperbolic group, the (minimal)
Martin boundary identifies with the Gromov boundary $\partial\GG$, and for every
$\xi\in\partial\GG$ the Martin kernel $K(\cdot,\xi)$ is minimal. Observe that for $x,y\in\GG$,
\begin{equation}\label{eq:K-green-horo}
K(x,y)=\exp\bigl(-(d_G(x,y)-d_G(e,y))\bigr).
\end{equation}
Now fix $\xi\in\partial\GG$ and choose any sequence $y_n\in\GG$ converging to $\xi$ in the Martin
compactification. By definition of the Martin boundary, $K(x,y_n)\to K(x,\xi)$ for each fixed $x$.
Taking limits in \eqref{eq:K-green-horo} yields
$$
K(x,\xi)=\exp\bigl(-b_\xi^G(x)\bigr), \quad\forall x\in\GG
$$
where $b_\xi^G:\GG \to \RR$ is the Busemann function corresponding to the Green metric:
$$
b_\xi^G(x):=\lim_{n\to\infty}\bigl(d_G(x,y_n)-d_G(e,y_n)\bigr)=-\log K(x,\xi) \quad x\in\GG.
$$
Since $b_\xi^G(e)=0$ and $b_\xi^G$ is $1$-Lipschitz with respect to $d_G$
(it is a limit of functions $x\mapsto d_G(x,y_n)-d_G(e,y_n)$, each of which is $1$-Lipschitz in $x$),
we have
$$
b_\xi^G(x)\ge -d_G(e,x)\qquad(\forall x\in\GG),
$$
and hence
$$
K(x,\xi)=e^{-b_\xi^G(x)}\le e^{d_G(e,x)}.
$$
Let $S:=\supp(\mu)$ and let $p_*:=\min_{s\in S}\mu(s)>0$. Since $S$ generates $\GG$, any $x\in\GG$ can be written as $x=s_1s_2\dots s_{n}$, with $s_i\in\GG$, $n=|x|$. So, $\PP(\tau_x<\infty)\ge \PP(X_1=s_1, X_2=s_1s_2, \dots, X_n =x)\geq p_*^{|x|}$ and
$$
d_G(e,x)=-\log \PP(\tau_x<\infty)\le (\log(1/p_*))|x|.
$$
Consequently, for $|x|\le r$,
$$
K(x,\xi)\le \exp(Cr)\qquad\text{with }C:=L\log(1/p_*),
$$
uniformly in $\xi$.

Fix $\xi\in\partial\GG$ and let $\gamma:\NN\to\GG$ be a word-geodesic ray from $e$ to $\xi$.
Uniform Ancona inequalities along geodesics imply that for all $m>n$,
$$
G(e,\gamma(m))\ \asymp\ G(e,\gamma(n))G(\gamma(n),\gamma(m)),
$$
with constants independent of $n,m$ and $\xi$.
Dividing by $G(e,\gamma(m))$ and letting $m\to\infty$ (so $\gamma(m)\to\xi$ in the Martin topology) yields
$$
K(\gamma(n),\xi)
=\lim_{m\to\infty}\frac{G(\gamma(n),\gamma(m))}{G(e,\gamma(m))}
 \asymp \frac{1}{G(e,\gamma(n))},
$$
again with constants independent of $n$ and $\xi$.

Since $\GG$ is non-elementary hyperbolic, it is non-amenable, and by Kesten's theorem  $\rho=\limsup_{n\to\infty}p_n(x,y)^{1/n}<1$ for all $x,y\in\GG$. Then $p_k(e,x)=0$ for $k<|x|$ (as $p_k(e,x)>0$ if and only if $x\in(\supp(\mu))^k$) and
$p_k(e,x)\le \rho^k$ for all $k$ (by Cauchy-Schwarz and $\|P^k\|_{2\to 2}=\rho^k$).
Hence
$$
G(e,x)=\sum_{k\ge 0}p_k(e,x) \le \sum_{k\ge |x|}\rho^k
 \le \frac{1}{1-\rho}\rho^{|x|}=\frac{1}{1-\rho}e^{-c_0|x|}
$$
for $c_0:=-(\log\rho)>0$.
Applying this to $x=\gamma(n)$ (where $|\gamma(n)|=n$) gives
$$
K(\gamma(n),\xi) \gtrsim e^{c_0 n},
$$
uniformly in $\xi$. Taking $n=\lfloor r\rfloor$ yields
$$
\sup_{|x|\le r}K(x,\xi) \ge K(\gamma(\lfloor r\rfloor),\xi)\gtrsim e^{c_0\lfloor r \rfloor} \geq e^{c r}
$$
where $c=\frac{c_0}{2}$.

\medskip
Combining everything gives $0<c\le C<\infty$ such that for all large $r$,
$$
\exp(cr)\ \preccurlyeq\ \sup_{|x|\le r}K(x,\xi)\ \preccurlyeq\ \exp(Cr),
$$
uniformly in $\xi\in\partial\GG$.
Distinct $\xi\in\partial\GG$ yield non-proportional kernels, and $\partial\GG$ is uncountable for
non-elementary hyperbolic groups. This proves the corollary.
\end{proof}

\begin{remark}
    On symmetric spaces, Martin kernels $K(x, \xi)$ are given by $e^{-hb}$, where $h$ is the mean curvature of the horospheres, and $b$ is a Busemann function (see \cite{BallmannGromovSchroeder1985}) corresponding to $\xi$ on the ideal boundary (which is identified with the Martin boundary).  
For more details, we refer the reader to \cite{BHM2011}, \cite{BallmannNandiPolymerakis2024}.
\end{remark}

\subsection{Decay of Martin kernels away from the pole}

\begin{lemma}[Green function and hitting probabilities]\label{lem:Green-hitting}
Let $(X_n)_{n\ge0}$ be the right $\mu$-random walk on a finitely generated group $\GG$,
started at $e$, with transition kernel $p_n(x,y)$ and Green function
$$
G(x,y) := \sum_{n\ge0} p_n(x,y).
$$
For $x\in\GG$ let
$$
T_x := \inf\{n\ge0 : X_n = x\}
$$
be the hitting time of $x$, and denote $F(e,x):=\PP_e(T_x<\infty)$.
Then
\begin{equation}\label{eq:Green-hitting-identity}
G(e,x) = F(e,x)G(e,e)\qquad\text{for all }x\in\GG.
\end{equation}
In particular,
\begin{equation}\label{eq:Green-uniform-bound}
0\ \le\ G(e,x)\le G(e,e)\qquad\text{for all }x\in\GG.
\end{equation}
\end{lemma}

\begin{proof}
Fix $x\in\GG$ and consider the local time at $x$:
$$
N_x := \sum_{n\ge0} \mathbf 1_{\{X_n = x\}}.
$$
By definition of the Green function,
$$
\EE_e[N_x] = \sum_{n\ge0} \PP_e(X_n = x) = G(e,x).
$$

Decompose according to whether $x$ is ever hit:
$$
\EE_e[N_x]
= \EE_e\big[N_x\mathbf 1_{\{T_x<\infty\}}\big]
= \PP_e(T_x<\infty)\EE_e\big[N_x\mid T_x<\infty\big].
$$
Let $\mathcal F_n:=\sigma(X_0,\dots,X_n)$ and note that $T_x$ is a stopping time.
On $\{T_x<\infty\}$, since $T_x$ is the \emph{first} visit to $x$, we have
$$
N_x=\sum_{n\ge0}\mathbf 1_{\{X_n=x\}}
=\sum_{n\ge T_x}\mathbf 1_{\{X_n=x\}}
=\sum_{k\ge0}\mathbf 1_{\{X_{T_x+k}=x\}}.
$$
For each $m\ge0$ set $N_x^{(m)}:=\sum_{k=0}^m \mathbf 1_{\{X_{T_x+k}=x\}}$.
By the strong Markov property at $T_x$,
$$
\EE_e\big[N_x^{(m)}\mid \mathcal F_{T_x}\big]
=\EE_{x}\Big[\sum_{k=0}^m \mathbf 1_{\{X_k=x\}}\Big]
=\sum_{k=0}^m p_k(x,x).
$$
Letting $m\to\infty$ and using monotone convergence theorem for conditional expectation yields
$$
\EE_e\Big[\sum_{k\ge0}\mathbf 1_{\{X_{T_x+k}=x\}}\ \Big|\ \mathcal F_{T_x}\Big]
=\sum_{k\ge0}p_k(x,x)=G(x,x)
\qquad\text{on }\{T_x<\infty\}.
$$
Taking expectations and dividing by $\PP_e(T_x<\infty)$ gives $\EE_e[N_x\mid T_x<\infty]=G(x,x)=G(e,e)$.

Substituting into the previous display yields
$$
G(e,x) = \PP_e(T_x<\infty)G(e,e) = F(e,x)G(e,e),
$$
which is \eqref{eq:Green-hitting-identity}. Since $0\le F(e,x)\le1$, we obtain the
uniform bound \eqref{eq:Green-uniform-bound}.
\end{proof}

\begin{prop}[Decay away from the pole under relaxed Ancona]\label{prop:decay-away}
Assume a geometric setting (e.g.\ hyperbolicity) in which the Martin boundary $\Cal M\GG$
is identified with a geometric boundary, and assume the random walk is \emph{transient}
(i.e.\ $G(e,e)<\infty$).

Assume moreover that there exists $\rho:\NN\to(0,\infty)$ with $\rho(R)\to0$ as $R\to\infty$
such that for all $x,y\in\GG$ and for \emph{every} geodesic segment $[e,y]$ from $e$ to $y$,
\beq\label{eq:relaxed-Ancona}
\frac{G(x,y)}{G(e,y)} \le \rho\big(\dist(x,[e,y])\big)G(e,x).
\eeq
Furthermore, assume that for distinct $\xi\neq\eta \in \Cal M\GG$ one has
$$
R(x):=\inf_{\gamma\in\Gamma(e,\xi)}\dist(x,\gamma) \to \infty
\qquad\text{as }x\to\eta \text{ in the Martin topology},
$$
where $\Gamma(e,\xi)$ denotes the set of all geodesic rays from $e$ to $\xi$.
Then
$$
\lim_{x\to \eta} K(x,\xi)=0,
$$
where $K(x,y):=\frac{G(x,y)}{G(e,y)}$ and $K(x,\xi)$ is the Martin limit.
\end{prop}

\begin{proof}
Fix $\xi(\ne\eta)\in\Cal M\GG$ and choose a geodesic ray $\gamma_\xi\in\Gamma(e,\xi)$.
Define
$$
R_\xi(x):=\dist(x,\gamma_\xi):=\inf\{d(x,z):z\in\gamma_\xi\}.
$$
Since $\rho$ may not be decreasing, replace it by
$$
\rho_\ast(R):=\sup_{t\ge R}\rho(t).
$$
Then $\rho_\ast$ is non-increasing, $\rho_\ast(R)\to0$ as $R\to\infty$, and \eqref{eq:relaxed-Ancona}
remains valid with $\rho$ replaced by $\rho_\ast$.

Let $y_k\in\gamma_\xi$ tend to $\xi$ with $[e,y_k]\subset[e,y_{k+1}]$, and interpret $[e,y_k]$
as the initial subsegment of the ray $\gamma_\xi$ from $e$ to $y_k$.
Set $a_k:=\dist(x,[e,y_k])$. Since $[e,y_k]\subset\gamma_\xi$, we have
$$
a_k=\dist(x,[e,y_k])\ \ge\ \dist(x,\gamma_\xi)=R_\xi(x)\qquad\forall k.
$$
Applying \eqref{eq:relaxed-Ancona} gives
$$
K(x,y_k)=\frac{G(x,y_k)}{G(e,y_k)}
 \le \rho_\ast(a_k)G(e,x)
 \le \rho_\ast(R_\xi(x))G(e,x)
\qquad\forall k\in\NN.
$$
Letting $k\to\infty$ along $y_k\to\xi$ and using \eqref{eq:Green-uniform-bound} yields
$$
0\le K(x,\xi)\le \rho_\ast(R_\xi(x))G(e,e).
$$
As $x\to\eta$, the assumption implies $R(x)\to\infty$, and since
$R_\xi(x)=\dist(x,\gamma_\xi)\ge \inf_{\gamma\in\Gamma(e,\xi)}\dist(x,\gamma)=R(x)$,
we have $R_\xi(x)\to\infty$. Therefore $\rho_\ast(R_\xi(x))\to0$, and thus $K(x,\xi)\to0$.
\end{proof}

\subsection{Further upgrades and extensions}

\subsubsection{Upper gauges with unit exponential rate}

We extend Theorem \ref{thm:subexp-upper-no-sym} to gauges $\psi$ with $\limsup_{r\to\infty}\psi(r)^{1/r}=1$.

\begin{lemma}[Unit rate $\Rightarrow$ uniform subexponential domination]\label{lem:unit-rate-to-subex}
Let $\psi:\NN\to(0,\infty)$ be a function such that $\limsup_{r\to\infty}\psi(r)^{1/r}=1$. Then for every $\varepsilon>0$ there exists $C_\varepsilon$ such that
$$
\psi(r)\ \le\ C_\varepsilon e^{\varepsilon r}\qquad\forall r\ge 0.
$$
\end{lemma}

\begin{proof}
Fix $\varepsilon>0$. Since $e^\varepsilon > 1$, by the definition of the limit superior, there exists $R_0$ such that for all $r \ge R_0$, $\psi(r)^{1/r} < e^{\varepsilon}$. This implies $\psi(r) < e^{\varepsilon r}$ for $r \ge R_0$.
Let
$$
C_\varepsilon = \max\Big(\{1\} \cup \{\psi(r)e^{-\varepsilon r} : 0 \le r < R_0\}\Big).
$$
Then $\psi(r) \le C_\varepsilon e^{\varepsilon r}$ holds for all $r\ge 0$.
\end{proof}

\begin{theorem}[Abundance under unit-rate upper gauges]\label{thm:unit-rate-abundance}
Let $\mu$ be finitely supported and non-degenerate (no symmetry assumed) on an infinite group $\GG$. Let $\psi$ be a growth gauge with 
$$
\displaystyle\limsup_{r\to\infty}\psi(r)^{1/r}=1.
$$ 
If there exists a non-constant minimal $\mu$-harmonic $h$ with upper growth $\preccurlyeq \psi$, then there are infinitely many minimal $\mu$-harmonic functions $h'$ having upper growth $\preccurlyeq \psi$.
\end{theorem}

\begin{proof}
By Lemma \ref{lem:unit-rate-to-subex}, the growth of $h$ is subexponential: for every $\varepsilon>0$ there exists $C_\varepsilon$ such that
$$
|h(x)|\ \le\ C_\varepsilon e^{\varepsilon|x|}\qquad\forall x\in\GG.
$$
Apply Theorem \ref{thm:subexp-upper-no-sym} to conclude the normalised orbit $\{\widehat{g.h}:g\in\GG\}$ is infinite. Growth preservation under translation follows from Lemma \ref{lem:equivariance}(2), since $\psi$ is a growth gauge.
\end{proof}

\subsubsection{Relative hyperbolicity and Floyd boundary: abundance in exponential bands}

We record an abstract transfer principle from a geometric compactification of the Cayley graph to the Martin boundary.

\begin{definition}[Geometric compactification / boundary]\label{def:geometric-boundary}
Let $d(x,y):=|x^{-1}y|$ be the word metric on $\GG$.
A \emph{geometric compactification} of $\GG$ is a compact metrizable space
$\overline\GG=\GG\sqcup\mathcal B$ such that
\begin{itemize}
\item $\GG$ is a dense open subset of $\overline\GG$ and the left action of $\GG$ on itself extends to a continuous action by homeomorphisms on $\overline\GG$;
\item for every $\xi\in\mathcal B$ there exists a \emph{geodesic} ray $\gamma_\xi:\NN\to\GG$ with $\gamma_\xi(0)=e$, $|\gamma_\xi(n)|=n$, and $\gamma_\xi(n)\to \xi$ in $\overline\GG$.
\end{itemize}
We call $\mathcal B$ the (geometric) boundary.
\end{definition}

\begin{prop}[Abstract boundary-to-kernel transfer]\label{prop:abstract-transfer}
Let $\mu$ be finitely supported and symmetric on $\GG$. Suppose there is a geometric compactification $\overline\GG=\GG\sqcup\mathcal B$
as in Definition \ref{def:geometric-boundary}, and a continuous $\GG$-equivariant map $\iota:\mathcal B\to\Cal M\GG$
whose image consists of minimal Martin points.
Assume there exist $A\ge 1$, $c,C>0$ and functions $b_\xi:\GG\to\RR$ such that:
\begin{itemize}
  \item (\emph{Uniform coarse $1$-Lipschitz}) there exists $C_{\mathrm{Lip}}\ge 0$ with
  $$
  |b_\xi(x)-b_\xi(y)|\ \le\ d(x,y)+C_{\mathrm{Lip}}
  \qquad\forall \xi\in\mathcal B,\ \forall x,y\in\GG,
  $$
  and $b_\xi(e)=0$ for all $\xi$;
  \item for all $\xi\in\mathcal B$ and $x\in\GG$,
  \begin{equation}\label{eq:kernel-busemann}
  A^{-1}e^{-Cb_\xi(x)} \le K(x,\iota(\xi)) \le A e^{-cb_\xi(x)};
  \end{equation}
  \item (\emph{Uniform coarse linearity along geodesic rays}) there exist constants $a_1,a_2>0$ and $B\ge 0$ such that for every $\xi\in\mathcal B$,
  along the geodesic ray $\gamma_\xi$ from Definition \ref{def:geometric-boundary} we have for all $n\in\NN$,
  $$
  -a_2 n - B  \le b_\xi(\gamma_\xi(n)) \le -a_1 n + B.
  $$
\end{itemize}
Then there exist $0<c'\le C'<\infty$ such that for all large $r$,
$$
\exp(c'r) \lesssim \sup_{|x|\le r}K(x,\iota(\xi)) \lesssim \exp(C' r),
$$
with constants independent of $\xi\in\mathcal B$.
If $\iota(\mathcal{B})$ is uncountable (for instance, if $\mathcal{B}$ is uncountable and $\iota$ is injective), this yields uncountably many pairwise non-proportional normalised minimals within a uniform exponential band.
\end{prop}

\begin{proof}
Fix $\xi\in\mathcal B$.

\emph{Upper bound.}
If $|x|\le r$, then using $b_\xi(e)=0$ and the uniform coarse Lipschitz property,
$$
b_\xi(x) \ge -|x|-C_{\mathrm{Lip}} \ge -r-C_{\mathrm{Lip}}.
$$
Hence, by the upper bound in \eqref{eq:kernel-busemann},
$$
K(x,\iota(\xi))
 \le A e^{-c b_\xi(x)}
 \le A e^{c(r+C_{\mathrm{Lip}})}
\lesssim e^{c r}
$$
independent of $\xi$.

\emph{Lower bound.}
Take $x=\gamma_\xi(r)$, so $|x|=r$. Using the coarse-linearity upper bound on $b_\xi(\gamma_\xi(r))$ and the lower bound in \eqref{eq:kernel-busemann},
$$
K(\gamma_\xi(r),\iota(\xi))
 \ge A^{-1} e^{-C b_\xi(\gamma_\xi(r))}
 \ge A^{-1} e^{-C(-a_1 r + B)}
\ge A^{-1} e^{-CB} e^{Ca_1 r}.
$$
Therefore,
$$
\sup_{|x|\le r} K(x,\iota(\xi))
 \ge K(\gamma_\xi(r),\iota(\xi))
 \gtrsim e^{c' r}
$$
for some $c'>0$ independent of $\xi$.

This gives the claimed uniform two-sided exponential bounds. Distinct points in $\iota(\mathcal B)$ lie in the minimal Martin boundary, hence their kernels are pairwise non-proportional.
\end{proof}

\subsubsection{Projective dynamics on the cone and eigencharacters}

Let $\mathcal C:=\{h>0: Ph=h\}$ and $\mathcal C_1:=\{h\in\mathcal C:\ h(e)=1\}$. For $g\in\GG$ set
$$
L_g:\ \mathcal C_1\to\mathcal C_1,\qquad L_g(h):=\widehat{g.h}.
$$

\begin{theorem}[Finite invariant family $\Rightarrow$ finite-state Markov factor]\label{thm:finite-stationary-factor}
Let $\mathcal M=\{h_1,\dots,h_m\}\subset\mathcal C_1$ be a finite $L$-invariant family of normalised minimals. Let $(X_n)_{n\ge 1}$ be the i.i.d.\ increments of the random walk. The process $(H_n)_{n\ge 0}$ defined by $H_0 \in \mathcal{M}$ and
$$
H_{n+1} = L_{X_{n+1}}(H_n)
$$
is a Markov chain on $\mathcal{M}$. This chain admits a stationary distribution.
\end{theorem}

\begin{proof}
\emph{Well-definedness.}
Since $\mathcal M$ is $L$-invariant, in particular $L_g(\mathcal M)\subseteq\mathcal M$ for every
$g\in\supp\mu$. Hence, starting from any deterministic $H_0\in\mathcal M$, the recursion
$H_{n+1}=L_{X_{n+1}}(H_n)$ stays in $\mathcal M$ almost surely.%
\footnote{If $H_0$ is random, it suffices to assume $H_0$ is independent of $(X_n)_{n\ge1}$
and take $\mathcal F_n:=\sigma(H_0,X_1,\dots,X_n)$ below.}

Let $\mathcal F_n:=\sigma(X_1,\dots,X_n)$. Then $H_n$ is $\mathcal F_n$-measurable for every $n$.

\emph{Markov property.}
Fix $n\ge0$ and states $h_i,h_j\in\mathcal M$. On the event $\{H_n=h_i\}$ we have
$H_{n+1}=L_{X_{n+1}}(h_i)$, so
$$
\PP\big(H_{n+1}=h_j \big| \mathcal F_n\big)
=\sum_{i=1}^m \mathbf 1_{\{H_n=h_i\}}
\PP\big(L_{X_{n+1}}(h_i)=h_j \big| \mathcal F_n\big).
$$
Since $X_{n+1}$ is independent of $\mathcal F_n$ and distributed as $\mu$, the conditional probability
is a constant given by
$$
P_{ij}:=\PP\big(L_{X_1}(h_i)=h_j\big)
=\sum_{g\in\GG:L_g(h_i)=h_j}\mu(g),
$$
and therefore
$$
\PP\big(H_{n+1}=h_j \big| \mathcal F_n\big)
=\sum_{i=1}^m \mathbf 1_{\{H_n=h_i\}}P_{ij}.
$$
This shows that the conditional law of $H_{n+1}$ given $\mathcal F_n$ depends only on $H_n$,
so $(H_n)$ is a time-homogeneous Markov chain on the finite state space $\mathcal M$
with transition matrix $P=(P_{ij})$.
Moreover, for each $i$,
$$
\sum_{j=1}^m P_{ij}
=\PP\big(L_{X_1}(h_i)\in\mathcal M\big)=1,
$$
using $L_g(h_i)\in\mathcal M$ for all $g\in\supp\mu$, so $P$ is stochastic.

\emph{Existence of a stationary distribution.}
Since the state space of $(H_n)$ is finite, it has a stationary distribution.
\end{proof}

\subsection{An explicit construction of a nontrivial Martin kernel with subexponential growth}\label{sec:Appendix_B}

We give an explicit construction of a positive harmonic function whose associated Martin kernel grows \emph{strictly subexponentially} in the space variable.
\medskip

\subsubsection{A finite-support example with killing: random walks in cones}\label{subsubsec:cone-killed}

Let $d\ge 2$. Consider a centered random walk $(S_n)_{n\ge 0}$ on $\Z^d$ with
\emph{finite support} increment distribution, strongly aperiodic in the sense of
\cite[(H3)]{DurajRaschelTarragoWachtel2022MartinCones}, and with nondegenerate
covariance matrix. Fix an \emph{open} convex cone $K\subset\RR^d$ with nonempty
interior, and kill the walk at the first exit time
\[
\tau_K:=\inf\{n\ge 0:\ S_n\notin K\}.
\]
Assume the cone-adapted hypotheses \cite[(H1)-(H5)]{DurajRaschelTarragoWachtel2022MartinCones}
(in particular, asymptotic strong irreducibility in $K$). Since the increments have
finite support, the moment hypothesis \cite[(M1)]{DurajRaschelTarragoWachtel2022MartinCones}
is automatically satisfied.

For $x,y\in K\cap\Z^d$, let
\[
G_K(x,y):=\sum_{n\ge 0}\PP_x\big(S_n=y,\ \tau_K>n\big)
\]
be the Green function of the killed walk.

\medskip\noindent
As in \cite[(H2) and the discussion following it]{DurajRaschelTarragoWachtel2022MartinCones},
one may apply an invertible linear change of variables to reduce to the
identity-covariance setting. This transforms the ambient lattice and the cone, but
does not affect the qualitative Martin-boundary conclusions. For simplicity we
describe the associated Brownian Dirichlet problem in this normalisation; without
it one should replace $\Delta$ below by the corresponding constant-coefficient
elliptic operator.

\begin{theorem}[Singleton Martin boundary and interior polynomial growth]\label{thm:cone}
Fix an interior basepoint $x_\ast\in K\cap\Z^d$ with $\dist(x_\ast,\partial K)>0$.
Then:
\begin{enumerate}
\item[(i)] \textbf{Singleton Martin boundary.}
The Martin boundary of the killed walk in $K$ reduces to a single point.
Consequently, there is a \emph{unique} (up to scale) positive harmonic function
$h:K\cap\Z^d\to(0,\infty)$ (with Dirichlet boundary condition $h\equiv0$ on $K^c$),
and for any sequence $(y_n)$ in $K\cap\Z^d$ with $\|y_n\|\to\infty$,
\[
\frac{G_K(x,y_n)}{G_K(x_\ast,y_n)}\longrightarrow \frac{h(x)}{h(x_\ast)}
\qquad\text{for every fixed }x\in K\cap\Z^d.
\]
In particular, the associated Martin kernel at infinity may be taken as
$K_K(x,\infty):=h(x)/h(x_\ast)$.

\item[(ii)] \textbf{Interior asymptotics and polynomial growth.}
Let $u$ be the (unique up to scale) positive solution of the Brownian Dirichlet problem
\[
\Delta u=0\ \text{ in }K,\qquad u=0\ \text{ on }\partial K.
\]
Then $u$ is homogeneous: $u(tx)=t^p u(x)$
for all $t>0$ and $x\in K$, for some $p>0$.

If the walk has bounded jumps (in particular, finite support), then there exists a
constant $c>0$ such that for every $\varepsilon>0$,
\[
h(x)=c\,u(x)\,(1+o(1))
\qquad\big(\|x\|\to\infty,\ \dist(x,\partial K)\ge \varepsilon\|x\|\big).
\]
Writing $\Sigma:=K\cap\mathbb S^{d-1}$ and $\varphi(\theta):=u(\theta)$, this is equivalently
\[
h(x)=c\,\|x\|^{p}\,\varphi\!\left(\frac{x}{\|x\|}\right)\,(1+o(1))
\qquad\big(\|x\|\to\infty,\ \dist(x,\partial K)\ge \varepsilon\|x\|\big),
\]
where $\varphi$ is strictly positive on $\Sigma$ (and vanishes on $\partial\Sigma$).
In particular $K_K(x,\infty)$ is nontrivial and grows polynomially (hence
subexponentially) in $\|x\|$ along interior directions.
\end{enumerate}
\end{theorem}

\begin{proof}[Proof sketch]
For (i), \cite[Theorem 1.4]{DurajRaschelTarragoWachtel2022MartinCones} identifies the
Martin boundary of the random walk killed on exiting $K$ with a single point,
corresponding to their harmonic function $V$. The ratio-limit statement for Green
functions in \cite[Corollary 4.2]{DurajRaschelTarragoWachtel2022MartinCones} yields the displayed limit with $h$ proportional to $V$; the normalisation by
$h(x_\ast)$ gives the Martin kernel $K_K(\cdot,\infty)$.

For (ii), the Brownian Dirichlet problem and the homogeneity $u(tx)=t^p u(x)$ are
recalled in \cite[(1.3)-(1.4) and the discussion following them]{DurajRaschelTarragoWachtel2022MartinCones}.
When the increments are bounded, Denisov-Wachtel show that the discrete harmonic
function satisfies $V(z)\sim u(z)$ \emph{uniformly} as $z\to\infty$ in interior
directions (i.e.\ $\dist(z,\partial K)\ge\varepsilon|z|$); see \cite[Lemma 13]{DenisovWachtel2015}.
Since $h$ is unique up to scale and proportional to $V$, this gives
$h(x)=c\,u(x)(1+o(1))$ and hence the stated polynomial growth.
\end{proof}

\section{Green geometry collapse, asymptotics beyond finite first moment and strong Liouville}\label{sec:speed}

Throughout this section, we assume that $\mu$ is a non-degenerate measure generating a transient random walk on a finitely generated group $\GG$. 
\subsection{Green speed and strong Liouville property}
First, we prove a technical lemma.
\begin{lemma}\label{lem:conv_prob}
    Let $(X_n)_{n\ge1}$ and $(Y_n)_{n\ge1}$ be sequences of real-valued random variables defined on the same probability space $(\Omega,\mathcal{F},\mathbb{P})$.  
Assume that for each $\omega \in \Omega$ there exists $k(\omega) \in \mathbb{N}$ such that
$$
0 \le X_n(\omega) \le Y_n(\omega), \quad \forall n \ge k(\omega).
$$
If $Y_n \to 0$ in probability, then $X_n \to 0$ in probability.
\end{lemma}
\begin{proof}
    Fix $\varepsilon > 0$. For each $n$,
$$
\{ X_n > \varepsilon \} \subset \{ Y_n > \varepsilon \} \ \cup \ \{ X_n > Y_n \}.
$$
Taking probabilities gives
\beq\label{eq:inequality_DCT}
\mathbb{P}(X_n > \varepsilon) \ \le \ \mathbb{P}(Y_n > \varepsilon) + \mathbb{P}(X_n > Y_n).
\eeq
Since $Y_n \to 0$ in probability, we have $\mathbb{P}(Y_n > \varepsilon) \to 0$ as $n \to \infty$. It remains to show that $\mathbb{P}(X_n > Y_n) \to 0$. Since for each $\omega$ there exists $k(\omega)$ such that $X_n(\omega) \le Y_n(\omega)$ for all $n \ge k(\omega)$, we have that
$$
\mathbf{1}_{\{ X_n > Y_n \}} \to 0 \quad \text{pointwise.}
$$
Moreover, $0 \le \mathbf{1}_{\{ X_n > Y_n \}} \le 1$, so by dominated convergence theorem,
$$
\mathbb{P}(X_n > Y_n) = \mathbb{E}\big[ \mathbf{1}_{\{ X_n > Y_n \}} \big] \to 0.
$$
Using \eqref{eq:inequality_DCT} we have
$$
\lim_{n \to \infty} \mathbb{P}(X_n > \varepsilon) = 0.
$$
Since $\varepsilon > 0$ is arbitrary, it follows that $X_n \to 0$ in probability.
\end{proof}

We now establish a relation between strong Liouville property and zero Green speed. Recall that the \emph{speed (or linear drift)} of a $\mu$-random walk on $\GG$ is defined by the almost sure limit
$$
    \lim_{n \to \infty} \frac{|X_n|}{n}.
$$
Note that the above limit exists almost surely whenever $\mu$ has finite first moment (see \cite{KarlLedrappier2007}).

\begin{definition}[Speed finite in probability]
We say that the (word) speed is \emph{finite in probability} if the sequence $(|X_n|/n)_{n\ge1}$ is
\emph{bounded in probability}, i.e.
$$
\forall\varepsilon>0 \text{ there exists } M>0 \text{ and } N\in\mathbb N \text{ such that } 
\mathbb P\left(\frac{|X_n|}{n}>M\right) < \varepsilon\qquad\forall n\ge N.
$$
Equivalently,
$$
\lim_{M\to\infty} \limsup_{n\to\infty} \mathbb P\left(\frac{|X_n|}{n}>M\right) = 0.
$$
\end{definition}

\noindent
(In particular, if $\frac{|X_n|}{n}\to v$ in probability or almost surely for some finite $v$, then the speed is finite in probability.)

It is well known that if $\mu$ is non-degenerate and $(\GG, \mu)$ is Liouville, then speed of the random walk corresponding to $\mu$ is $0$ (\cite[Corollary 2]{KarlLedrappier2007}). In spirit of the latter result, we have the following theorem:

\begin{theorem}\label{thm:form_entropy_harmonic}
Let $\mu$ be a symmetric, non-degenerate measure on $\GG$ which generates a transient random walk $(X_n)$. If the Martin boundary of $\GG$ is trivial (i.e. consists of a single point), then we have 
\beq\label{eq:Green_metric_rate}
\lim_{r \to \infty} \sup_{|x| = r}\dfrac{d_G(e, x)}{r} = 0.
\eeq   
Furthermore, we also have 
    \beq\label{eq:Green_speed_zero}
    l_{G}(\omega) = 0,
    \eeq
    for every path $\omega$ in the Wiener space of the $\mu$-random walk $(X_n)$ on $\GG$ along which the speed $\displaystyle \lim_{n \to \infty} \frac{|X_n(\omega)|}{n}$ is finite (here $l_G$ is the {\em Green speed} of the random walk generated by $\mu$ on $\GG$). As a consequence, if speed is finite (in probability), then Green speed converges to zero (in probability). In particular, if $(\GG, \mu)$ is strong Liouville and $\mu$ is a symmetric, non-degenerate measure with superexponential moments that generates a transient random walk, then \eqref{eq:Green_metric_rate} and \eqref{eq:Green_speed_zero} hold.
\end{theorem}

For the notion of superexponential moments, see Definition \ref{def:smooth_measure}
(in this general context, we also refer the reader to \cite{Gouezel2015}). Here, the
Green speed is defined by the limit (whenever it exists almost surely)
\beq\label{def:green_speed}
l_{G} = \lim_{n \to \infty} \frac{d_G(e, X_n)}{n},
\eeq
where 
\beq\label{def:Green_dist}
d_G(x, y) := \log G(e, e) - \log G(x, y)
\eeq
is the Green distance. 

Let $H(\mu)$ denote the (Shannon) entropy of the random walk generated by $\mu$: 
$$
    H(\mu) = -\sum_{x \in \GG}\mu(x)\log \mu(x).    
$$
If $H(\mu) < \infty$, then via Kingman's subadditive lemma \cite{Kingman1968},
$\lim_{n\to\infty} \frac{H(\mu^{(n)})}{n}$ exists ($\mu^{(n)}$ denotes the $n$-fold
convolution of $\mu$ with itself). In such a case, the {\em asymptotic entropy} (see
\cite{KV1983}, \cite[Section~9]{pete2014probability}) of the random walk on $\GG$
generated by $\mu$ is defined as
$$
    \rho(\mu) := \lim_{n\to\infty} \frac{H(\mu^{(n)})}{n}.
$$

Note that under the assumption $H(\mu) < \infty$, the Green speed $l_G$ exists almost
surely and one has $\rho(\mu) = l_G$ \cite[Theorem~1.1]{BlachHai}. Hence, if $\GG$ is
Liouville and $H(\mu) < \infty$, then $l_G = \rho(\mu) = 0$. However, in
Theorem \ref{thm:form_entropy_harmonic} we do \emph{not} impose any restriction on
$H(\mu)$. Moreover, the proof of Theorem \ref{thm:form_entropy_harmonic} below indicates
that $H(\mu) < \infty$ is not a necessary condition for $l_G$ to exist almost surely.

\cite[Corollary~2]{KarlLedrappier2007} also has a converse (positive speed implies
non-Liouville). We cannot expect a converse of Theorem \ref{thm:form_entropy_harmonic}
to hold. For instance, take a finitely supported non-degenerate measure $\mu$ on the
Lamplighter group $\ZZ_2 \wr \ZZ$. Then $(\ZZ_2 \wr \ZZ, \mu )$ is Liouville, since it
has trivial Poisson boundary \cite[Theorem~3.3]{Kaimanovich1991}. Since $H(\mu) < \infty$,
we have $l_{\ZZ_2 \wr \ZZ} = \rho(\mu) = 0$. However, $(\ZZ_2 \wr \ZZ, \mu)$ cannot have
the strong Liouville property as it has exponential growth (see the discussion below
Theorem \ref{prop:exp_growth_delta_zero}). Having said that, we have the following related
question:

\begin{question}
Can we impose some extra condition(s) along with $l_G = 0$ to ensure that $(\GG, \mu)$ is strong Liouville?
\end{question}

\begin{proof}[Proof of Theorem \ref{thm:form_entropy_harmonic}]
Fix a finite symmetric generating set $S$ for $\GG$, and let $|\cdot|$ denote the word
length with respect to $S$.

\medskip\noindent
Recall that for $x,y\in\GG$,
$$
G(x,y) = \sum_{n\ge0} p_n(x,y),\qquad
K(x,y):=\frac{G(x,y)}{G(e,y)}
$$
is the Martin kernel based at $e$. The Martin compactification $\widehat\GG$ is the
minimal compactification of $\GG$ to which all $K(x,\cdot)$ extend continuously. By
assumption, the Martin boundary $\partial_M\GG:=\widehat\GG\setminus\GG$ consists of a
single point, say $\{\xi\}$.

For each fixed $x\in\GG$, the extension $K(x,\cdot)$ is a positive harmonic function on
$\GG$ (in the first variable) normalised by $K(e,\cdot)\equiv 1$. Since the minimal
Martin boundary is reduced to a single point, the only extremal positive harmonic function
normalised at $e$ is the constant function $1$. Hence
$$
K(x,\xi) = 1\qquad\text{for every }x\in\GG.
$$
By continuity at $\xi$, this implies
\begin{equation}\label{eq:K-to-1}
K(x,y) = \frac{G(x,y)}{G(e,y)} \longrightarrow 1\qquad\text{as }y\to\xi,
\text{ for every fixed }x\in\GG.
\end{equation}

Now fix $s\in S$ and consider the function
$$
f_s(z) := \frac{G(e,zs)}{G(e,z)}\qquad(z\in\GG).
$$
Using symmetry $G(a,b)=G(b,a)$ and left-invariance, one checks that $f_s$ can be written
in terms of the Martin kernel $K(\cdot,\cdot)$, and \eqref{eq:K-to-1} shows that $f_s(z)\to1$
whenever $z\to\xi$ in $\widehat\GG$. Hence, for every
$\varepsilon>0$ there exists $R_s(\varepsilon)>0$ such that
\begin{equation}\label{eq:fs-uniform}
|f_s(z)-1|\ \le\ \varepsilon\qquad\text{for all }z\in\GG\text{ with }|z|\ge R_s(\varepsilon).
\end{equation}
As $S$ is finite, taking $R(\varepsilon):=\max_{s\in S} R_s(\varepsilon)$ yields
\begin{equation}\label{eq:fs-uniform-all}
|f_s(z)-1|\ \le\ \varepsilon\qquad\text{for all }s\in S\text{ and all }z\in\GG\text{ with }|z|\ge R(\varepsilon).
\end{equation}

Define the Green increment along oriented edges
$$
\phi_s(z) := \log\frac{G(e,z)}{G(e,zs)} = -\log f_s(z).
$$
From \eqref{eq:fs-uniform-all} and the fact that $\log$ is Lipschitz near $1$ we deduce:
for every $\varepsilon>0$ there exists $R(\varepsilon)$ such that
\begin{equation}\label{eq:phi-small}
|\phi_s(z)|\ \le\ \varepsilon\qquad\text{for all }s\in S\text{ and all }z\text{ with }|z|\ge R(\varepsilon).
\end{equation}

\medskip\noindent
Let $x\in\GG$ with $|x|=n$, and fix a word geodesic decomposition
$$
x = s_1 s_2\cdots s_n,\qquad s_i\in S.
$$
Set $z_0=e$ and $z_k:=s_1\cdots s_k$ for $1\le k\le n$. Then for each $k$,
$$
G(e,z_{k+1}) = G(e,z_k s_{k+1}),\qquad
\phi_{s_{k+1}}(z_k) = \log\frac{G(e,z_k)}{G(e,z_{k+1})}.
$$
Telescoping gives
$$
\log G(e,x) - \log G(e,e)
= \sum_{k=0}^{n-1} \bigl(\log G(e,z_{k+1}) - \log G(e,z_k)\bigr)
= -\sum_{k=0}^{n-1}\phi_{s_{k+1}}(z_k),
$$
hence
\begin{equation}\label{eq:dG-telescope}
d_G(e,x) = \log G(e,e)-\log G(e,x) = \sum_{k=0}^{n-1}\phi_{s_{k+1}}(z_k).
\end{equation}

Fix $\varepsilon>0$ and let $R:=R(\varepsilon)$ be as in \eqref{eq:phi-small}. Define
$$
M_R := \max\{|\phi_s(z)|:\ s\in S,\ |z|\le R\} < \infty.
$$
Splitting the sum in \eqref{eq:dG-telescope} around $R$ gives
$$
|d_G(e,x)|\le \sum_{0\le k\le R}|\phi_{s_{k+1}}(z_k)|+ \sum_{R<k\le n-1}|\phi_{s_{k+1}}(z_k)|\le M_R (R+1) + \varepsilon (n-R).
$$
Therefore
$$
\frac{d_G(e,x)}{|x|} = \frac{d_G(e,x)}{n}
\ \le\ \varepsilon\left( 1-\frac{R}{n}\right) + \frac{M_R (R+1)}{n}.
$$
Taking the supremum over $|x|=n$ and letting $n\to\infty$ yields
$$
\limsup_{n\to\infty}\ \sup_{|x|=n}\frac{d_G(e,x)}{n} \ \le\ \varepsilon.
$$
Since $\varepsilon>0$ was arbitrary, this proves \eqref{eq:Green_metric_rate}.

For later use, define
$$
\rho(r) := \sup_{|x|=r}\frac{d_G(e,x)}{r}\qquad(r\in\NN),
$$
so that \eqref{eq:Green_metric_rate} is equivalent to $\rho(r)\to 0$ as $r\to\infty$.
Note also that $0\le d_G(e,x)\le \log G(e,e)$, so $\rho(r)$ is bounded and non-negative.

\medskip\noindent
\emph{Green speed along paths of finite word speed.}
Let $\omega$ be a path of the random walk and write $X_n=X_n(\omega)$. Define
$$
v(\omega) := \lim_{n\to\infty}\frac{|X_n(\omega)|}{n},
$$
assuming the limit exists and is finite. For each $n$,
$$
0 \le \frac{d_G(e,X_n)}{n}
   = \frac{d_G(e,X_n)}{|X_n|}\cdot\frac{|X_n|}{n}
    \le \rho(|X_n|)\frac{|X_n|}{n}.
$$

Fix $\varepsilon>0$ and choose $R$ so large that $\rho(r)\le\varepsilon$ for all $r\ge R$.
Write
$$
Y_n(\omega) := \frac{d_G(e,X_n(\omega))}{n}.
$$
For large $n$ (so that $||X_n|/n - v(\omega)|\le 1$) we have:
\begin{itemize}
  \item If $|X_n(\omega)|\ge R$, then
  $$
  Y_n(\omega) \le \rho(|X_n(\omega)|)\frac{|X_n(\omega)|}{n}
   \le \varepsilon(v(\omega)+1).
  $$
  \item If $|X_n(\omega)|<R$, then $d_G(e,X_n(\omega))$ ranges over a finite set, so
  $Y_n(\omega)\le C_R/n$ for some $C_R<\infty$, and hence $Y_n(\omega)\to 0$ along such indices.
\end{itemize}
Combining the two cases, we see that $Y_n(\omega)\to 0$ as $n\to\infty$ for every path
$\omega$ along which $v(\omega)$ exists and is finite. This proves \eqref{eq:Green_speed_zero}.

\medskip\noindent
\emph{Green speed converges to $0$ in probability when speed is finite in probability.}
Assume now that the word speed is finite in probability: for every $\eta>0$ there exists
$M>0$ and $N_0\in\NN$ such that
$$
\PP\Big(\frac{|X_n|}{n}>M\Big) \le \eta\qquad\text{for all }n\ge N_0.
$$
Define
$$
X_n' := \frac{d_G(e,X_n)}{n},\qquad
Y_n' := \rho(|X_n|)\cdot\frac{|X_n|}{n}.
$$
By construction, $0\le X_n'\le Y_n'$ for all $n$ and all sample points. We claim that
$Y_n'\to 0$ in probability. Once this holds, Lemma \ref{lem:conv_prob} applied to
$(X_n')$ and $(Y_n')$ yields $X_n'\to 0$ in probability, i.e.\ the Green speed converges
to zero in probability.

Fix $\varepsilon>0$ and choose $M$ and $N_0$ as above with $\eta:=\varepsilon/4$.
Next, choose $R$ large enough so that $\rho(r)\le \varepsilon/(2M)$ for all $r\ge R$.
Finally, since $\rho(r)$ is bounded on $r\le R$, let $L_R:=\max_{1\le r\le R}\rho(r)<\infty$.
Pick $N_1\ge N_0$ so large that $L_R R/N_1\le \varepsilon/4$.

For $n\ge N_1$ we have
$$
\PP(Y_n'>\varepsilon)\le\PP\big(|X_n|\ge R, Y_n'>\varepsilon\big),
$$
since on $\{|X_n|<R\}$,
$$
Y_n' = \rho(|X_n|)\frac{|X_n|}{n} \le L_R\frac{R}{n} \le \frac{\varepsilon}{4}.
$$
On the event $\{|X_n|\ge R,\ |X_n|/n\le M\}$ we have
$$
Y_n' = \rho(|X_n|)\frac{|X_n|}{n} \le \frac{\varepsilon}{2M}M = \frac{\varepsilon}{2}
<\varepsilon.
$$
Hence,
$$
\PP(Y_n'>\varepsilon)\le \PP\Big(\frac{|X_n|}{n}>M\Big)\le \frac{\varepsilon}{4}
$$
for all $n\ge N_1$. This shows $Y_n'\to 0$ in probability, as claimed.

\medskip\noindent
\emph{Strong Liouville + superexponential moments.}
If $\mu$ has superexponential moments, then for every $\xi\in\Cal M(\GG)$ the Martin
kernel $K(\cdot,\xi)$ is $\mu$-harmonic (see \cite[Lemma~7.1]{GekhtmanGerasimovPotyagailoYang2021}). Hence, if $(\GG,\mu)$ is strong Liouville, every positive harmonic function is constant, so the only normalised positive harmonic function is the constant~$1$. In particular, each $K(\cdot,\xi)$ is constant, so the Martin boundary reduces to a single point. We are therefore in the setting of the first part of the theorem, and \eqref{eq:Green_metric_rate} and \eqref{eq:Green_speed_zero} follow.
\end{proof}

\subsection{Proofs of \Cref{thmG} and \Cref{thmH}}

Let $(X_n)_{n\ge 0}$ denote the $\mu$-random walk on $\GG$, and $S_0, S_1, \dots, S_n, \dots$ denote the independent jumps of $(X_n)_{n\ge 0}$. We claim that $\sum_{x \in \GG}|x|\mu(x) < \infty$ if and only if $\lim_{n \to \infty} \frac{\mathbb{E}[|X_n|]}{n} < \infty$. Note that $\lim_{n \to \infty} \frac{\mathbb{E}[|X_n|]}{n}$ always exists as a finite or infinite limit via an application of Fekete's lemma. Suppose $\sum_{x \in \GG}|x|\mu(x) < \infty$. Then, for all $n \in \NN$,
\begin{align*}
    \frac{\mathbb{E}[|X_n|]}{n} 
    &\leq \frac{1}{n}\sum_{k =0}^{n-1} \mathbb{E}[|S_k|]  
     = \sum_{x \in \GG} |x| \mu(x) < \infty. 
\end{align*}

Conversely, assume $\sum_{x \in \GG} |x| \mu(x) = \infty$. Then, $\mathbb{E}[|S_1 x|] = \sum_{y \in \GG} |yx|\mu(y) \geq \sum_{y \in \GG} (|y| - |x| )\mu(y) = \infty$ for all $x \in \GG$. Hence for all $n \in \NN$,
\begin{align*}
 \mathbb{E}[|X_n|] &= \sum_{x_1, x_2, \dots, x_n \in \GG}|x_1x_2\dots x_n |\mu(x_1)\mu(x_2)\dots\mu(x_n) \\
 &=\sum_{x \in \GG}{\mathbb{E}[|S_1 x|]}\sum_{x_2x_3\dots x_n = x}\mu(x_2)\mu(x_3)\dots\mu(x_n) = \infty.
\end{align*}
Therefore, $\lim_{n \to \infty} \frac{\mathbb{E}[|X_n|]}{n} = \infty$. Note, however, that this does not imply that $\frac{|X_n|}{n}$ has a limit along a set of random paths of non-zero measure, as the usual application of Kingman's subadditivity lemma is not possible here. What we can show, however, is that infinite first moment rules out an almost sure finite speed. 

We now show that the speed converges to zero in probability under rather benign assumptions on the heat kernel. Assume there exist constants $C,c>0$ and exponents $\alpha,\beta>0$ such that for all $n\ge1$ and all $x\in\GG$,
\begin{equation}\label{HK}
p_n(e,x)\le C\exp\Big(-c\frac{|x|^{\alpha}}{n^{\beta}}\Big).
\end{equation}

\begin{prop}[Subexponential volume: threshold $\alpha\ge 1+\beta$]\label{prop:subexp}
Assume $\GG$ has subexponential growth, i.e.\ for every $\delta>0$ there exists $R_\delta$ such that
$$
|B(e,r)|\le \exp(\delta r)\qquad\text{for all }r\ge R_\delta,
$$
and assume \eqref{HK} holds with $\alpha\ge 1+\beta$.
Then, for every $\varepsilon>0$,
$$
\PP\left(\frac{|X_n|}{n}>\varepsilon\right)\to0 \text{ as }n\to\infty.
$$
\end{prop}

\begin{proof}
Fix $\varepsilon>0$. Using \eqref{HK} and $|S(e,r)|\le |B(e,r)|$,
\begin{align*}
\PP\left(\tfrac{|X_n|}{n}>\varepsilon\right)
&= \sum_{|x|> \varepsilon n} p_n(e,x)
 \le \sum_{r>\varepsilon n} \Big(\sup_{|x|=r} p_n(e,x)\Big)|S(e,r)| \\
&\le \sum_{r>\varepsilon n} C\exp\Big(-c\frac{r^\alpha}{n^\beta}\Big)|B(e,r)|.
\end{align*}
Fix $M>1$ and split the sum as $I_1+I_2$ with
$$
I_1=\sum_{\varepsilon n < r \le Mn} Ce^{-c r^\alpha/n^\beta}|B(e,r)|,\qquad
I_2=\sum_{r>Mn} Ce^{-c r^\alpha/n^\beta}|B(e,r)|.
$$

By subexponential growth, for any $\delta>0$ and all sufficiently large $n$,
$
|B(e,r)|\le \exp(\delta r)\le \exp(\delta M n)
$
for every $\varepsilon n\le r\le Mn$. Therefore,
$$
I_1 \le (Mn)C\exp\Big(-c\frac{(\varepsilon n)^\alpha}{n^\beta}\Big)\exp(\delta M n)
=(Mn)C\exp\Big(-c\varepsilon^\alpha n^{\alpha-\beta} + \delta M n\Big).
$$
If $\alpha-\beta>1$, then $-c \varepsilon^\alpha n^{\alpha-\beta}+\delta M n\to -\infty$ for any fixed $\delta>0$, hence $I_1\to0$. 
For $r\ge Mn$ we have
$$
\frac{r^\alpha}{n^\beta} \ge \frac{(Mn)^{\alpha-1}}{n^\beta}r = M^{\alpha-1}n^{\alpha-1-\beta}r.
$$
With this fixed $\delta>0$, take $n$ large so that $|B(e,r)|\le e^{\delta r}$ for all $r\ge Mn$. Then
$$
-c\frac{r^\alpha}{n^\beta} + \log|B(e,r)|
\le -\big(c M^{\alpha-1}n^{\alpha-1-\beta}-\delta\big)r.
$$
Since $\alpha-\beta>1$, $c M^{\alpha-1}n^{\alpha-1-\beta}\to\infty$ so for large $n$ the coefficient is $\ge \delta$, giving
$
I_2\le \sum_{r>Mn} Ce^{-\delta r}\le C'e^{-\delta Mn}\to 0.
$

If $\alpha-\beta=1$, choose $\delta>0$ so small that $\delta M < \tfrac12 c \varepsilon^\alpha$ and $\delta<c M^{\alpha-1}$; then $I_1\le (Mn)C\exp(-\tfrac12 c\varepsilon^\alpha n)\to0$. Also, $c M^{\alpha-1}n^{\alpha-1-\beta}=c M^{\alpha-1}$ is constant and strictly larger than $\delta$ by our choice of $\delta$, yielding again $I_2\le C'e^{-(c M^{\alpha-1}-\delta)Mn}\to0$.

Combining the two cases completes the proof.
\end{proof}

\begin{remark}[On the threshold]
The condition $\alpha\ge 1+\beta$ is optimal at this level of generality: for groups of subexponential growth arbitrarily close to exponential (e.g.\ of intermediate growth), one needs a negative exponent at least linear in $r$ at the scale $r\asymp n$ to dominate $\log|B(e,r)|=o(r)$ uniformly.
\end{remark}

\begin{prop}[Polynomial volume: $\alpha>\beta$ suffices]\label{prop:poly}
Assume $|B(e,r)|\le C_0 r^{D}$ for some $D\ge0$ and all $r\ge1$, and that \eqref{HK} holds with $\alpha>\beta$. Then, for every $\varepsilon>0$,
$$
\PP\left(\frac{|X_n|}{n}>\varepsilon\right)\to 0 \text{ as } n\to\infty.
$$
\end{prop}

\begin{proof}
Using $|S(e,r)|\le C_1 r^{D-1}$ and \eqref{HK},
$$
\PP\left(\tfrac{|X_n|}{n}>\varepsilon\right)
\le C \sum_{r>\varepsilon n} r^{D-1}e^{-c r^\alpha/n^\beta}
\le C \int_{\varepsilon n-1}^{\infty} r^{D-1}e^{-c r^\alpha/n^\beta}dr.
$$
Make the change of variables $u=r^\alpha/n^\beta$; then $r=(un^\beta)^{1/\alpha}$ and $dr=\frac{n^{\beta/\alpha}}{\alpha} u^{1/\alpha-1}du$. Therefore, for large $n$
\begin{align*}
\int_{\varepsilon n-1}^{\infty} r^{D-1}e^{-c r^\alpha/n^\beta}dr
&\lesssim n^{\beta D/\alpha}\int_{A(n)}^{\infty} u^{D/\alpha - 1} e^{-c u}du \\
&\lesssim n^{\beta D/\alpha} \exp\big(-cA(n)\big),
\end{align*}
where $A(n)=\frac{(\varepsilon n-1)^\alpha}{n^\beta}$ and the last bound uses the standard tail estimate $\int_R^\infty u^{\theta-1} e^{-c u}du \le C'e^{-cR}$ for some $C'>0$, for large $R$ and $\theta>0$ (for some possibly smaller $c>0$). Since $\alpha>\beta$, the right-hand side tends to $0$.
\end{proof}

We now show that on nilpotent groups, the speed can go to zero in probability even when the first moment is infinite. This complements the discussion at the beginning of this subsection.

\begin{lemma}[Mal'cev coordinates ann associated estimates]\label{lem:coords}
Let $\Gamma$ be a finitely generated torsion-free nilpotent group of step $c$, endowed with a word
metric $|\cdot|$ coming from a finite symmetric generating set.
Let $N$ be its Mal'cev completion, a connected simply connected nilpotent Lie group, and let $\mathfrak n$ be the Lie algebra of $N$.
Fix a vector-space decomposition
$$
   \mathfrak n = \mathfrak n_1 \oplus \mathfrak n_2 \oplus \cdots \oplus \mathfrak n_c
$$
adapted to the lower central series (so that $[\mathfrak n_i,\mathfrak n_j]\subset \mathfrak n_{\ge i+j}$). Fix a Mal'cev basis adapted to this decomposition and write, for each $g\in\Gamma$,
$$
   g = \exp\Big(\sum_{k=1}^c x^{(k)}(g)\Big),
   \qquad x^{(k)}(g)\in \mathfrak n_k.
$$
Then there exists a constant $C\ge1$ such that for all $g\in\Gamma$ and all $1\le k\le c$,
\begin{equation}\label{eq:coord-growth}
   \|x^{(k)}(g)\|_1 \le C|g|^k.
\end{equation}
\end{lemma}

\begin{proof}
Define 
$$
N(g) := \sum_{k=1}^c \|x^{(k)}(g)\|_1^{1/k} \quad g\in\Gamma.
$$
By construction, for each $k$,
\begin{equation}\label{eq:coord-bound}
   \|x^{(k)}(g)\|_1 \le N(g)^k.
\end{equation}

It is a classical theorem of Guivarc'h \cite{Guivarch1970} (see also \cite{Bass1972}, \cite[Theorem 2.7, Remark 2.10]{Breuillard2014}) that
word length $|\cdot|$ is coarsely bi-Lipschitz equivalent to $N$.
That is, there exist constants $A,B\ge 1$ such that
\begin{equation}\label{eq:equivalence}
   A^{-1} N(g) - B \le |g| \le A N(g) + B,
   \qquad \forall g\in \Gamma.
\end{equation}

From \eqref{eq:equivalence} we obtain $N(g) \le A(|g| + B) \le A'|g|$
for some constant $A'\ge 1$ depending only on $A$, $B$.

Substituting this into \eqref{eq:coord-bound} gives for $1\le k\le c$,
$$
   \|x^{(k)}(g)\|_1 \le N(g)^k
   \le (A'|g|)^k.
$$
Taking $C = (A')^c$ we obtain the desired inequality
$$
   \|x^{(k)}(g)\|_1 \le C |g|^k
$$
for all $g\in \Gamma$ and for $1\le k\le c$.
\end{proof}

\begin{theorem}\label{thm:nilpotent-wlln}
Let $\Gamma$ be a finitely generated torsion-free nilpotent group equipped with a word metric $|\cdot|$ associated with a fixed finite symmetric generating set. Let $\mu$ be a symmetric probability measure on $\Gamma$. Let $(S_k)_{k\ge1}$ be i.i.d.\ $\Gamma$-valued random variables with law $\mu$, and let $X_n = S_1 S_2 \cdots S_n$ be the associated random walk.

Assume there exist constants $r_0 \ge 3$ and $0 < c_1 \le c_2 < \infty$ such that for all integers $r \ge r_0$, the probability mass $p_r := \sum_{|g|=r} \mu(g)$ satisfies
\begin{equation}\label{eq:shell-condition}
  \frac{c_1}{r^2 \log r} \le p_r \le \frac{c_2}{r^2 \log r}.
\end{equation}
Then $\E|S_1| = \infty$, but
$$
  \frac{|X_n|}{n} \to 0 \text{ as } n\to\infty \text{ in probability}.
$$
\end{theorem}

\begin{lemma}\label{lem:shell-consequences}
Under the shell condition \eqref{eq:shell-condition}:
\begin{enumerate}
\item $\E|S_1| = \infty$.
\item There exists a constant $C>0$ such that for all $x\ge r_0$,
\begin{equation}\label{eq:tail-bound}
\PP(|S_1|>x) \le \frac{C}{x\log x}.
\end{equation}
In particular, $\lim_{x\to\infty} x\PP(|S_1|>x) = 0$.
\item For every integer $m \ge 2$, there exist constants $R_0 = R_0(r_0,m) \ge r_0$ and $C_m>0$ such that for all $R\ge R_0$,
\begin{equation}\label{eq:truncated-moment-S}
\E\bigl[|S_1|^m \Ind_{\{|S_1|\le R\}}\bigr]
\le C_m \frac{R^{m-1}}{\log R}.
\end{equation}
\end{enumerate}
\end{lemma}

\begin{proof}
(1) Using the lower bound, $\E|S_1| = \sum_{r\ge 1} r p_r \ge c_1 \sum_{r\ge r_0} \frac{1}{r \log r} = \infty$.

(2) For $x\ge r_0$, $\PP(|S_1|>x) = \sum_{r>x} p_r \le c_2 \sum_{r>x} \frac{1}{r^2\log r}$. Comparing the sum to the integral $\int_x^\infty \frac{dt}{t^2 \log t}$ yields the bound \eqref{eq:tail-bound}.

(3) For $m \ge 2$,
$$
  \E\bigl[|S_1|^m \Ind_{\{|S_1|\le R\}}\bigr] \le \sum_{r=1}^{r_0-1}r^mp_r +c_2 \sum_{r=r_0}^R \frac{r^m}{r^2 \log r}.
$$
If $m>2$, we can compare $\sum_{r=r_0}^R \frac{r^m}{r^2 \log r}$ to the integral $\int_{r_0}^{R+1}\frac{t^{m-2}}{\log t} dt$, and if $m=2$, we compare it to the integral $\int_{r_0-1}^{R}\frac{1}{\log t} dt$.
Using integration by parts and standard estimation techniques involving integrals we get the desired upperbound.
\end{proof}

\begin{proof}[Proof of Theorem \ref{thm:nilpotent-wlln}]
Fix $n \in \mathbb{N}$ and define the truncated increments $\bar{S}_k^{(n)}$ (for $k\in \NN$) by
$$
  \bar{S}_k^{(n)} :=
  \begin{cases}
    S_k, & |S_k|\le n,\\
    e,   & |S_k|>n.
  \end{cases}
$$
Let $\bar{X}_k^{(n)} := \bar{S}_1^{(n)} \cdots \bar{S}_k^{(n)}$.

The coupling probability is controlled by the tail bound:
\begin{equation}\label{eq:coupling-estimate}
  \PP(X_n \neq \bar X_n^{(n)})
  \le \sum_{k=1}^n \PP(|S_k|>n)
  = n \PP(|S_1|>n)
  \le \frac{C}{\log n} \to 0 \text{ as } n\to\infty.
\end{equation}
Thus, it suffices to prove convergence for the truncated walk $\bar X_n^{(n)}$.
Since $\mu$ is symmetric and the truncation event $\{|g|\le n\}$ is invariant under inversion, the distribution of $\bar S_k^{(n)}$ is symmetric.

We embed $\Gamma$ as a cocompact lattice in a simply connected nilpotent Lie group $G$ (Malcev completion). We identify $G$ with its Lie algebra $\mathfrak{g}\cong\mathbb{R}^d$ via exponential coordinates. For $g\in G$, let $x(g)$ denote the image of $g$ in $\RR^d$ via exponential coordinates, and $x_j(g)$ denote the $j$-th coordinate of $g$. Fix a vector-space decomposition
$$
   \mathfrak g = \mathfrak g_1 \oplus \mathfrak g_2 \oplus \cdots \oplus \mathfrak g_c
$$
adapted to the lower central series, so that $[\mathfrak n_i,\mathfrak n_j]\subset \mathfrak n_{\ge i+j}$. This yields integer weights $w_j\in\mathbb{Z}_{\ge1}$ for coordinates $x_j$ ($1\le j\le d$), with the property that the Baker-Campbell-Hausdorff product is polynomial and triangular with respect to these weights.

By \eqref{eq:coord-growth} there is $C\ge1$ such that for all $g\in\Gamma$ and each $j$,
\begin{equation}\label{eq:coord-growth-2}
  |x_j(g)| \le C|g|^{w_j}.
\end{equation}
To prove $|\bar X_n^{(n)}|/n \to 0$ in probability, it suffices to show that for every coordinate $j$,
$$
  \frac{x_j(\bar X_n^{(n)})}{n^{w_j}} \to 0 \text{ as } n\to\infty \text{ in probability}.
$$
Let $\xi_k^{(n)} = x(\bar S_k^{(n)}) \in \mathbb{R}^d$ and $\xi_{k,j}^{(n)}=x_j(\bar{S^{(n)}_k})$. By the symmetry of $\bar S_k^{(n)}$, we have $\xi_k^{(n)} = -\xi_k^{(n)}$ in distribution, hence $\E[\xi_k^{(n)}] = 0$.
Applying \eqref{eq:coord-growth-2} and Lemma \ref{lem:shell-consequences}(3) with $R=n$, for any even integer $m \ge 2$ the moments of the coordinates satisfy:
\begin{equation}\label{eq:xi-moment}
  \E[|\xi_{k,j}^{(n)}|^m]
  \le C \E\bigl[|S_1|^{m w_j}\Ind_{\{|S_1|\le n\}}\bigr]
  \le C_{j,m} \frac{n^{m w_j - 1}}{\log n},
\end{equation}
for all $n$ sufficiently large (depending on $j,m$) and for some $C_{j,m}>0$ depending on $j$, $m$.

We prove convergence by establishing bounds on all even moments of the coordinates.

\begin{prop}\label{prop:inductive-moments}
For every coordinate $j$ with weight $w_j$ and every even integer $m \ge 2$, there exist constants $C_{j,m}>0$ and $N_{j,m}\in\mathbb{N}$ such that for all $n \ge N_{j,m}$,
\begin{equation}\label{eq:strong-hypothesis}
  \sup_{1 \le k \le n} \E[|x_j(\bar X_k^{(n)})|^m] \le C_{j,m} \frac{n^{m w_j}}{\log n}.
\end{equation}
\end{prop}

\begin{proof}
We proceed by induction on the weight $w_j$. First, let us prove the following useful lemma:

\begin{lemma}[Polynomial moment lemma]\label{lem:polynomial-growth}
Suppose that for every coordinate $i$ with weight $< w$ and every even integer $M\ge2$, there exist constants $C_{i,M}$ and $N_{i,M}$ such that \eqref{eq:strong-hypothesis} holds (for that coordinate and moment order) for all $n\ge N_{i,M}$.
Let $P(u)$ be a monomial in coordinates of weights $< w$, say $P(u) = \prod_{\ell=1}^r u_{i_\ell}^{\alpha_\ell}$, with total homogeneous weight $D = \sum_{\ell=1}^r \alpha_\ell w_{i_\ell}$.
Then for any even integer $m \ge 2$ there exist $C>0$ and $N\in\mathbb{N}$ (depending on $P$) such that for all $n\ge N$,
$$
  \sup_{1 \le k \le n} \E[|P(x(\bar X_k^{(n)}))|^m] \le C \frac{n^{m D}}{\log n}.
$$   
\end{lemma}

\begin{proof}
Apply the generalised H\"older inequality with equal exponents $p_\ell = r$ for all $\ell=1,\dots,r$ (so $\sum_{\ell=1}^r 1/p_\ell = 1$):
$$
  \E[|P(x(\bar X_k^{(n)}))|^m]
  = \E\left[ \prod_{\ell=1}^r |x_{i_\ell}(\bar X_k^{(n)})|^{m \alpha_\ell} \right]
  \le \prod_{\ell=1}^r \left( \E\left[ |x_{i_\ell}(\bar X_k^{(n)})|^{r m \alpha_\ell} \right] \right)^{1/r}.
$$
The exponent $M_\ell = r m \alpha_\ell$ is an even integer (since $m$ is even). By our assumption, for each $\ell$ there exist $C_\ell$ and $N_\ell$ such that for all $n\ge N_\ell$,
$$
  \sup_{1\le k\le n} \E\left[|x_{i_\ell}(\bar X_k^{(n)})|^{M_\ell}\right]
  \le C_\ell \frac{n^{M_\ell w_{i_\ell}}}{\log n}.
$$
Let $N:=\max_{1\le\ell \le r} N_\ell$ and $C:=\max_{1\le\ell\le r} C_\ell$. Then, for all $n\ge N$ and $1\le k \le n$,
\begin{align*}
  \E[|P(x(\bar X_k^{(n)}))|^m]
  &\le \prod_{\ell=1}^r \left( C \frac{n^{(r m \alpha_\ell) w_{i_\ell}}}{\log n} \right)^{1/r}
   = C' \frac{n^{\sum_\ell m \alpha_\ell w_{i_\ell}}}{(\log n)^{\sum_\ell 1/r}}
   = C' \frac{n^{m D}}{\log n}.
\end{align*}
where $C'=C^r$.
\end{proof}
For a coordinate $x_j$ of weight $1$, the group law is linear in that coordinate:
$$
  x_j(\bar X_k^{(n)}) = \sum_{i=1}^k \xi_{i,j}^{(n)}.
$$
By Rosenthal's inequality \cite{Rosenthal1970}, for even $m \ge 2$,
$$
  \E\left[ \left| \sum_{i=1}^k \xi_{i,j}^{(n)} \right|^m \right]
  \le C_m \left( \sum_{i=1}^k \E\left[|\xi_{i,j}^{(n)}|^m\right]
  + \left( \sum_{i=1}^k \E\left[|\xi_{i,j}^{(n)}|^2\right] \right)^{m/2} \right).
$$
Using \eqref{eq:xi-moment} with $w_j = 1$ and $1\le k \le n$ (for $n$ sufficiently large), there exists $C''>0$ depending on $j$, $m$ such that
$$
\sum_{i=1}^k \E|\xi_{i,j}^{(n)}|^m \le C'' \frac{n^m}{\log n},
\qquad
\left(\sum_{i=1}^k \E|\xi_{i,j}^{(n)}|^2\right)^{m/2}
\le \left(C''\frac{n^2}{\log n}\right)^{m/2}
\le C'' \frac{n^m}{\log n}.
$$
Therefore,
$$
  \sup_{1\le k\le n}\E|x_j(\bar X_k^{(n)})|^m \le C'' \frac{n^m}{\log n},
$$
which is \eqref{eq:strong-hypothesis} for $w_j=1$.

\smallskip\noindent
\emph{Inductive Step:}
Assume \eqref{eq:strong-hypothesis} holds for all coordinates of weight $< w$ and for all even moment orders. Let $x_j$ be a coordinate of weight $w_j = w$, and fix $n$ sufficiently large.

For each $k\ge1$ we write
$$
  \bar X_k^{(n)} = \bar X_{k-1}^{(n)} \bar S_k^{(n)}. 
$$
By the BCH formula (in these weighted coordinates), there exist finitely many pairs of monomials
$$
  \{(P_\alpha,Q_\alpha)\}_{\alpha\in A_j}
$$
and constants $c_\alpha$ such that for all $g,h\in G$,
\begin{equation}\label{eq:BCH-coordinate}
  x_j(gh)
  = x_j(g) + x_j(h)
    + \Phi_j(x(g),x(h)),
  \qquad
  \Phi_j(u,v) := \sum_{\alpha\in A_j} c_\alpha P_\alpha(u) Q_\alpha(v),
\end{equation}
where each term corresponds to a single Baker-Campbell-Hausdorff monomial contributing to the $j$th coordinate, and every such monomial has total homogeneous weight $w$:
\begin{equation}\label{eq:weight-balance}
  w(P_\alpha) + w(Q_\alpha) = w,
  \qquad
  w(P_\alpha),w(Q_\alpha)\ge 1.
\end{equation}
Moreover, every coordinate variable appearing in $P_\alpha$ or $Q_\alpha$ has weight strictly less than $w$ (equivalently, $P_\alpha$ and $Q_\alpha$ only involve lower-weight coordinates). This follows because any BCH correction monomial is at least bilinear in $(u,v)$, and if a variable of weight $w$ were to appear then the total weight would exceed $w$.

Applying \eqref{eq:BCH-coordinate} with $g=\bar X_{k-1}^{(n)}$ and $h=\bar S_k^{(n)}$ gives
$$
  x_j(\bar X_k^{(n)})
  = x_j(\bar X_{k-1}^{(n)}) + x_j(\bar S_k^{(n)})
    + \Phi_j(x(\bar X_{k-1}^{(n)}), x(\bar S_k^{(n)}))
  = x_j(\bar X_{k-1}^{(n)}) + \xi_{k,j}^{(n)}
    + \Phi_j(x(\bar X_{k-1}^{(n)}), \xi_k^{(n)}).
$$
Consider $K\in\{ 1,2,\dots, n\}$. Summing from $k=1$ to $K$ (with $\bar X_0^{(n)} = e$) yields
$$
  x_j(\bar X_K^{(n)})
  = \sum_{k=1}^K \xi_{k,j}^{(n)}
    + \sum_{k=1}^K \Phi_j(x(\bar X_{k-1}^{(n)}),\xi_k^{(n)}).
$$
Define
$$
  L_K := \sum_{k=1}^K \xi_{k,j}^{(n)},
  \qquad
  R_K := \sum_{k=1}^K \Phi_j(x(\bar X_{k-1}^{(n)}),\xi_k^{(n)}),
$$
so that $x_j(\bar X_K^{(n)}) = L_K + R_K$.

\medskip\noindent
\emph{Estimate for $L_K$.}
Using \eqref{eq:xi-moment} with $w_j=w$ and Rosenthal's inequality, for every even $m\ge2$ there exists $C_{j,m}>0$ such that
\begin{equation}\label{eq:L_K-estimate}
  \sup_{1\le K\le n} \E|L_K|^m \le C_{j,m} \frac{n^{mw}}{\log n}.
\end{equation}

\medskip\noindent
We now decompose $R_K$ suitably. Let $(\mathcal{F}_k)_{k\ge0}$ be the natural filtration
$$
  \mathcal{F}_k := \sigma(\bar S_1^{(n)},\dots,\bar S_k^{(n)}) =\sigma(\xi_1^{(n)},\dots,\xi_k^{(n)}),\quad \mathcal{F}_0 \text{ trivial}.
$$
Then for each $k$, $x(\bar X_{k-1}^{(n)})$ (and hence each $P_\alpha(x(\bar X_{k-1}^{(n)}))$) is $\mathcal{F}_{k-1}$-measurable, while $\xi_k^{(n)}$ (and thus each $Q_\alpha(\xi_k^{(n)})$) is independent of $\mathcal{F}_{k-1}$.

Define
\begin{align*}
  d_k &:= \Phi_j(x(\bar X_{k-1}^{(n)}),\xi_k^{(n)})
       = \sum_{\alpha\in A_j}c_\alpha P_\alpha(x(\bar X_{k-1}^{(n)})) Q_\alpha(\xi_k^{(n)}), \\
  a_k &:= \E[d_k \mid \mathcal{F}_{k-1}]
       = \sum_{\alpha\in A_j}
         c_\alpha P_\alpha(x(\bar X_{k-1}^{(n)}))
                  \E[Q_\alpha(\xi_1^{(n)})],\\
  \tilde d_k &:= d_k - a_k
              = \sum_{\alpha\in A_j}
                c_\alpha P_\alpha(x(\bar X_{k-1}^{(n)}))
                \Bigl(Q_\alpha(\xi_k^{(n)}) - \E[ Q_\alpha(\xi_1^{(n)})]\Bigr).
\end{align*}
Set
$$
  D_K := \sum_{k=1}^K a_k,
  \qquad
  M_K := \sum_{k=1}^K \tilde d_k,
$$
so that $R_K = D_K + M_K$.

By construction, $\E[\tilde d_k \mid \mathcal{F}_{k-1}]=0$ and $\tilde d_k$ is $\mathcal{F}_k$-measurable, hence $(M_K)_{K\ge0}$ is a martingale with respect to $(\mathcal{F}_K)_{K\ge0}$.

\medskip\noindent
\emph{Estimating $D_K$.}
For each $\alpha\in A_j$, the monomial $Q_\alpha$ has either odd or even total degree.

\emph{Case 1: $Q_\alpha$ has odd degree.}
By symmetry, $\xi_1^{(n)} = -\xi_1^{(n)}$ in distribution, so for any integrable $f$ we have $\E[ f(\xi_1^{(n)})]=\E [f(-\xi_1^{(n)})]$. Since $Q_\alpha$ is odd, $Q_\alpha(-v)=-Q_\alpha(v)$ and therefore $\E[Q_\alpha(\xi_1^{(n)})]=0$. Hence, such terms do not contribute to $D_K$.

\emph{Case 2: $Q_\alpha$ has even degree $d\ge2$.}
Write $Q_\alpha(v)=\prod_i v_{j_i}^{\beta_i}$ with $\sum_i \beta_i=d\ge2$. Since $w_{j_i}\ge1$,
$$
  w(Q_\alpha)=\sum_i \beta_i w_{j_i} \ge \sum_i \beta_i=d\ge2.
$$
Using \eqref{eq:coord-growth} and Lemma \ref{lem:shell-consequences}(3) (with exponent $w(Q_\alpha)$),
\begin{align*}
\bigl|\E[Q_\alpha(\xi_1^{(n)})]\bigr|\le\E\bigl[|Q_\alpha(\xi_1^{(n)})|\bigr]\le C'_{j,m}\E\bigl[|S_1|^{w(Q_\alpha)}\Ind_{\{|S_1|}\le n\}\bigr]\le C'_{j,m} \frac{n^{w(Q_\alpha)-1}}{\log n}.
\end{align*}
Now estimate $\|D_K\|_m$ for even $m\ge2$. By Minkowski and the preceding bound,
\begin{align*}
  \|D_K\|_m
  &= \left\|\sum_{k=1}^K a_k\right\|_m
  \le \sum_{k=1}^K \|a_k\|_m \\
  &\le \sum_{k=1}^K \sum_{\alpha\in A_j} |c_\alpha|
      \bigl|\E[Q_\alpha(\xi_1^{(n)})]\bigr|\|P_\alpha(x(\bar X_{k-1}^{(n)}))\|_m \\
  &\le C'_{j,m} \sum_{\alpha\in A_j} K \cdot \frac{n^{w(Q_\alpha)-1}}{\log n}
      \cdot \sup_{1\le k\le n}\|P_\alpha(x(\bar X_{k-1}^{(n)}))\|_m.
\end{align*}
By Lemma \ref{lem:polynomial-growth} (since $w(P_\alpha)<w$) we have
$$
  \sup_{1\le k\le n}\|P_\alpha(x(\bar X_{k-1}^{(n)}))\|_m
  \le \left(L\frac{n^{m w(P_\alpha)}}{\log n}\right)^{1/m}
  = L' \frac{n^{w(P_\alpha)}}{(\log n)^{1/m}}.
$$
Since $1\le K\le n$ and $w(P_\alpha)+w(Q_\alpha)=w$, this yields
$$
  \|D_K\|_m \le C'_{j,m} \sum_{\alpha\in A_j} \frac{n^{w}}{(\log n)^{1+1/m}}
  \le C''_{j,m} \frac{n^{w}}{(\log n)^{1+1/m}}.
$$
where $C''_{j,m}=C'_{j,m}|A_j|$. Raising to power $m$,
\begin{equation}\label{eq:D_K-estimate}
  \sup_{1\le K\le n}\E|D_K|^m \le C'''_{j,m} \frac{n^{mw}}{(\log n)^{m+1}}
  \le C'''_{j,m} \frac{n^{mw}}{\log n}, \quad C'''_{j,m}=(C''_{j,m})^m.
\end{equation}

\medskip\noindent
\emph{Estimating $M_K$.}
We use the Burkholder-Rosenthal inequality (see \cite{Burkholder1973}) for the martingale $M_K = \sum_{k=1}^K \tilde d_k$. For even $m \ge 2$, there exists $C'_m>0$ such that
\begin{equation}\label{eq:Burkholder-Rosenthal}
  \E|M_K|^m \le C'_m \left( \sum_{k=1}^K \E|\tilde d_k|^m
  + \E\left[ \left( \sum_{k=1}^K \E[\tilde d_k^2 \mid \mathcal{F}_{k-1}] \right)^{m/2} \right] \right).
\end{equation}

Let $\tilde{Q}_\alpha(\xi_k^{(n)}) := Q_\alpha(\xi_k^{(n)}) - \E [Q_\alpha(\xi_1^{(n)})]$.
By the triangle inequality and Jensen's inequality, for every integer $p\ge 1$, we have
$$
  \E\bigl[|\tilde{Q}_\alpha(\xi_k^{(n)})|^p\bigr]
  \le 2^p \E\bigl[|Q_\alpha(\xi_k^{(n)})|^p\bigr].
$$
In particular for $p=m$ (with $m$ even), using \eqref{eq:coord-growth} and Lemma \ref{lem:shell-consequences}(3),
$$
  \E\bigl[|\tilde{Q}_\alpha(\xi_k^{(n)})|^p\bigr]
  \le C\E\bigl[|S_1|^{p w(Q_\alpha)}\Ind_{\{|S_1|\le n\}}\bigr]
  \le C \frac{n^{p w(Q_\alpha)-1}}{\log n},
$$
for all $n$ sufficiently large (depending on $\alpha,m$).

Let $N_j:=|A_j|<\infty$. Using the inequality
$\bigl|\sum_{\alpha\in A_j} z_\alpha\bigr|^m \le N_j^{m-1}\sum_{\alpha\in A_j} |z_\alpha|^m$,
we obtain
\begin{align*}
  \E|\tilde d_k|^m
  &= \E\left|\sum_{\alpha\in A_j} c_\alpha P_\alpha(x(\bar X_{k-1}^{(n)}))\tilde Q_\alpha(\xi_k^{(n)})\right|^m \\
  &\le N_j^{m-1}\sum_{\alpha\in A_j}
      \E\left[ |P_\alpha(x(\bar X_{k-1}^{(n)}))|^m |\tilde Q_\alpha(\xi_k^{(n)})|^m \right].
\end{align*}
Since $P_\alpha(x(\bar X_{k-1}^{(n)}))$ is $\mathcal F_{k-1}$-measurable and $\xi_k^{(n)}$ is independent of $\mathcal F_{k-1}$, the factors are independent, hence
$$
  \E\left[ |P_\alpha(x(\bar X_{k-1}^{(n)}))|^m|\tilde Q_\alpha(\xi_k^{(n)})|^m \right]
  = \E[|P_\alpha(x(\bar X_{k-1}^{(n)}))|^m]\cdot \E[|\tilde Q_\alpha(\xi_1^{(n)})|^m].
$$
Using Lemma \ref{lem:polynomial-growth} for $P_\alpha$  and the bound above for $\tilde Q_\alpha$ with $p=m$, there exists $L'>0$ such that
$$
  \sup_{1\le k\le n}\E[|P_\alpha(x(\bar X_{k-1}^{(n)}))|^m] \le L'\frac{n^{m w(P_\alpha)}}{\log n},
  \qquad
  \E[|\tilde Q_\alpha(\xi_1^{(n)})|^m] \le L'\frac{n^{m w(Q_\alpha)-1}}{\log n}.
$$
Therefore,
\begin{align}\label{eq:M_K-subestimate1}
  \sum_{k=1}^K \E|\tilde d_k|^m
  &\le N_j^{m-1}L'^2\sum_{\alpha\in A_j}\sum_{k=1}^K
      \left(\frac{n^{m w(P_\alpha)}}{\log n}\right)\left(\frac{n^{m w(Q_\alpha)-1}}{\log n}\right) \\
  &\le N_j^{m-1}L'^2\sum_{\alpha\in A_j} K \cdot \frac{n^{m(w(P_\alpha)+w(Q_\alpha))-1}}{(\log n)^2}
  \le N_j^{m}L'^2 \frac{n^{mw}}{(\log n)^2},
\end{align}
since $1\le K\le n$ and $w(P_\alpha)+w(Q_\alpha)=w$.

Define
$$
  Z_K := \sum_{k=1}^K \E[\tilde d_k^2\mid \mathcal F_{k-1}].
$$
Set $M_j := \max_{\alpha\in A_j}|c_\alpha|$. By Cauchy-Schwarz inequality and the finiteness of $A_j$,
\begin{align*}
  \E[\tilde d_k^2\mid \mathcal F_{k-1}]
  &= \E\left[\left(\sum_{\alpha\in A_j} c_\alpha P_\alpha(x(\bar X_{k-1}^{(n)}))\tilde Q_\alpha(\xi_k^{(n)})\right)^2 \Big| \mathcal F_{k-1}\right] \\
  &\le N_j M_j^2 \sum_{\alpha\in A_j} P_\alpha(x(\bar X_{k-1}^{(n)}))^2
       \E\left[\tilde Q_\alpha(\xi_1^{(n)})^2\right].
\end{align*}
Let
$$
  V_\alpha^{(n)} := \E\left[\tilde Q_\alpha(\xi_1^{(n)})^2\right]
  \le C \frac{n^{2w(Q_\alpha)-1}}{\log n}.
$$
Then
$$
  Z_K \le C \sum_{\alpha\in A_j} V_\alpha^{(n)} \sum_{k=1}^K P_\alpha(x(\bar X_{k-1}^{(n)}))^2.
$$

We estimate $\|Z_K\|_{m/2}$ (note $m/2\ge1$). By Minkowski's inequality in $L^{m/2}$,
\begin{align*}
  \|Z_K\|_{m/2}
  &\le C \sum_{\alpha\in A_j} V_\alpha^{(n)}
        \left\|\sum_{k=1}^K P_\alpha(x(\bar X_{k-1}^{(n)}))^2\right\|_{m/2} \\
  &\le C \sum_{\alpha\in A_j} V_\alpha^{(n)} \sum_{k=1}^K \|P_\alpha(x(\bar X_{k-1}^{(n)}))^2\|_{m/2}.
\end{align*}
But $\|P_\alpha(\cdot)^2\|_{m/2} = \|P_\alpha(\cdot)\|_m^2$, and by Lemma \ref{lem:polynomial-growth},
$$
  \|P_\alpha(x(\bar X_{k-1}^{(n)}))\|_m
  \le \left(C\frac{n^{m w(P_\alpha)}}{\log n}\right)^{1/m}
  = C \frac{n^{w(P_\alpha)}}{(\log n)^{1/m}}.
$$
Therefore,
$$
  \sum_{k=1}^K \|P_\alpha(x(\bar X_{k-1}^{(n)}))^2\|_{m/2}
  \le K \cdot C \frac{n^{2w(P_\alpha)}}{(\log n)^{2/m}}
  \le C \frac{n^{2w(P_\alpha)+1}}{(\log n)^{2/m}}.
$$
Combining with the bound on $V_\alpha^{(n)}$ yields
\begin{align*}
  \|Z_K\|_{m/2}
  &\le C \sum_{\alpha\in A_j}
     \left(\frac{n^{2w(Q_\alpha)-1}}{\log n}\right)
     \left(\frac{n^{2w(P_\alpha)+1}}{(\log n)^{2/m}}\right) \\
  &= C \sum_{\alpha\in A_j} \frac{n^{2(w(P_\alpha)+w(Q_\alpha))}}{(\log n)^{1+2/m}}
  \le C \frac{n^{2w}}{(\log n)^{1+2/m}}.
\end{align*}
Hence
\begin{equation}\label{eq:M_K-subestimate2}
  \E\left[Z_K^{m/2}\right] = \|Z_K\|_{m/2}^{m/2}
  \le C \frac{n^{mw}}{(\log n)^{1+m/2}}.
\end{equation}

Putting \eqref{eq:M_K-subestimate1} and \eqref{eq:M_K-subestimate2} into \eqref{eq:Burkholder-Rosenthal} gives
\begin{equation}\label{eq:M_K-estimate}
  \sup_{1\le K\le n}\E|M_K|^m
  \le C_m\left(\frac{n^{mw}}{(\log n)^2} + \frac{n^{mw}}{(\log n)^{1+m/2}}\right)
  \le C \frac{n^{mw}}{(\log n)^2}.
\end{equation}

\medskip\noindent
From \eqref{eq:L_K-estimate}, \eqref{eq:D_K-estimate} and \eqref{eq:M_K-estimate}, have a large enough $C_1>0$ such that
\begin{align*}
  \sup_{1 \le K \le n} \E|L_K|^m &\le C_1 \frac{n^{mw}}{\log n}, \\
  \sup_{1 \le K \le n} \E|D_K|^m &\le C_1 \frac{n^{mw}}{(\log n)^{m+1}}, \\
  \sup_{1 \le K \le n} \E|M_K|^m &\le C_1 \frac{n^{mw}}{(\log n)^2}.
\end{align*}
Since $x_j(\bar X_K^{(n)}) = L_K + D_K + M_K$, taking $m$-th roots and using Minkowski gives
\begin{align*}
  \|x_j(\bar X_K^{(n)})\|_m
  &\le \|L_K\|_m + \|D_K\|_m + \|M_K\|_m \\
  &\le C_1^{1/m} n^{w} (\log n)^{-1/m}
     + C_1^{1/m} n^{w} (\log n)^{-1-1/m}
     + C_1^{1/m} n^{w} (\log n)^{-2/m}
  \le 3C_1^{1/m} n^{w} (\log n)^{-1/m}.
\end{align*}
Raising to the $m$-th power yields
$$
  \sup_{1 \le K \le n} \E|x_j(\bar X_K^{(n)})|^m \le 3^mC_1 \frac{n^{mw}}{\log n}.
$$
This completes the inductive step and hence the proof of Proposition \ref{prop:inductive-moments}.
\end{proof}

By Proposition \ref{prop:inductive-moments} with $m=2$, for every coordinate $j$ there exist $C>0$ and $N \in \NN$ (depending on $j$) such that for all $n\ge N$,
$$
  \PP\left( \frac{|x_j(\bar X_n^{(n)})|}{n^{w_j}} > \epsilon \right)
  \le \frac{\E[|x_j(\bar X_n^{(n)})|^2]}{\epsilon^2 n^{2w_j}}
  \le \frac{C}{\epsilon^2 \log n} \to 0 \text{ as } n\to\infty.
$$
This also implies $\frac{|x_j(\bar X_n^{(n)})|^{1/w_j}}{n}\to0$ in probability for each $j$, and by \eqref{eq:equivalence}
$$
  \frac{|\bar X_n^{(n)}|}{n} \to 0 \text{ in probability}.
$$
Finally, using \eqref{eq:coupling-estimate} we get
$$
  \frac{|X_n|}{n} \to 0 \text{ in probability}.
$$
This completes the proof of Theorem \ref{thm:nilpotent-wlln}.
\end{proof}

\begin{remark}[Virtually nilpotent reduction]\label{rem:virt-nilp}
Suppose $\Gamma$ is virtually nilpotent, so there exists a torsion-free nilpotent subgroup $K\le \Gamma$ of finite index. Fix a word metric $|\cdot|_\Gamma$ on $\Gamma$
from a finite generating set, and let $|\cdot|_K$ be a word metric on $K$.
Since $K$ has finite index, the inclusion $K\hookrightarrow \Gamma$ is a
quasi-isometry, so there exist constants $A\ge1$ and $B\ge0$ such that
$$
A^{-1}|h|_K - B \le |h|_\Gamma \le A|h|_K + B
\qquad\text{for all }h\in K.
$$

Assume in addition that the support of $\mu$ is contained in $K$. Then the
random walk $(X_n)$ with increments $S_i\sim\mu$ actually takes values in $K$.
We may therefore regard $(X_n)$ as a random walk on the nilpotent group $K$
with the same increment law $\mu$.

If $\mu$ satisfies the shell condition \eqref{eq:shell-condition} with respect to the
metric $|\cdot|_K$, then Theorem \ref{thm:nilpotent-wlln} applies to the $\mu$-walk on $K$ and yields
$|X_n|_K/n\to0$ in probability. By the quasi-isometry inequality above we have
$$
\frac{|X_n|_\Gamma}{n}\le \frac{A|X_n|_K+B}{n} \longrightarrow 0\text{ in probability as }n\to\infty,
$$
and similarly the converse implication (from $|X_n|_\Gamma/n\to0$ to
$|X_n|_K/n\to0$) holds. In particular, the conclusion of
Theorem \ref{thm:nilpotent-wlln} extends to virtually nilpotent groups in the case where
the random walk is supported in a finite-index nilpotent subgroup.
\end{remark}

Now we give a strengthened version of the observation at the beginning of this subsection, in a general group-theoretic setting.

\begin{theorem}[Finite speed $\Rightarrow$ finite first moment]\label{thm:finite-speed-finite-moment}
Let $\GG$ be a group endowed with a length function $|\cdot|$ satisfying
$$
|gh|\le |g|+|h|\quad\text{and}\quad |g|=|g^{-1}|
$$
(e.g.\ a word metric for a symmetric generating set). Let $(S_n)_{n\ge1}$ be i.i.d.\ $\GG$-valued increments with law $\mu$, and set $X_n=S_1\cdots S_n$.
If there exist $v<\infty$ and an event $E$ with $\mathbb P(E)>0$ such that
$$
\frac{|X_n|}{n}\longrightarrow v
\quad\text{on }E,
$$
then $\mathbb E|S_1|<\infty$.
\end{theorem}

\begin{proof}
We argue by contrapositive. Assume $\mathbb E|S_1|=\infty$.

Fix any $\lambda>0$ and set $Y:=|S_1|/\lambda$. Then $Y\ge0$ and $\mathbb E Y=\infty$.
For any non-negative random variable $Y$ one has
$$
\mathbb E\lfloor Y\rfloor = \sum_{n=1}^\infty \mathbb P(Y\ge n).
$$
Now, $\mathbb E Y=\infty$ implies $\mathbb E\lfloor Y\rfloor=\infty$ and hence
$$
\sum_{n=1}^\infty \mathbb P(Y\ge n)
= \sum_{n=1}^\infty \mathbb P(|S_1|\ge \lambda n)
= \infty.
$$
Define $A_n(\lambda):=\{|S_n|\ge \lambda n\}$ for $n\ge1$. Since the $(S_n)$ are i.i.d., the
events $(A_n(\lambda))_{n\ge1}$ are independent, and
$$
\sum_{n=1}^\infty \PP\big(A_n(\lambda)\big)
= \sum_{n=1}^\infty \PP\big(|S_1|\ge \lambda n\big)
= \infty.
$$
By the second Borel-Cantelli lemma,
$$
\PP\big(|S_n|\ge \lambda n\text{ i.o.}\big)=1
\qquad\text{for every }\lambda>0.
$$
In particular, for each integer $m\ge1$ we have
$$
\PP\big(|S_n|\ge m n\text{ i.o.}\big)=1.
$$
Taking the intersection over all $m\in\mathbb N$, we get
$$
\PP\Big(\forall m\in\mathbb N,\ \ |S_n|\ge mn\ \text{i.o.}\Big)=1,
$$
and on this event
$$
\limsup_{n\to\infty}\frac{|S_n|}{n} = \infty.
$$

On the other hand, for all $n$ we have, by subadditivity and symmetry of the length,
$$
|S_n| = |X_{n-1}^{-1}X_n|\le |X_{n-1}|+|X_n|.
$$
Hence on any event where $\frac{|X_n|}{n}\to v<\infty$, we get
$$
\limsup_{n\to\infty}\frac{|S_n|}{n}
\ \le\ \lim_{n\to\infty}\Big(\frac{|X_{n-1}|}{n}+\frac{|X_n|}{n}\Big)
= 2v < \infty.
$$

Now suppose there exists an event $E$ with $\PP(E)>0$ and $\frac{|X_n|}{n}\to v<\infty$ on $E$.
Then on $E$ we must have
$$
\limsup_{n\to\infty}\frac{|S_n|}{n} < \infty,
$$
whereas by the previous paragraph
$$
\PP\Big(\limsup_{n\to\infty}\frac{|S_n|}{n}=\infty\Big)=1.
$$
This is impossible unless $\PP(E)=0$, contradicting the assumption that $\PP(E)>0$.
Therefore our contrapositive assumption $\EE|S_1|=\infty$ is incompatible with the existence
of such an $E$, and we conclude $\EE|S_1|<\infty$.
\end{proof}

\begin{remark}
The argument uses only two structural inputs:
\begin{itemize}
\item subadditivity and symmetry of the length, $|gh|\le |g|+|h|$ and $|g|=|g^{-1}|$;
\item independence of the increments, to invoke the second Borel-Cantelli lemma.
\end{itemize}
Consequently, the result holds in any countable group endowed with such a symmetric length function;
finite generation plays no role. No property of the law $\mu$ beyond independence is needed:
in particular, the measure need not be symmetric, supported near the identity, or have any
finite moment a priori.

Moreover, we did not really use the existence of the \emph{limit} $\lim_{n\to\infty}|X_n|/n$ on $E$,
only that
$$
\limsup_{n\to\infty}\frac{|X_n|}{n} < \infty
$$
on an event $E$ of positive probability. In that situation, the same proof shows that
$\EE|S_1|<\infty$. Indeed, if $\limsup |X_n|/n < L <\infty$ on $E$, then
$$
\limsup_{n\to\infty}\frac{|S_n|}{n}
\ \le\ \limsup_{n\to\infty}\Big(\frac{|X_{n-1}|}{n}+\frac{|X_n|}{n}\Big)
\ \le\ 2L
$$
on $E$, contradicting the almost sure divergence of $\limsup |S_n|/n$ under $\EE|S_1|=\infty$.
\end{remark}

\begin{prop}[Finite first moment and symmetry $\Rightarrow$ zero speed]\label{prop:stepB-all}
Let $\GG$ be a finitely generated group of polynomial growth (equivalently, virtually nilpotent),
and let $\mu$ be a symmetric, non-degenerate probability measure on $\GG$ with finite first moment
$\EE|S_1|<\infty$. Let $X_n=S_1\cdots S_n$ be the associated random walk. Then
$$
\frac{|X_n|}{n}\longrightarrow 0\qquad\text{almost surely as }n\to\infty.
$$
\end{prop}

\begin{proof}
By integrability, Kingman's subadditive ergodic theorem applies to the process
$n\mapsto |X_n|$ and yields the almost sure limit
$$
\frac{|X_n|}{n}\longrightarrow v
$$
for some constant $v\ge0$. Thus it remains to show $v=0$.

For symmetric random walks on groups of polynomial growth, it is a classical result that
the drift $v$ must vanish. There are several equivalent formulations and proofs in the
literature. One route, in the nilpotent (and more generally polynomial-growth) setting,
is to use sharp heat kernel and Green-function estimates (Guivarc'h, Varopoulos, Hebisch-Saloff-Coste)
to show that linear escape is incompatible with the Gaussian-type bounds and volume growth.
Another route is via the study of the abelianisation and higher commutator coordinates,
using martingale methods and the central limit behaviour on nilpotent groups.

We refer, for example, to Guivarc'h \cite{Guivarch1973}, Varopoulos \cite{Varopoulos1985},
Hebisch-Saloff-Coste \cite{HS-C93}, and to Woess's monograph \cite{Woess2000} for detailed
treatments. In any of these approaches, one concludes that $v=0$, and hence $|X_n|/n\to0$
almost surely.
\end{proof}

\subsection{Dispersion, heavy tails, and abelian Poisson boundaries}\label{sec:dispersion-Margulis}

In this subsection we isolate a simple \emph{dispersion} property of random walks and explain how,
on nilpotent groups of class at most~$2$, it forces every bounded harmonic function to factor
through the abelianisation. 

\subsubsection{Dispersion along the centre and a rigidity theorem}

Let $\GG$ be a finitely generated group and let $\mu$ be a non-degenerate probability measure on $\GG$.
We write $\mu^{(n)}$ for the $n$-fold convolution power. For two probability measures $\nu_1,\nu_2$
on $\GG$, we use the total variation distance
$$
\|\nu_1-\nu_2\|_{\mathrm{TV}}
:=\frac12\sum_{g\in \GG}|\nu_1(g)-\nu_2(g)|.
$$

\begin{definition}[Dispersion along the centre]\label{def:dispersion_center}
We say that $(\GG,\mu)$ has \emph{dispersion property along the centre} if for every $z\in Z(\GG)$,
\begin{equation}\label{eq:dispersion_center}
\big\|\mu^{(n)}-\mu^{(n)}*\delta_{z}\big\|_{\mathrm{TV}}\to 0 \text{ as } n\to\infty.
\end{equation}
\end{definition}

\begin{lemma}[Central invariance from dispersion]\label{lem:center_invariance}
Let $\GG$ be a finitely generated group and $\mu$ a non-degenerate, symmetric probability measure on $\GG$
satisfying \eqref{eq:dispersion_center}. If $h:\GG\to\RR$ is bounded and $\mu$-harmonic, then
$$
h(gz)=h(g)\qquad\text{for all }g\in \GG, z\in Z(\GG).
$$
\end{lemma}

\begin{proof}
Fix $g\in \GG$ and $z\in Z(\GG)$, and set $M:=\sup_{x\in G}|h(x)|<\infty$.
Iterating harmonicity gives, for every $n\in\NN$,
$$
h(g)=\sum_{x\in \GG}\mu^{(n)}(x)h(gx).
$$
Similarly,
$$
h(gz)=\sum_{x\in \GG}\mu^{(n)}(x)h(gzx).
$$
Since $z\in Z(\GG)$, $gzx=gxz$ for all $x$. Put $y=xz$ (so $x=yz^{-1}$) to obtain
$$
h(gz)=\sum_{y\in \GG}\mu^{(n)}(yz^{-1})h(gy)
=\sum_{y\in \GG}(\mu^{(n)}*\delta_z)(y)h(gy),
$$
because $(\mu^{(n)}*\delta_z)(y)=\mu^{(n)}(yz^{-1})$.
Subtracting and using $|h|\le M$,
\begin{align*}
|h(g)-h(gz)|
&=\Big|\sum_{x\in \GG} h(gx)\Big(\mu^{(n)}(x)-(\mu^{(n)}*\delta_z)(x)\Big)\Big|\\
&\le \sum_{x\in \GG}|h(gx)|\big|\mu^{(n)}(x)-(\mu^{(n)}*\delta_z)(x)\big|\\
&\le M\sum_{x\in \GG}\big|\mu^{(n)}(x)-(\mu^{(n)}*\delta_z)(x)\big|
= 2M\big\|\mu^{(n)}-\mu^{(n)}*\delta_z\big\|_{\mathrm{TV}}.
\end{align*}
By \eqref{eq:dispersion_center}, the right-hand side tends to $0$ as $n\to\infty$, while the
left-hand side is independent of $n$. Hence, $|h(g)-h(gz)|=0$.
\end{proof}

\begin{prop}[Dispersion and bounded harmonic functions in nilpotency class two]\label{prop:Margulis_class2}
Let $\GG$ be a finitely generated nilpotent group of nilpotency class at most $2$, and let $\mu$
be a symmetric, non-degenerate probability measure on $\GG$ satisfying \eqref{eq:dispersion_center}.
Then every bounded $\mu$-harmonic function on $\GG$ is constant on cosets of the commutator subgroup $[\GG,\GG]$.

Equivalently, every bounded $\mu$-harmonic function factors through the abelianisation
$\pi:\GG\to \GG^{\mathrm{ab}}:=\GG/[\GG,\GG]$, and pullback by $\pi$ identifies the spaces of bounded
harmonic functions on $(\GG,\mu)$ and $(\GG^{\mathrm{ab}},\pi_*\mu)$, where $\pi_*\mu$ is defined on $\GG^\mathrm{ab}$ by $\pi_*\mu(g[\GG,\GG])=\sum_{h\in[\GG,\GG]}\mu(gh)$.
\end{prop}

\begin{proof}
Since $\GG$ has nilpotency class at most~$2$, we have $[\GG,\GG]\subset Z(\GG)$.
Let $h$ be bounded and $\mu$-harmonic. By Lemma \ref{lem:center_invariance},
$h(gz)=h(g)$ for all $z\in Z(\GG)$, hence in particular for all $z\in[\GG,\GG]$.
Thus $h$ is constant on each coset $g[\GG,\GG]$, so there exists a bounded function
$h^{\mathrm{ab}}:\GG^{\mathrm{ab}}\to\RR$ with $h=h^{\mathrm{ab}}\circ\pi$.

To see that $h^{\mathrm{ab}}$ is $\pi_*\mu$-harmonic, let $\bar x=\pi(x)$ and compute:
\begin{align*}
h^{\mathrm{ab}}( x[\GG,\GG])
=h(x)
=\sum_{g\in \GG}\mu(g)h(xg)
=\sum_{g[\GG,\GG]\in \GG^\mathrm{ab}}\pi_*\mu(g[\GG,\GG])h^{\mathrm{ab}}(x[\GG,\GG]g[\GG,\GG]).
\end{align*}
\end{proof}

\subsubsection{Verifying dispersion for heavy-tailed lifts on $H\times\ZZ$}

We now record a class of examples where dispersion along the centre can be checked directly.

\begin{definition}[Heavy-tailed product lift on $H\times\ZZ$]\label{def:heavy_lift}
Let $H$ be a finitely generated group of polynomial growth and let $\nu$ be a non-degenerate,
aperiodic probability measure on $H$ (e.g.\ $\nu(e_H)>0$).

Let $\mu_{\ZZ}$ be a non-degenerate probability measure on $\ZZ$ such that the greatest common divisor of $\{x-y : x,y \in \supp(\mu_\ZZ) \}$ is $1$, and which lies in the domain of attraction of a non-degenerate strictly $\alpha$-stable law for some $\alpha\in(0,2)$, meaning that if
$Z_1,Z_2,\dots$ are i.i.d.\ with law $\mu_{\ZZ}$ and $S_n:=Z_1+\cdots+Z_n$. Then there exist
sequences $a_n\to\infty$ and $b_n\in\RR$ and a non-degenerate strictly $\alpha$-stable random
variable $Y$ such that
$$
\frac{S_n-b_n}{a_n}\to Y \text{ in distribution}.
$$
Set
$$
G:=H\times\ZZ,\qquad \mu:=\nu\otimes\mu_{\ZZ}.
$$
\end{definition}

\begin{remark}[Product decomposition]
Let $(X_n)_{n\ge0}$ be the $\nu$-random walk on $H$ and $(S_n)_{n\ge0}$ the $\mu_{\ZZ}$-random walk on $\ZZ$,
driven by independent increments. Then $(X_n,S_n)$ is the $\mu$-random walk on $G=H\times\ZZ$.
\end{remark}

\begin{lemma}[Dispersion for stable-domain walks on $\ZZ$]\label{lem:dispersion_Z}
Let $\mu_{\ZZ}$ be as in Definition \ref{def:heavy_lift}. Then for every fixed $k\in\ZZ$,
$$
\big\|\mu_{\ZZ}^{(n)}-\mu_{\ZZ}^{(n)}*\delta_{k}\big\|_{\mathrm{TV}}\to 0 \text{ as } n\to\infty.
$$
\end{lemma}

\begin{proof}
Let $Z_1,Z_2,\dots$ be i.i.d.\ with law $\mu_{\ZZ}$ and $S_n:=Z_1+\cdots+Z_n$, and write
$p_n(x):=\PP(S_n=x)=\mu_{\ZZ}^{(n)}(x)$.
By Definition \ref{def:heavy_lift} there exist $a_n\to\infty$ and $b_n\in\RR$ such that
$(S_n-b_n)/a_n\to Y$ in distribution for a non-degenerate strictly $\alpha$-stable random variable $Y$.
Such stable laws admit a continuous density $f$ on $\RR$.

Moreover, the asumptions on $\mu_\ZZ$ yield the \emph{uniform stable local limit theorem} \cite[Section 50]{GnedenkoKolmogorov1968}:
\begin{equation}\label{eq:stableLLT_uniform_dispersion}
\sup_{x\in\ZZ}\left|a_n p_n(x)-f\left(\frac{x-b_n}{a_n}\right)\right|\to 0 \text{ as } n\to\infty.
\end{equation}
Fix $k\in\ZZ$. Then
$$
\big\|\mu_{\ZZ}^{(n)}-\mu_{\ZZ}^{(n)}*\delta_{k}\big\|_{\mathrm{TV}}
=\frac12\sum_{m\in\ZZ}|p_n(m)-p_n(m-k)|.
$$
Fix $M>0$ and split the sum as follows:
$$
\sum_{m\in\ZZ}|p_n(m)-p_n(m-k)|
=\underbrace{\sum_{|m-b_n|\le Ma_n}|p_n(m)-p_n(m-k)|}_{B_{n,M}}
 +\underbrace{\sum_{|m-b_n|>Ma_n}|p_n(m)-p_n(m-k)|}_{T_{n,M}}.
$$

\smallskip\noindent
Let $\varepsilon_n:=k/a_n$ and define 
$$
r_n(x):=p_n(x)-\frac1{a_n}f\left(\frac{x-b_n}{a_n}\right),
$$
so that \eqref{eq:stableLLT_uniform_dispersion} implies $\sup_{x\in\ZZ}|r_n(x)|=o(1/a_n)$.
For $|m-b_n|\le Ma_n$, the triangle inequality gives
$$
|p_n(m)-p_n(m-k)|
\le \frac1{a_n}\left|f\left(\frac{m-b_n}{a_n}\right)-f\left(\frac{m-b_n}{a_n}-\varepsilon_n\right)\right|
     +|r_n(m)|+|r_n(m-k)|.
$$
Summing over $|m-b_n|\le Ma_n$ yields
$$
B_{n,M}
\le \frac1{a_n}\sum_{|m-b_n|\le Ma_n}
\left|f\left(\frac{m-b_n}{a_n}\right)-f\left(\frac{m-b_n}{a_n}-\varepsilon_n\right)\right| + o(1).
$$
Set $g_n(y):=|f(y)-f(y-\varepsilon_n)|$. For $n$ large we have $|\varepsilon_n|\le 1$ and $g_n$ evaluates $f$ on the compact
interval $[-M-1,M+1]$. Define the modulus of continuity of $f$ on this compact set by
$$
\omega_M(\delta):=\sup\bigl\{|f(u)-f(v)|: u,v\in[-M-1,M+1], |u-v|\le\delta\bigr\},
$$
so $\omega_M(\delta)\to0$ as $\delta\downarrow0$. Then for $y,y'\in[-M,M]$ and $n$ large,
$$
|g_n(y)-g_n(y')|
\le |f(y)-f(y')|+|f(y-\varepsilon_n)-f(y'-\varepsilon_n)|
\le 2\omega_M(|y-y'|),
$$
so $(g_n)$ is uniformly equicontinuous on $[-M,M]$. Also
$\sup_{y\in[-M,M]}g_n(y)\le 2\sup_{[-M-1,M+1]}f=:C_M<\infty$.
With mesh size $h_n:=1/a_n\to0$ and $y_m:=(m-b_n)/a_n$, a standard Riemann-sum estimate gives
$$
\left|\frac1{a_n}\sum_{|m-b_n|\le Ma_n} g_n\left(\frac{m-b_n}{a_n}\right)
-\int_{-M}^{M} g_n(y)dy\right|
\le 4M\to 0,
$$
as $n\to\infty$. Therefore
$$
\limsup_{n\to\infty}B_{n,M}
\le \limsup_{n\to\infty}\int_{-M}^{M}|f(y)-f(y-\varepsilon_n)|dy
\le \limsup_{n\to\infty}\|f(\cdot)-f(\cdot-\varepsilon_n)\|_{L^1(\RR)}=0,
$$
using $L^1$-continuity of translations and $\varepsilon_n\to0$.

\smallskip\noindent
By the triangle inequality,
$$
T_{n,M}\le \PP(|S_n-b_n|>Ma_n)+\PP(|S_n-b_n-k|>Ma_n).
$$
For $n$ large, $|k|\le (M/2)a_n$, hence
$$
\PP(|S_n-b_n-k|>Ma_n)\le \PP(|S_n-b_n|>(M/2)a_n),
$$
so
$$
\limsup_{n\to\infty}T_{n,M}
\le 2\limsup_{n\to\infty}\PP\left(\left|\frac{S_n-b_n}{a_n}\right|\ge\frac{M}{2}\right).
$$
By $(S_n-b_n)/a_n\to Y$ in distribution and Portmanteau's theorem,
$$
\limsup_{n\to\infty}\PP\left(\left|\frac{S_n-b_n}{a_n}\right|\ge\frac{M}{2}\right)
\le \PP(|Y|\ge M/2).
$$
Since $\{|Y|\ge t\}\downarrow\varnothing$ as $t\to\infty$, we have $\PP(|Y|\ge t)\to0$. Fix $\varepsilon>0$ and choose $M$ large enough so that $\limsup_{n\to\infty}T_{n,M} < \varepsilon$. This implies there exists $N_1(\varepsilon,M) \in \NN$ such that 
$$
T_{n,M} < \frac{\varepsilon}{2} \text{ for all } n \geq N_1.
$$
Since $\lim_{n\to\infty} B_{n,M}=0$, choose $N_2=N_2(\varepsilon,M)\in\NN$ such that $B_{n,M} < \frac{\varepsilon}{2}$ for all $n\geq N_2$.
Combining, we get 
$$
\sum_{m\in\ZZ}|p_n(m)-p_n(m-k)|\le T_{n,M} + B_{n,M}<\varepsilon \text{ for all } n\geq \max\{N_1, N_2\},
$$ 
i.e. the desired dispersion on $\ZZ$.
\end{proof}

\begin{lemma}[Dispersion for the product lift]\label{lem:dispersion_lift}
Let $H$ and $\nu$ be as in Definition \ref{def:heavy_lift}, and assume that for every $z_H\in Z(H)$,
$$
\big\|\nu^{(n)}-\nu^{(n)}*\delta_{z_H}\big\|_{\mathrm{TV}}\to 0 \text{ as }n\to\infty.
$$
Let $\mu_{\ZZ}$ be as in Definition \ref{def:heavy_lift} and set $G:=H\times\ZZ$, $\mu:=\nu\otimes\mu_{\ZZ}$.
Then, for every $z=(z_H,z_{\ZZ})\in Z(G)=Z(H)\times\ZZ$,
$$
\big\|\mu^{(n)}-\mu^{(n)}*\delta_{z}\big\|_{\mathrm{TV}}\to 0 \text{ as }n\to\infty.
$$
In particular, $(G,\mu)$ satisfies the dispersion property along its centre.
\end{lemma}

\begin{proof}
Write $p_n:=\nu^{(n)}$ on $H$ and $q_n:=\mu_{\ZZ}^{(n)}$ on $\ZZ$. Then $\mu^{(n)}=p_n\otimes q_n$ on $G$.
For $z=(z_H,z_{\ZZ})\in Z(G)$ we have
$$
\mu^{(n)}*\delta_z = (p_n*\delta_{z_H})\otimes(q_n*\delta_{z_{\ZZ}}).
$$
For probability measures $\alpha_1,\alpha_2$ on $H$ and $\beta_1,\beta_2$ on $\ZZ$,
$$
\alpha_1\otimes\beta_1-\alpha_2\otimes\beta_2
=(\alpha_1-\alpha_2)\otimes\beta_1+\alpha_2\otimes(\beta_1-\beta_2),
$$
and hence 
$$
\|\alpha_1\otimes\beta_1-\alpha_2\otimes\beta_2\|_{\mathrm{TV}}
\le \|\alpha_1-\alpha_2\|_{\mathrm{TV}}+\|\beta_1-\beta_2\|_{\mathrm{TV}}.
$$
Applying this with $\alpha_1=p_n$, $\alpha_2=p_n*\delta_{z_H}$ and $\beta_1=q_n$, $\beta_2=q_n*\delta_{z_{\ZZ}}$ yields
$$
\big\|\mu^{(n)}-\mu^{(n)}*\delta_z\big\|_{\mathrm{TV}}
\le \big\|p_n-p_n*\delta_{z_H}\big\|_{\mathrm{TV}}+\big\|q_n-q_n*\delta_{z_{\ZZ}}\big\|_{\mathrm{TV}}.
$$
The first term tends to $0$ by hypothesis, and the second tends to $0$ by Lemma \ref{lem:dispersion_Z}.
\end{proof}

In particular, when $H$ is nilpotent of class at most~$2$, the product $G=H\times\ZZ$ is again of class at most~$2$.
Combining Lemma \ref{lem:dispersion_lift} with Proposition \ref{prop:Margulis_class2} shows that, for such heavy-tailed
product lifts, every bounded $\mu$-harmonic function on $G$ factors through the abelianisation $G^{\mathrm{ab}}$.

\subsection{Acknowledgements}  The research of the first author was partially supported by SEED Grant RD/0519-IRCCSH0-024. During the initial stages of the preparation of the manuscript the first author was a visitor at MPIM Bonn. The second author would like to thank the PMRF for partially supporting his work. 
The third author was partially supported by IIT Bombay IRCC fellowship, TIFR Mumbai post-doctoral fellowship and NBHM postdoctoral fellowship (Sr. No. 0204/17/2025/R\&D-II/12398) during this work. The authors are deeply grateful to Gideon Amir, S\'{e}bastien Gou\"{e}zel and Debanjan Nandi for very insightful correspondence. All three authors would like to thank IIT Bombay for providing ideal working conditions.

\bibliographystyle{alpha}
\bibliography{references1}

\end{document}